\documentclass[a4paper,11pt]{article}

\usepackage{amsmath}
\usepackage{amsfonts}
\usepackage{amssymb}
\usepackage{latexsym}
\usepackage{epsfig}
\usepackage{graphicx}
\usepackage{oldgerm}
\usepackage{theorem}

\def\la{\lambda}
\def\La{\Lambda}
\def\C{\mathbb{C}}
\def\R{\mathbb{R}}
\def\N{\mathbb{N}}
\def\Z{\mathbb{Z}}
\def\H{\mathbb{H}}
\def\S{\mathbb{S}}
\def\B{\mathbb{B}}
\def\T{\mathbb{T}}

\def\ker{{\rm ker} \,}
\def\Ai{{\rm Ai}}
\def\gO{{\rm O}}
\def\U{{\rm U}}
\def\SO{{\rm SO}}
\def\Sp{{\rm Sp}}

\def\rR{{\rm R}}
\def\rI{{\rm I}}
\def\rN{{\rm N}}
\def\bE{{\bf E}}
\def\bP{{\bf P}}
\def\bp{{\bf p}}
\def\x{\mib{x}}
\def\y{\mib{y}}
\def\z{\mib{z}}
\def\RN{{\it R}_N}
\def\AN{{\it A}_{N-1}}
\def\BN{{\it B}_N}
\def\BNv{{\it B}^{\vee}_N}
\def\CN{{\it C}_N}
\def\CNv{{\it C}^{\vee}_N}
\def\BCN{{\it BC}_N}
\def\DN{{\it D}_N}

\def\cC{{\cal C}}
\def\cS{{\cal S}}
\def\cP{{\cal P}}
\def\cH{{\cal H}}
\def\cY{{\cal Y}}
\def\cO{{\cal O}}
\def\cK{{\cal K}}
\def\cN{{\cal N}}
\def\cW{{\cal W}}
\def\cF{{\cal F}}
 
\def\cA{{\cal A}}
\def\cB{{\cal B}}
\def\cS{{\cal S}}
\def\cU{{\cal U}}

\def\law{\stackrel{\rm (law)}{=}}
\def\weak{\stackrel{N \to \infty}{\, \Longrightarrow \,}}
\def\weakT{\stackrel{\Im \tau \to \infty}{\, \Longrightarrow \,}}
\def\weaka{\stackrel{\alpha \to \infty}{\, \Longrightarrow \,}}
\def\weakp{\stackrel{p \to \infty}{\, \Longrightarrow \,}}
\def\Det{\mathop{\mathrm{Det}}}
\def\tr{{\rm Tr} \,}
\def\locally{\mathcal{I}_{1,\mathrm{loc}}}

\def\etabar{\overline{\eta}}
\def\zetabar{\overline{\zeta}}

\theorembodyfont{\itshape}

\newtheorem{thm}{Theorem}[section]
\newtheorem{lem}[thm]{Lemma}
\newtheorem{cor}[thm]{Corollary}
\newtheorem{prop}[thm]{Proposition}

\newtheorem{df}[thm]{Definition}

\newcommand{\mib}[1]{\mbox{\boldmath $#1$}}
\newcommand{\SSC}[1]{\section{#1}\setcounter{equation}{0}}
\newcommand{\qed}{\hbox{\rule[-2pt]{3pt}{6pt}}}

\begin{document}

\title{{\bf 
Partial Isometries, Duality, \\
and Determinantal Point Processes}
\\
\vskip 0.5cm
{\normalsize {\it Dedicated to Professor Hirofumi Osada 
on the occasion of his 60th birthday}}
}

\author{
Makoto Katori \\
Department of Physics,
Faculty of Science and Engineering, \\
Chuo University,
Kasuga, Bunkyo-ku, Tokyo 112-8551, Japan \\
e-mail: katori@phys.chuo-u.ac.jp \\
\\
Tomoyuki Shirai \\
Institute of Mathematics for Industry, 
Kyushu University, \\
744 Motooka, Nishi-ku,
Fukuoka 819-0395, Japan \\
e-mail: shirai@imi.kyushu-u.ac.jp
}

\date{7 September 2021}
\pagestyle{plain}
\maketitle

\begin{abstract}

A determinantal point process (DPP) is 
an ensemble of random nonnegative-integer-valued 
Radon measures $\Xi$ on 
a space $S$ with measure $\lambda$,
whose correlation functions are all
given by determinants specified by an
integral kernel $K$ called the correlation kernel.
We consider a pair of Hilbert spaces, 
$H_{\ell}, \ell=1,2$, which are assumed to
be realized as $L^2$-spaces,
$L^2(S_{\ell}, \lambda_{\ell})$, $\ell=1,2$, 
and introduce 
a bounded linear operator ${\cal W} : H_1 \to H_2$
and its adjoint ${\cal W}^{\ast} : H_2 \to H_1$.
We show that if ${\cal W}$ is a partial isometry 
of locally Hilbert--Schmidt class, then we have a unique DPP 
$(\Xi_1, K_1, \lambda_1)$ associated with $\cW^* \cW$. 
In addition, if $\cW^*$ is also 
of locally Hilbert--Schmidt class, then we have a unique pair
of DPPs, $(\Xi_{\ell}, K_{\ell}, \lambda_{\ell})$, $\ell=1,2$.
We also give a practical framework 
which makes ${\cal W}$ and ${\cal W}^{\ast}$ satisfy the above conditions.
Our framework to construct pairs of DPPs implies
useful duality relations between DPPs making pairs.
For a correlation kernel of a given DPP
our formula can provide plural different expressions, 
which reveal different aspects of the DPP.
In order to demonstrate these advantages of 
our framework as well as to show that
the class of DPPs obtained by this method is large enough to
study universal structures in a variety of DPPs, 
we report plenty of examples of DPPs 
in one-, two-, and higher-dimensional spaces $S$,
where several types of weak convergence from finite DPPs
to infinite DPPs are given.
One-parameter ($d \in \mathbb{N}$) series of infinite DPPs
on $S=\mathbb{R}^d$ and $\mathbb{C}^d$ are discussed,
which we call the Euclidean and the Heisenberg
families of DPPs, respectively, following 
the terminologies of Zelditch.

\vskip 0.3cm


\vskip 0.3cm


\end{abstract}

\SSC
{Introduction} \label{sec:Introduction}
Let $S$ be a base space, 
which is a locally compact Hausdorff space
with countable base, 
and $\lambda$ be a Radon measure on $S$.
The configuration space over $S$ is given by
the set of nonnegative-integer-valued Radon measures; 
\[
{\rm Conf}(S)
=\left\{ \xi = \sum_j \delta_{x_j} : \mbox{$x_j \in S$,
$\xi(\Lambda) < \infty$ for all bounded set $\Lambda \subset S$} \right\}.
\]
${\rm Conf}(S)$ is equipped with the topological Borel $\sigma$-fields
with respect to the vague topology; 
we say $\xi_n, n \in \N :=\{1, 2, \dots\}$ converges
to $\xi$ in the vague topology, if
$\int_{S} f(x) \xi_n(dx) \to
\int_{S} f(x) \xi(dx)$,
$\forall f \in \cC_{\rm c}(S)$, 
where $\cC_{\rm c}(S)$ is the set of
all continuous real-valued functions 
with compact support.
A {\it point process} on $S$ is
a ${\rm Conf}(S)$-valued random variable
$\Xi=\Xi(\cdot, \omega)$ on a probability space
$(\Omega, \cF, \bP)$.
If $\Xi(\{x\}) \in \{0, 1\}$ for any point $x \in S$,
then the point process is said to be {\it simple}.

Assume that $\Lambda_j, j =1, \dots, m$, $m \in \N$ are
disjoint bounded sets in $S$ and
$k_j \in \N_0 :=\{0,1, \dots\}, j=1, \dots, m$ satisfy
$\sum_{j=1}^m k_j = n \in \N_0$.
A symmetric measure $\lambda^n$ on $S^n$
is called the $n$-th {\it correlation measure},
if it satisfies
\[
\bE \left[
\prod_{j=1}^m \frac{\Xi(\Lambda_j)!}
{(\Xi(\Lambda_j)-k_j)!} \right]
=\lambda^n(\Lambda_1^{k_1} \times \cdots
\times \Lambda_m^{k_m}),
\]
where if $\Xi(\Lambda_j)-k_j <0$,
we interpret $\Xi(\Lambda_j)!/(\Xi(\Lambda_j)-k_j)!=0$.
If $\lambda^n$ is absolutely continuous 
with respect to the $n$-product measure $\lambda^{\otimes n}$,
the Radon--Nikodym derivative
$\rho^n(x_1, \dots, x_n)$ is called
the {\it $n$-point correlation function}
with respect to the background measure $\lambda$;
\[
\lambda^n(dx_1 \cdots dx_n)
=\rho^n(x_1, \dots, x_n) 
\lambda^{\otimes n}(dx_1 \cdots dx_n).
\]
Determinantal point process (DPP) 
is defined as follows 
\cite{Macchi75,ST00,Sos00,ST03a,ST03b,HKPV06,HKPV09}. 

\begin{df}
A simple point process $\Xi$ on $(S, \lambda)$ is said to be
a determinantal point process 
{\rm ({\it DPP})} with 
correlation kernel $K : S \times S \to \C$ if it has correlation functions 
$\{\rho^n \}_{n \in \N}$, 
and they are given by
\begin{equation}
\rho^n(x_1, \dots, x_n) = \det_{1 \leq j, k \leq n}
[ K(x_j, x_k) ]
\quad \mbox{for every $n \in \N$, 
and $x_1, \dots, x_n \in S$}.
\label{eqn:DPP}
\end{equation}
The triplet $(\Xi, K, \lambda(d x))$ 
denotes the DPP;
$\Xi \in {\rm Conf}(S)$,  
specified by the correlation kernel $K$
with respect to the measure $\lambda(d x)$.
\end{df}

If the integral projection operator $\cK$ on $L^2(S,\lambda)$ 
with a kernel $K$ is of rank $N \in \N$,
then the number of points is $N$ a.s.
If $N < \infty$ (resp. $N=\infty$),
we call the system a {\it finite DPP} (resp. an {\it infinite DPP}).
The density of points with respect to the background measure $\lambda(dx)$
is given by
\[
\rho(x) :=\rho^1(x) = K(x,x).
\]
The DPP is negatively correlated as shown by
\begin{align}
\rho^2(x, x') &= \det \left[
\begin{array}{ll}
K(x, x) & K(x, x') \cr
K(x', x) & K(x', x') 
\end{array}
\right]
\nonumber\\
&= K(x, x) K(x', x') - |K(x,x')|^2 \leq \rho(x) \rho(x'),
\quad x, x' \in S,
\label{eqn:rho2}
\end{align}
provided that $K$ is Hermitian.

Let $H$ be a separable Hilbert space. 
For linear operators $\cA, \cB$ on $H$, we say
that $\cA$ is positive definite and write 
$\cA \geq O$ if $\langle \cA f, f \rangle_H \geq 0$
for any $f \in H$, and write
$\cA \geq \cB$ if $\cA-\cB \geq O$. 
The operator $\cA^* \cA$ is
positive definite and it admits a unique 
positive-definite square-root
$\sqrt{\cA^* \cA}$ which is denoted by $|\cA|$. 
Let $\{\phi_n\}_{n \ge 1}$ be an orthonormal basis of $H$. 
For $\cA \ge O$, we define the trace of $\cA$ by 
\[
 \tr \cA := \sum_{n=1}^{\infty} \langle \cA \phi_n, \phi_n
 \rangle_H, 
\]
which does not depend on the choice of an orthonormal basis. 
A bounded linear operator
$\cA$ is said to be \textit{of trace class} or a
\textit{trace class operator} if the trace norm $\|\cA\|_1 := \tr |\cA|$ is
finite. The trace $\tr \cA$ is defined whenever $\|\cA\|_1 < \infty$. 

Now, we consider the case $H = L^2(S, \la)$.
For a compact set $\Lambda \subset S$,
the projection from $L^2(S, \lambda)$ to 
the space of all functions vanishing outside $\Lambda$ $\lambda$-a.e.
is denoted by $\cP_{\Lambda}$.
$\cP_{\Lambda}$ is the operation of multiplication of
the indicator function ${\bf 1}_{\Lambda}$ of the set $\Lambda$;
${\bf 1}_{\Lambda}(x)=1$ if $x \in \Lambda$, and
${\bf 1}_{\Lambda}(x)=0$ otherwise.
We say that a bounded linear operator 
$\cA$ on $L^2(S, \lambda)$ is 
{\it of locally trace class} or a {\it locally trace class operator},
if the restriction of $\cA$ to each compact subset $\Lambda$
is of trace class; that is, 
\[
\mbox{
$\| \cA_{\Lambda} \|_1 < \infty$ \,
with \, $\cA_{\Lambda} :=\cP_{\Lambda} \cA \cP_{\Lambda}$ \,
for any compact set $\Lambda \subset S.$
}
\]
The totality of locally trace
class operators on $L^2(S, \lambda)$ 
is denoted by $\locally(S, \lambda)$. 
It is known that \cite{ST00,Sos00,ST03a,ST03b}, 
if $\cK \in \locally(S, \lambda)$ and
$O \leq \cK \leq I$, 
where $I$ is the identity operator, then
we have a unique DPP on $S$ with 
the determinantal correlation functions (\ref{eqn:DPP})
with respect to $\lambda$ and
the correlation kernel $K$ is given by the Hermitian
integral kernel for $\cK$ 
(see Section \ref{sec:DPP} below). 

In the present paper, we consider the case in which 
\[
\cK f = f \quad
\mbox{for all $f \in (\ker \cK)^{\perp} \subset L^2(S, \lambda)$}.
\]
Here $(\ker \cK)^{\perp}$ denotes the
{\it orthogonal complement of the kernel space} of $\cK$.
That is, $\cK$ is an {\it orthogonal projection}.
By definition, it is obvious that the condition
$O \leq \cK \leq I$ is satisfied.
The purpose of the present paper is to
propose a useful method to provide
orthogonal projections $\cK$ and
DPPs whose correlation kernels are
given by the Hermitian integral kernels 
of $\cK$, $K(x, x'), x, x' \in S$.

We consider a pair of Hilbert spaces, 
$H_{\ell}, \ell=1,2$, which are assumed to
be realized as $L^2$-spaces,
$L^2(S_{\ell}, \lambda_{\ell})$, $\ell=1,2$.
We introduce 
a bounded linear operator $\cW$ and its adjoint $\cW^{\ast}$, 
\begin{equation}
\cW : H_1 \to H_2,
\quad
\cW^{\ast} : H_2 \to H_1.
\label{eqn:W_W*}
\end{equation}
Then, we have the following basic existence theorem of DPP
via a partial isometry $\cW$ of locally
Hilbert--Schmidt class. 
\begin{thm}\label{thm:exist_DPP_HSversion}
Assume that $\cW : L^2(S_1,\lambda_1) \to
 L^2(S_2,\lambda_2)$ is a partial isometry of locally
 Hilbert--Schmidt class. Then, there exists a
unique DPP $(\Xi_1, K_{S_1}, \la_1)$ on $S_1$ 
with 
\begin{equation}
 K_{S_1}(x,x') = \int_{S_2} \overline{W(y,x)} W(y,x')
  \la_2(dy),  
\label{eq:KS1_kernel} 
\end{equation}
where $\cW$ admits a measurable kernel $W : S_2 \times S_1
 \to \C$ such that $\Psi_1 \in L^2_{\mathrm{loc}}(S_1,
 \la_1)$ with $\Psi_1(x) := \|W(\cdot, x)\|_{L^2(S_2,
 \la_2)} \ (x \in S_1)$. 
\end{thm}
The definitions of partial isometries, 
locally Hilbert--Schmidt operators,
and $L^2_{\mathrm{loc}}(S, \la)$ will be
given in Section~\ref{sec:isometry}. 
There we will show
basic properties of them
and Theorem~\ref{thm:exist_DPP_HSversion} will be concluded from
the well-known existence theorem
of DPP \cite{ST00,Sos00,ST03a,ST03b} (Theorem~\ref{thm:exist_DPP}).

We assume that 
\begin{description}
\item{(i)} \quad ${\cal W}$ is a {\it partial isometry}, 
\item{(ii)} \quad both $\cW$ and $\cW^{\ast}$ are 
{\it locally Hilbert--Schmidt operators}. 
\end{description}
Under the assumption (i), the adjoint ${\cal W}^*$ is also a
partial isometry. 
If two conditions (i) and (ii) are satisfied, 
by Theorem~\ref{thm:exist_DPP_HSversion}, 
we have a unique pair
of DPPs, $(\Xi_{\ell}, K_{\ell}, \lambda_{\ell})$, $\ell=1,2$,
where the correlation kernel $K_1$
(resp. $K_2$) is given by the integral kernel of
$\cW^{\ast} \cW$ (resp. $\cW \cW^{\ast}$)
(Theorem \ref{thm:main1}) as in \eqref{eq:KS1_kernel}. 
We give a practical framework which makes 
${\cal W}$ and ${\cal W}^{\ast}$ satisfy the above two assumptions
(Corollaries \ref{thm:main2} and \ref{thm:main3}).

One of the advantages of our framework is that
the obtained pairs of DPPs satisfy useful duality relations,
which will be reported in 
Sections \ref{sec:duality}, \ref{sec:duality_application1}, 
\ref{sec:duality_application2}.
As mentioned above, 
one of a pair of DPPs discussed here is associated with 
a Hilbert space $H_1$ having
an orthogonal projection $\cK_1$, 
and $\cK_1$ is given in the form $\cK_1=\cW^{\ast} \cW$.
This equality can be regarded as a {\it decomposition formula}
of $\cK_1$ by a product of an operator 
$\cW$ and its dual $\cW^{\ast}$ acting as (\ref{eqn:W_W*})
provided that another Hilbert space $H_2$ is chosen.
We note the fact that for a given DPP associated with $H_1$
and $\cK_1$, 
choice of $H_2$ is not unique.
As demonstrated in Sections
\ref{sec:three_Ginibre}--\ref{sec:Fock} 
using the Ginibre DPPs on $\C$,
such multivalency in our framework 
can give plural different expressions
for one correlation kernel $K_1$ and they 
will help us to study different aspects of the DPP 
which we consider. 

In order to demonstrate the class of DPPs 
obtained by our framework is large enough to
study a variety of DPPs and universal structures behind them,
we show plenty of examples of DPPs in one- and two-dimensional spaces.
In particular, we use the symbols of classical and affine
roots systems (e.g., $\AN, \BN, \CN, \DN, N \in \N$) to classify
finite DPPs.
Several types of weak convergence theorems of finite DPPs
to infinite DPPs are given.
We will show that in the one-dimensional space,
there are three universal DPPs with an infinite number of points
specified by the correlation kernels,
\begin{align*}
K_{\rm sinc}(x, x') &:=
\frac{\sin(x-x')}{\pi (x-x')}
=\frac{1}{2 \pi} \int_{-1}^1 e^{i \gamma(x-x')} d \gamma,
\quad x, x' \in \R,
\nonumber\\
K_{\rm Bessel}^{(1/2)}(x, x')
&:= \frac{\sin(x-x')}{\pi(x-x')} - \frac{\sin(x+x')}{\pi (x+x')}
= \frac{1}{\pi} \int_{-1}^1 \sin(\gamma x) \sin(\gamma x') d \gamma,
\quad x, x' \in [0, \infty),
\nonumber\\
K_{\rm Bessel}^{(-1/2)}(x, x')
&:= \frac{\sin(x-x')}{\pi(x-x')} + \frac{\sin(x+x')}{\pi (x+x')}
= \frac{1}{\pi} \int_{-1}^1 \cos(\gamma x) \cos(\gamma x') d \gamma,
\quad x, x' \in [0, \infty),
\end{align*}
where $i :=\sqrt{-1}$.
$K_{\rm sinc}$ is usually called the sine kernel
in random matrix theory \cite{Meh04}, but it shall be called 
the {\it sinc kernel}.
$K_{\rm Bessel}^{(1/2)}$ and $K_{\rm Bessel}^{(-1/2)}$
are special cases of the {\it Bessel kernels} 
$K_{\rm Bessel}^{(\nu)}$, $\nu \in (-1, \infty)$ with
indices $\nu=1/2$ and $-1/2$, respectively \cite{For10}. 
Note that 
$K_{\rm sinc}(x, x')=\{K_{\rm Bessel}^{(1/2)}(x, x')+K_{\rm Bessel}^{(-1/2)}(x, x')\}/2$,
$x, x' \in [0, \infty)$.
Corresponding to the threefold of DPPs with the correlation kernels,
$K_{\rm sinc}$, $K_{\rm Bessel}^{(1/2)}$, 
$K_{\rm Bessel}^{(-1/2)}$,
we also show the three universal DPPs on $\C$,
whose correlation kernels are given by
\begin{align*}
K_{\rm Ginibre}^A(x, x')
&:= e^{x \overline{x'}}
=\sum_{n=0}^{\infty} \frac{(x \overline{x'})^n}{n!},
\nonumber\\
K_{\rm Ginibre}^C(x, x')
&:= \sinh(x \overline{x'})
=\sum_{n=0}^{\infty} \frac{(x \overline{x'})^{2n+1}}{(2n+1)!},
\nonumber\\
K_{\rm Ginibre}^D(x, x')
&:= \cosh(x \overline{x'})
=\sum_{n=0}^{\infty} \frac{(x \overline{x'})^{2n}}{(2n)!},
\quad x, x' \in \C,
\end{align*}
where $\overline{x'}$ denotes the complex conjugate
of $x'$.
$K_{\rm Ginibre}^{A}$ is known as the correlation kernel
of the {\it Ginibre ensemble} in random matrix theory
\cite{Gin65,For10}, and
$K_{\rm Ginibre}^{C}$ and $K_{\rm Ginibre}^{D}$
were studied in \cite{Kat19b}.
Note that
$K_{\rm Ginibre}^{A}(x,x')=K_{\rm Ginibre}^{C}(x,x')+K_{\rm Ginibre}^{D}(x,x')$,
$x, x' \in \C$.

Our method to generate DPPs is also valid
in higher dimensional spaces.
We will state that the DPP with the sinc kernel
$K_{\rm sinc}$ is the lowest-dimensional ($d=1$) 
example of the one-parameter ($d \in \N$)
family of DPPs on $\R^d$, whose correlation kernels
are given by
\begin{align*}
K_{\rm Euclid}^{(d)}(x, x')
&:= \frac{1}{(2 \pi)^{d/2}}
\frac{ J_{d/2}(\|x-x'\|_{\R^d}) }{ \|x-x'\|_{\R^d}^{d/2} }
\nonumber\\
&=\frac{1}{(2\pi)^d} \int_{\B^d}
e^{i \gamma \cdot (x-x')} d \gamma,
\quad x, x' \in \R^d,
\nonumber
\end{align*}
with respect to the Lebesgue measures of $\R^d$,
$\lambda(dx)=dx$, 
where $J_{\nu}$ is the Bessel function of the first kind,
$\|x-x'\|_{\R^d}$ is the Euclidean distance between
$x$ and $x'$ in $\R^d$,
and $\B^d$ is the unit ball in $\R^d$ centered at the origin.
We also claim that the Ginibre ensemble is the lowest-dimensional example
$(d=1)$ of another one-parameter $(d \in \N)$
family of DPPs on $\C^d$, whose correlation kernel is
given by
\[
K_{\rm Heisenberg}^{(d)}(x, x') :=
e^{x \cdot \overline{x'}}, \quad x, x' \in \C^d,
\]
where the background measure $\lambda$ is
assumed to be the $d$-dimensional 
complex normal distribution.
We call these two families of DPPs 
the {\it Euclidean family} of DPPs and the
{\it Heisenberg family} of DPPs, respectively,
following the terminologies by Zelditch \cite{Zel00}. 
See also \cite{BSZ00,SZ02,Zel09,CH15}.

The paper is organized as follows.
In Section \ref{sec:Main} we give main theorems
which enable us to generate DPPs.
Sections \ref{sec:examples_1d} and \ref{sec:examples_2d}
are devoted to a variety of examples
of DPPs obtained by our framework
for the one-dimensional and
the two-dimensional spaces, respectively.
Examples in spaces with arbitrary dimensions $d \in \N$
are given in Section \ref{sec:example_d_dim}.
We list out open problems in Section \ref{sec:remarks}.
Appendices \ref{sec:Weyl_denominators}
and \ref{sec:Macdonald_denominators}
are used to explain useful multivariate functions
and determinantal formulas
associated with the classical
and the affine root systems, respectively.
The definitions and basic properties of
the Jacobi theta functions are summarized in
Appendix \ref{sec:Jacobi}. 

\SSC
{Main Theorems} \label{sec:Main}

\subsection{Existence theorem of DPPs}
\label{sec:DPP}

We recall the existence theorem for DPPs. 
Let $(S, \lambda)$ be a $\sigma$-finite measure space. 
We assume that $\cK \in \locally(S, \lambda)$.
If, in addition, $\cK \ge O$, then it  
admits a Hermitian integral kernel $K(x, x')$ such that (cf. \cite{GY05})
\begin{description}
\item{(i)} \quad
$\displaystyle{\det_{1 \leq j, k \leq n}[K(x_j, x_k)] \ge 0}$ for 
$\lambda^{\otimes n}$-a.e.\,$(x_1,\dots, x_n)$ 
with every $n \in \N$, 

\item{(ii)} \quad
$K_{x'} := K(\cdot, x') \in L^2(S, \lambda)$ for
$\lambda$-a.e.\,$x'$, 

\item{(iii)} \quad
$\tr \cK_{\Lambda} = \int_{\Lambda} K(x, x) \lambda(d x)$, 
$\Lambda \subset S$ and 
\[
 \tr (\cP_{\Lambda} \cK^n \cP_{\Lambda}) = \int_{\Lambda} \langle K_{x'}, \cK^{n-2} K_{x'}
 \rangle_{L^2(S, \lambda)} \lambda(d x'), 
\quad \forall n \in \{2, 3, \dots\}. 
\]
\end{description}
This is based on the fact that every positive definite
trace class operator has the form $\cB^*\cB$ of a
Hilbert--Schmidt operator $\cB$
together with a similar idea of 
the proof of Proposition~\ref{prop:kernel-thm}
mentioned below. 

\begin{thm}[\cite{ST00,Sos00,ST03a,ST03b}]
\label{thm:exist_DPP}
Assume that $\cK \in \locally(S,\lambda)$ and $O \le \cK \le I$. 
Then there exists a
unique DPP $(\Xi, K, \la)$ on $S$.  
\end{thm}

If $\cK \in \locally(S, \lambda)$ is 
a projection onto a closed subspace 
$H \subset L^2(S, \lambda)$, 
one has {\it the DPP associated with $K$ and
$\lambda$}, or one may say 
{\it the DPP associated with the subspace $H$}. 
This situation often appears in the setting of
{\it reproducing kernel Hilbert space} \cite{Aro50}. 
Let $\cF = \cF(S)$ be 
a Hilbert space of complex functions on $S$ with inner
product $\langle \cdot, \cdot \rangle_{\cF}$.  
A function $K(x, x')$ on $S \times S$ is said
to be a {\it reproducing kernel} of $\cF$ if 
\begin{enumerate}
 \item For every $x^{\prime} \in S$, the function
       $K(\cdot, x^{\prime})$ belongs to $\cF$. 
\item The function $K(x, x^{\prime})$ has reproducing kernel
      property; that is, for any $f \in \cF$, 
\[
 f(x') = \langle f(\cdot), K(\cdot, x^{\prime}) \rangle_{\cF}. 
\] 
\end{enumerate}
A reproducing kernel of $\cF$ is unique if exists, 
and it exists if and only if 
the point evaluation map $\cF \ni f \mapsto f(x) \in \C$ is
bounded for every $x \in S$. 
The Moore--Aronszajn theorem states that if a kernel 
$K(\cdot, \cdot)$ on $S \times S$ is positive definite
in the sense that for any $n \geq 1$,
$x_1, \dots, x_n \in S$, 
the matrix $(K(x_j, x_k))_{j, k \in \{1, \dots, n\}}$ is positive definite, then  
there exists a unique Hilbert space $H_K$ of functions with 
inner product in which $K(x, x')$ is a reproducing
kernel \cite{Aro50}. If $H_K$ is realized in $L^2(S, \lambda)$ for some
measure $\lambda$, the kernel $K(x, x')$ defines a projection onto $H_K$. 

\subsection{Partial isometries, locally Hilbert--Schmidt operators, and DPPs}
\label{sec:isometry}

First we recall the notion of partial isometries between
Hilbert spaces \cite{Hal74,HL63}. 
Let $H_{\ell}, \ell=1,2$ be separable Hilbert spaces
with inner products $\langle \cdot, \cdot \rangle_{H_{\ell}}$.
For a bounded linear operator
$\cW : H_1 \to H_2$,
the adjoint of $\cW$ is defined as the operator
$\cW^{\ast} : H_2 \to H_1$, 
such that
\begin{equation}
\langle \cW f, g \rangle_{H_2} = \langle f, \cW^{\ast} g \rangle_{H_1}
\quad \mbox{for all $f \in H_1$ and $g \in H_2$}.
\label{eqn:ast1}
\end{equation}
A linear operator $\cW$ is called
an {\it isometry}  if
\[
\|\cW f\|_{H_2} = \|f\|_{H_1} \quad \mbox{for all $f \in H_1$}.
\]The kernel space of $\cW$ is denoted as $\ker \cW$
and its orthogonal complement is written as $(\ker \cW)^{\perp}$.
A linear operator $\cW$ is called 
a {\it partial isometry}, if 
\[
\|\cW f\|_{H_2} = \|f\|_{H_1} \quad \mbox{for all $f \in (\ker \cW)^{\perp}$}.
\]
For the partial isometry $\cW$, 
$(\ker \cW)^{\perp}$ is called the {\it initial space}
and the range of $\cW$, ${\rm ran} \cW$, is called the {\it final space}. 
By the definition (\ref{eqn:ast1}), 
$\|\cW f\|_{H_2}^2=\langle \cW f, \cW f \rangle_{H_2}
=\langle f, \cW^{\ast} \cW f \rangle_{H_1}$. As is suggested 
from this equality, we have the following fact for partial isometries. 
Although this might be known, we give a proof below. 
\begin{lem} \label{lem:A_Astar_partial_isometry2}
Let $H_1$ and $H_2$ 
be separable Hilbert spaces and $\cW : H_1 \to H_2 $ be a bounded
operator. Then, the following are equivalent. 
\begin{description}
\item{\upshape{(i)}} \quad $\cW$ is a partial isometry. 
\item{\upshape{(ii)}} \quad $\cW^* \cW$ is a projection on $H_1$, 
which acts as the identity on $(\ker \cW)^{\perp}$. 
\item{\upshape{(iii)}} \quad $\cW = \cW \cW^* \cW$. 
\end{description}
Moreover, $\cW$ is a partial isometry if and only if so is
 $\cW^*$. 
\end{lem}
\noindent{\it Proof} \,
 When $H_1=H_2$, this fact is well-known
 (cf.~\cite{Hal74}). If we apply it to $H = H_1 \oplus H_2$
 and $\widetilde{\cW} : H \to H$ defined by 
\[
 \widetilde{\cW} := 
\begin{pmatrix}
 O & O \\
 \cW & O
\end{pmatrix}, 
\]
the assertion is followed by verifying 
$\widetilde{\cW^*} \widetilde{\cW} = \cW^* \cW \oplus O$, 
$\widetilde{\cW} \widetilde{\cW^*} = O \oplus \cW \cW^*$ and
\[
\widetilde{\cW} \widetilde{\cW}^* \widetilde{\cW} =
\begin{pmatrix}
 O & O \\
 \cW \cW^* \cW & O
\end{pmatrix}. 
\]
\qed
\vskip 0.3cm
We note that the conditions (i), (ii) and (iii) above, and the conditions 
(i)',  (ii)' and (iii)' obtained 
by applying Lemma~\ref{lem:A_Astar_partial_isometry2} to the
adjoint $\cW^* : H_2 \to H_1$ are all equivalent. 
\vskip 0.3cm
\noindent{\bf Assumption 1} \,
$\cW$ is a partial isometry.
\vskip 0.3cm

By Lemma~\ref{lem:A_Astar_partial_isometry2}, 
under Assumption 1, $\cW^*$ is also a partial
isometry and hence the operator $\cW^{\ast} \cW$ (resp. $\cW\cW^{\ast}$)
is the projection onto the initial space of $\cW$
(resp. the final space of $\cW$).

Now we assume that $H_1$ and $H_2$ are realized as $L^2$-spaces, 
$L^2(S_1, \lambda_1)$ and $L^2(S_2, \lambda_2)$,
respectively. 

A bounded linear operator $\cA : L^2(S_1, \lambda_1) \to
L^2(S_2, \lambda_2)$ 
is a {\it Hilbert--Schmidt operator} 
if Hilbert--Schmidt norm is finite; 
$\|\cA \|_{\mathrm{HS}}^2 :=\tr(\cA^{\ast} \cA) < \infty$.
We say that $\cA$ is a {\it locally Hilbert--Schmidt operator} 
or {\it of locally Hilbert--Schmidt class}, 
if $\cA \cP_{\Lambda}$ is
a Hilbert--Schmidt operator for any compact set
$\Lambda \subset S$. 
It is known as the \textit{kernel theorem} that every
Hilbert--Schmidt operator $\cA : L^2(S_1, \la_1) \to
L^2(S_2, \la_2)$ is defined as an integral operator with kernel 
$A \in L^2(S_1 \times S_2, \la_1 \otimes \la_2)$
 (cf.~Theorem 12.6.2 \cite{Aub00}). In Proposition~\ref{prop:kernel-thm},
we prove a local version of the kernel theorem. 

We put the second assumption.

\vskip 0.3cm
\noindent{\bf Assumption 2} \,
(i) $\cW$ is a locally Hilbert--Schmidt 
operator, and (ii) $\cW^*$ is a locally Hilbert--Schmidt operator. 
\vskip 0.3cm

We note that for any compact set $\La_1 \subset S_1$, 
the operator $\cW \cP_{\La_1}$ is of Hilbert--Schmidt class 
if and only if 
the operator $\cP_{\La_1} \cW^* \cW \cP_{\La_1}$
is of trace class since 
\[
 \|\cW \cP_{\La_1}\|_{\mathrm{HS}}^2 =
\tr\Big( (\cW \cP_{\La_1})^* \cW \cP_{\La_1} \Big)
=\tr\Big( \cP_{\La_1} \cW^* \cW \cP_{\La_1} \Big) < \infty. 
\]
Therefore, Assumption 2 (i) (resp. Assumption 2 (ii)) 
is equivalent to the following
Assumption 2' (i) (resp. Assumption 2' (ii)), 
which guarantees the existence of DPP
associated with $\cW^{\ast} \cW$ (resp. $\cW \cW^{\ast}$).  
\vskip 0.3cm
\noindent{\bf Assumption 2'} \,
(i) $\cW^{\ast} \cW \in \locally(S_1, \lambda_1)$ and 
(ii) $\cW \cW^{\ast} \in \locally(S_2, \lambda_2)$. 
\vskip 0.3cm

Given a measure space $(S, \lambda)$,  
if $f \in L^2(\Lambda, \lambda)$ for all compact
subsets $\Lambda$ of $S$, 
then $f$ is said to be locally $L^2$-integrable.
The set of all such functions is denoted by
$L^2_{\mathrm{loc}}(S, \la)$. 
By this definition if $\cP_{\Lambda} f \in L^2(S, \lambda)$
for any compact set $\Lambda \subset S$, then
$f \in L^2_{\mathrm{loc}}(S, \lambda)$. 
The following proposition is a \textit{local} version of
the kernel theorem for Hilbert--Schmidt operators. 

\begin{prop}
\label{prop:kernel-thm}
Suppose Assumption 
{\rm 2 (i)} holds.
Then, $\cW$ is regarded as an integral operator associated with
a kernel $W : S_2 \times S_1 \to \C$; 
\begin{equation}
(\cW f)(y) = \int_{S_1} W(y, x) f(x) \lambda_1(dx),
\quad f \in L^2(S_1, \lambda_1), 
\label{eqn:W1}
\end{equation}
such that 
$\Psi_1 \in L^2_{\mathrm{loc}}(S_1, \la_1)$, where
$\Psi_1(x) := \|W(\cdot, x)\|_{L^2(S_2, \la_2)}, x \in S_1$.
\end{prop}
\noindent{\it Proof} \,
From Assumption 2 (i)
and the kernel theorem 
for Hilbert--Schmidt operators
(cf. Theorem 12.6.2 \cite{Aub00}), for each compact set $\La
\subset S_1$, there exists 
a kernel $W_{\La} \in L^2(S_2 \times S_1, \la_2 \otimes
\la_1)$ such that 
\[
 \cW\cP_{\La}f(y) = \int_{S_1} W_{\La}(y,x) f(x) \la_1(dx). 
\]
Since $\cW\cP_{\La}f(y) =0$ for 
all $f \in L^2(S_1, \la_1)$ whose support is
contained in $\La^{\mathrm{c}}$, $W_{\La}(y,x) = 0$ on $S_2 \times
\La^{\mathrm{c}}$ for $\la_2 \otimes \la_1$-a.e.$(y,x)$. 
We take two compact sets $\La$ and $\La'$ with $\La \subset
\La' \subset S_1$. For any $f \in L^2(S_1, \la_1)$ whose 
support is contained in $\La$, we see that $\cW\cP_{\La'} f =
\cW\cP_{\La'} \cP_{\La} f = \cW\cP_{\La} f$. Hence, for any
$\La \subset \La' \subset S_1$, 
\[
W_{\La'}(y,x) {\bf 1}_{\Lambda}(x) = W_{\La}(y,x) \quad 
\text{on $S_2 \times S_1$ for $\la_2 \otimes \la_1$-a.e.$(y,x)$}. 
\]
From this consistency, we can define $W(y, x) \in
L^2_{\mathrm{loc}}(S_2\times S_1, \la_2 \otimes \la_1)$ so
that for any compact set $\La \subset S_1$,  
\begin{equation}
 W_{\La}(y,x) = W(y,x) {\mathbf 1}_{\La}(x) 
\quad \text{on $S_2 \times S_1$ for $\la_2
\otimes \la_1$-a.e.$(y,x)$}. 
\label{eq:consistency}
\end{equation}
Since $W_{\La} \in L^2(S_2 \times S_1, \la_2
\otimes\la_1)$, by \eqref{eq:consistency}, 
$\|W(\cdot,x) {\mathbf 1}_{\La}(x)\|_{L^2(S_2, \la_2)} =
\Psi_1(x){\mathbf 1}_{\La}(x)$ is finite $\la_1$-a.e.$x$,
and also 
$\Psi_1 {\mathbf 1}_{\La} \in L^2(S_1, \la_1)$. 
This means that $\Psi_1 \in L^2_{\mathrm{loc}}(S_1,\la_1)$.   
This completes the proof. 
\qed
\vskip 0.3cm
From Proposition~\ref{prop:kernel-thm}, 
under Assumption 2 (ii), the dual operator
$\cW^{\ast}$ also admits an integral kernel 
$W^* : S_1 \times S_2 \to \C$ such that 
$\Psi_2 \in L^2_{\mathrm{loc}}(S_2, \la_2)$, where
$\Psi_2(y) := \|W^*(\cdot, y)\|_{L^2(S_1, \la_1)}, y \in S_2$.
It is easy to see that 
$W^{\ast}(x,y) = \overline{W(y,x)}$ for $\la_1 \otimes\la_2$-a.e.$(x,y)$. 
Then
\begin{equation}
(\cW^{\ast} g)(x)=\int_{S_2} \overline{W(y, x)} g(y) \lambda_2(dy),
\quad g \in L^2(S_2, \lambda_2).
\label{eqn:W1*}
\end{equation}

Following (\ref{eqn:W1}) and (\ref{eqn:W1*}), 
we have
\begin{align*}
(\cW^{\ast} \cW f)(x) &= \int_{S_1} K_{S_1}(x, x') f(x') \lambda_1(dx'),
\quad f \in L^2(S_1, \lambda_1),
\nonumber\\
(\cW \cW^{\ast} g)(y) &= \int_{S_2} K_{S_2}(y, y') g(y') \lambda_2(dy'),
\quad g \in L^2(S_2, \lambda_2), 
\end{align*}
with the integral kernels, 
\begin{align}
K_{S_1}(x, x') &= \int_{S_2} \overline{W(y, x)} W(y, x') \lambda_2(dy)
=\langle W(\cdot, x'), W(\cdot, x) \rangle_{L^2(S_2, \lambda_2)},
\nonumber\\
K_{S_2}(y,y') &= \int_{S_1} W(y,x) \overline{W(y', x)} \lambda_1(dx)
=\langle W(y, \cdot), W(y', \cdot) \rangle_{L^2(S_1, \lambda_1)}.
\label{eqn:K1}
\end{align}
We see that $\overline{K_{S_1}(x', x)}=K_{S_1}(x, x')$
and $\overline{K_{S_2}(y', y)}=K_{S_2}(y, y')$.

Under Assumptions 1 and 2 (i), we obtain 
Theorem~\ref{thm:exist_DPP_HSversion} in Introduction 
as an immediate consequence of the well-known existence theorem
of DPP (Theorem \ref{thm:exist_DPP}). We also state the 
following theorem to emphasize duality of DPPs, and 
it is a starting-point for our discussion in the present paper. 
\begin{thm}\label{thm:main1}
Under Assumptions {\rm 1} and {\rm 2}, 
associated with $\cW^{\ast} \cW$ and $\cW \cW^{\ast}$,
there exists a unique pair of DPPs; 
$(\Xi_1, K_{S_1}, \lambda_1(d x))$
on $S_1$ and
$(\Xi_2, K_{S_2}, \lambda_2(d y))$
on $S_2$.
The correlation kernels 
$K_{S_{\ell}}, \ell=1,2$ are 
Hermitian and given by {\rm (\ref{eqn:K1})}.
\end{thm}
\vskip 0.3cm
Note that the densities of
the DPPs, $(\Xi_1, K_{S_1}, \lambda_1(d x))$ 
and $(\Xi_2, K_{S_2}, \lambda_2(d y))$,
are given by
\begin{align*}
\rho_1(x) &=K_{S_1}(x,x) 
=\int_{S_2} |W(y,x)|^2 \lambda_2(dy)
=\|W(\cdot, x)\|_{L^2(S_2,\lambda_2)}^2, 
\quad x \in S_1, 
\nonumber\\
\rho_2(y) &= K_{S_2}(y, y) 
=\int_{S_1} |W(y,x)|^2 \lambda_1(dx)
=\|W(y, \cdot)\|_{L^2(S_1,\lambda_1)}^2,
\quad y \in S_2,
\end{align*}
with respect to the background measures $\lambda_1(dx)$
and $\lambda_2(dy)$, respectively.

We say that 
a pair of DPPs $(\Xi_1, K_{S_1}, \lambda_1(d x))$
on $S_1$ and $(\Xi_2, K_{S_2}, \lambda_2(d y))$
on $S_2$ is associated with $\cW$.
One of the advantages of our framework is that 
the obtained pairs of DPPs satisfy 
useful duality relations,
which will be reported in Sections~\ref{sec:duality}, 
\ref{sec:duality_application1}, and \ref{sec:duality_application2}.
Now we concentrate on one of a pair of DPPs
constructed in our framework, $(\Xi_1, K_{S_1}, \lambda_1)$.
The correlation kernel $K_{S_1}$ is given by the first equation
of (\ref{eqn:K1}), which is an integral kernel for $f \in L^2(S_1, \lambda_1)$.
We can regard this equation as a {\it decomposition formula}
of $K_{S_1}$ by a product of $W$ and $\overline{W}$.
Since $W$ is an integral kernel for an isometry
$L^2(S_1, \lambda_1) \to L^2(S_2, \lambda_2)$, 
as a matter of course, it depends on a choice of
another Hilbert space $L^2(S_2, \lambda_2)$.
We note that a given DPP, $(\Xi_1, K_{S_1}, \lambda_1)$, 
choice of $L^2(S_2, \lambda_2)$ is not unique. 
Such multivalency gives plural different expressions
for one correlation kernel $K_{S_1}$ and they reveal
different aspects of the DPP 
as demonstrated in Sections
\ref{sec:three_Ginibre}--\ref{sec:Fock}.

\subsection{Basic properties of DPPs}
\label{sec:basic_DPP}

For $v=(v^{(1)}, \dots, v^{(d)}) \in \R^d$,
$y=(y^{(1)}, \dots, y^{(d)}) \in \R^d$, $d \in \N$,
the inner product of them is given by
$v \cdot y = y \cdot v
: =\sum_{a=1}^d v^{(a)} y^{(a)}$,
and $|v|^2 := v \cdot v$.
When $S \subset \C^d, d \in \N$,
$x \in S$ has $d$ complex components;
$x=(x^{(1)}, \dots, x^{(d)})$ with
$x^{(a)}=\Re x^{(a)}+i \Im x^{(a)}$, $a=1, \dots, d$.
In order to describe clearly such a complex structure, 
we set
$x_{\rR}=(\Re x^{(1)}, \dots, \Re x^{(d)}) \in \R^d$,
$x_{\rI}=(\Im x^{(1)}, \dots, \Im x^{(d)}) \in \R^d$,
and write $x=x_{\rR}+i x_{\rI}$ in this paper.
The Lebesgue measure is written as
$d x = d x_{\rR} d x_{\rI} := \prod_{a=1}^d d \Re x^{(a)} d \Im x^{(a)}$.
The complex conjugate of $x=x_{\rR}+i x_{\rI}$
is defined as $\overline{x}=x_{\rR}-i x_{\rI}$.
For $x = x_{\rR}+ i x_{\rI}$,
$x'=x'_{\rR}+ i x'_{\rI} \in \C^d$,
we use the {\it Hermitian inner product}; 
\[
x \cdot \overline{x'} := (x_{\rR}+i x_{\rI}) \cdot (x'_{\rR}-i x'_{\rI})
=(x_{\rR} \cdot x'_{\rR}+x_{\rI} \cdot x'_{\rI})-i(x_{\rR} \cdot x'_{\rI}- x_{\rI} \cdot x'_{\rR})
\]
and define
\[
|x|^2 := x \cdot \overline{x} = |x_{\rR}|^2+|x_{\rI}|^2,
\quad x \in \C^d.
\]

For $(\Xi, K, \lambda(d x))$ defined 
on $S=\R^d$, $S=\C^d$, 
or on a space having appropriate 
periodicities or symmetries, 
we write $\Xi = \sum_j \delta_{X_j}$ and
introduce the following operations. 
\begin{description}
\item{\bf (Shift)} 
For $u \in S$, 
$\cS_{u} \Xi :=\sum_j \delta_{X_j - u} $,
\[
\cS_{u} K(x, x')=K(x+u, x'+u), 
\]
and $\cS_{u} \lambda(d x)=\lambda(u+d x)$.
We write $(\cS_{u} \Xi, \cS_{u} K, \cS_{u} \lambda(d x))$ 
simply as $\cS_{u} (\Xi, K, \lambda(d x))$.
\item{\bf (Dilatation)} 
For $c>0$, we set
$c \circ \Xi :=\sum_j \delta_{c X_j}$
\[
c \circ K(x,x') :=
K \left( \frac{x}{c}, \frac{x'}{c} \right), 
\quad x, x' \in c S:=\{ c x : x \in S\}, 
\]
and $c \circ \lambda(d x):=\lambda(d x/c)$.
We define $c \circ (\Xi, K, \lambda(d x)) := (c \circ \Xi, c \circ K, c \circ \lambda(d x))$.
\end{description}

Moreover, we also consider the following operations.
\begin{description}
\item{\bf (Square root)}
For $(\Xi, K, \lambda(dx))$ on $S=[0, \infty)$, we put
$\Xi^{\langle 1/2 \rangle} :=\sum_j \delta_{\sqrt{X_j}}$, 
$K^{\langle 1/2 \rangle}(x, x') :=K(x^2, {x'}^2)$,
and $\lambda^{\langle 1/2 \rangle}(d x) :=(\lambda \circ v^{-1}) (d x)$,
where $v(x)=\sqrt{x}$.
We define $(\Xi, K, \lambda(d x))^{\langle 1/2 \rangle}$
$:= (\Xi^{\langle 1/2 \rangle}, K^{\langle 1/2 \rangle}, 
\lambda^{\langle 1/2 \rangle}(d x))$
on $[0, \infty)$.
\item{\bf (Gauge transformation)} 
For non-vanishing $u: S \to \C$, 
a gauge transformation of $K$ by $u$ is defined as
\[
K(x, x') \mapsto \widetilde{K}_u(x, x') :=u(x) K(x, x') u(x')^{-1}.
\]
In particular, when $u: S \to \U(1)$,
the $\U(1)$-gauge transformation of $K$ is given by
\[
K(x, x') \mapsto \widetilde{K}_u(x, x') 
= u(x) K(x, x') \overline{u(x')}.
\]
\end{description}

We will use the following basic properties of DPP.
\begin{description}
\item{\bf [Gauge invariance]} 
For any $u: S \to \C$, a gauge transformation
does not change the probability law of DPP;
\[
(\Xi, K, \lambda(d x)) \law (\Xi, \widetilde{K}_u, \lambda(d x)).
\]

\item{\bf [Measure change]}
For a measurable function $g : S \to [0, \infty)$,
\begin{equation}
(\Xi, K(x, x^{\prime}), g(x) \lambda(dx))
\law
(\Xi, \sqrt{g(x)} K(x, x^{\prime}) \sqrt{g(x^{\prime})}, \lambda(dx)).
\label{eqn:measure_change}
\end{equation}

\item{\bf [Mapping and scaling]} 
For a one-to-one measurable mapping
$h: S \to \widehat{S}$, if we set
\[
\widehat{\Xi} :=\sum_j \delta_{h(X_j)},
\quad
\widehat{K}(x, x') :=K(h^{-1}(x), h^{-1}(y)),
\quad 
\widehat{\lambda}(dx) :=(\lambda \circ h^{-1})(dx), 
\]
then 
$(\widehat{\Xi}, \widehat{K}, \widehat{\lambda}(dx))$
is a DPP on $\widehat{S}$.
In particular, 
when $h(x)=x-u, u \in S$,
$(\widehat{\Xi}, \widehat{K}, \widehat{\lambda}(dx))
=\cS_u (\Xi, K, \lambda(dx))$,
when $h(x)=cx, c > 0$,
$(\widehat{\Xi}, \widehat{K}, \widehat{\lambda}(dx))
=c \circ (\Xi, K, \lambda(dx))$,
and when $h(x)=\sqrt{x}$ for $S=[0, \infty)$, 
$(\widehat{\Xi}, \widehat{K}, \widehat{\lambda}(dx))
=(\Xi, K, \lambda(dx))^{\langle 1/2 \rangle}$.
If $c \circ \lambda(d x) = c^{-d} \lambda(d x)$, 
then (\ref{eqn:measure_change}) with $g(x) \equiv c >0$ gives 
\[
c \circ (\Xi, K, \lambda(d x)) \law (c \circ \Xi, K_c, \lambda(d x)),
\quad c >0, 
\]
with
\[
K_c(x, x') := 
\frac{1}{c^d} K \left( \frac{x}{c}, \frac{x'}{c} \right),
\]
where the base space is given by $cS$.
\end{description}

We will give some limit theorems for DPPs in this paper.
Consider a DPP which depends on a continuous
parameter, or a series of DPPs labeled by
a discrete parameter (e.g., the number of points $N \in \N$),
and describe the system by 
$(\Xi, K_p, \lambda_p(dx))$ with 
the continuous or discrete parameter $p$.
If $(\Xi, K_p, \lambda_p(dx))$ converges
to a DPP, $(\Xi, K, \lambda(dx))$, as $p \to \infty$, 
weakly in the vague topology,
we write this limit theorem as
$(\Xi, K_p, \lambda_p(dx)) \weakp (\Xi, K, \lambda(d x))$.
The weak convergence of DPPs is verified by 
the uniform convergence of the kernel $K_p \to K$
on each compact set $C \subset S \times S$ \cite{ST03a}.

\subsection{Duality relations}
\label{sec:duality}
For $f \in \cC_{\rm c}(S)$, the Laplace transform
of the probability measure $\bP$ for a point 
process $\Xi$ is defined as
\begin{equation}
\Psi[f] := \bE \left[
\exp \left( \int_{S} f(x) \Xi(d x) \right) \right].
\label{eqn:Laplace1}
\end{equation}
For the DPP, $(\Xi, K, \lambda(d x))$,
this is given by the Fredholm determinant
on $L^2(S, \lambda)$ \cite{Sim05},
\begin{align*}
\Det_{L^2(S, \lambda)} [ I - (1-e^f) \cK]
& := 1+\sum_{n \in \N}
\frac{(-1)^n}{n!} \int_{S^n} \det_{1 \leq j, k \leq n} [K(x_j, x_k)]
\prod_{\ell=1}^n (1-e^{f(x_{\ell})}) \lambda^{\otimes n}(d \x).
\end{align*}

\begin{lem}
\label{thm:duality}
Between two DPPs, 
$(\Xi_1, K_{S_1}, \lambda_1(d x))$
on $S_1$ and
$(\Xi_2, K_{S_2}, \lambda_2(d y))$
on $S_2$, given by Theorem {\rm \ref{thm:main1}},
the following equality holds with an arbitrary parameter $\alpha \in \C$,
\begin{equation}
\Det_{L^2(S_1, \lambda_1)} [I+\alpha \cK_{S_1} ]
=\Det_{L^2(S_2, \lambda_2)} [I+\alpha \cK_{S_2} ].
\label{eqn:duality}
\end{equation}
\end{lem}
\noindent{\it Proof} \quad 
We recall that if $\cA \cB$ and $\cB \cA$ are trace class operators on
a Hilbert space $H$ then \cite{Sim05}
\begin{equation}
 \Det_H [I+ \cB \cA] = \Det_H [I + \cA \cB]. 
\label{eqn:commutative} 
\end{equation}
Now we have $\cA : H_1 \to H_2$ and $\cB : H_2 \to H_1$ between  
two Hilbert spaces $H_1$ and $H_2$. 
Let $\widetilde{\cA}$ and $\widetilde{\cB}$ be two operators on 
$H_1\oplus H_2$ defined by 
\[
 \widetilde{\cA}
:=
\begin{pmatrix}
O & O \\
\cA & O \\
\end{pmatrix}, \quad 
\widetilde{\cB} 
:= \begin{pmatrix}
O & \cB \\
O & O \\
\end{pmatrix}
\]
Then, $\widetilde{\cA} \widetilde{\cB}$ and $\widetilde{\cB} \widetilde{\cA}$ are diagonal operators 
$O \oplus \cA \cB$ and $\cB \cA \oplus O$, respectively, and hence 
also they are trace class operators. 
By applying (\ref{eqn:commutative}) to $\widetilde{\cA}$ and
$\widetilde{\cB}$ with $H := H_1 \oplus H_2$, we obtain 
\[
 \Det_{H_1} [I+\cB \cA] = \Det_{H_2} [I+ \cA \cB]. 
\]
Consequently, taking $\cA = \sqrt{\alpha} \cW$, $\cB=\sqrt{\alpha} \cW^*$, 
$H_1=L^2(S_1, \lambda_1)$, and $H_2=L^2(S_2, \lambda_2)$ yields 
(\ref{eqn:duality}). 
\qed
\vskip 0.3cm

For $\Lambda_{\ell} \subset S_{\ell}, \ell=1,2$, let
\begin{equation}
\widetilde{\cW} := \cP_{\Lambda_2} \cW \cP_{\Lambda_1},
\quad
\cK_{S_1}^{(\Lambda_2)}
:= \cW^{\ast} \cP_{\Lambda_2} \cW,
\quad
\cK_{S_2}^{(\Lambda_1)}
:= \cW \cP_{\Lambda_1} \cW^{\ast}.
\label{eqn:defK_A}
\end{equation}
They admit the following integral kernels, 
\begin{align}
\widetilde{W}(y, x)
&={\bf 1}_{\Lambda_2}(y) W(y, x) {\bf 1}_{\Lambda_1}(x),
\nonumber\\
K_{S_1}^{(\Lambda_2)}(x, x^{\prime})
&= \int_{\Lambda_2} \overline{W(y, x)} W(y, x^{\prime}) \lambda_2(d y),
\nonumber\\
K_{S_2}^{(\Lambda_1)}(y, y^{\prime})
&= \int_{\Lambda_1} W(y, x) \overline{W(y^{\prime}, x)} \lambda_1(d x).
\label{eqn:K_A}
\end{align}
Using Lemma \ref{thm:duality}, the following theorem is proved.
\begin{thm}
\label{thm:duality2}
Let $(\Xi_1^{(\Lambda_2)}, K_{S_1}^{(\Lambda_2)}, \lambda_1(d x))$
and $(\Xi_2^{(\Lambda_1)}, K_{S_2}^{(\Lambda_1)}, \lambda_2(d y))$ 
be DPPs associated with the kernels
$K_{S_1}^{(\Lambda_2)}$ and $K_{S_2}^{(\Lambda_1)}$
given by {\rm (\ref{eqn:K_A})}, respectively.
Then, $\Xi_1^{(\Lambda_2)}(\Lambda_1) \law
 \Xi_2^{(\Lambda_1)}(\Lambda_2)$, i.e., 
\[
\bP( \Xi_1^{(\Lambda_2)}(\Lambda_1)=m)
= \bP( \Xi_2^{(\Lambda_1)}(\Lambda_2)=m),
\quad \forall m \in \N_0.
\]
\end{thm}
\noindent{\it Proof} \,
As a special case of (\ref{eqn:Laplace1})
with $f(x)={\bf 1}_{\Lambda_1}(x) \log z$
for $\Xi=\Xi_1^{(\Lambda_2)}$,
$z \in \C$, we have the equality,
\begin{equation}
\bE\left[ z^{\Xi_1^{(\Lambda_2)}(\Lambda_1)} \right]
=\Det_{L^2(S_1, \lambda_1)}
[I-(1-z) \cP_{\Lambda_1} \cK_{S_1}^{(\Lambda_2)} \cP_{\Lambda_1} ],
\label{eqn:generating1}
\end{equation}
where $\cK_{S_1}^{(\Lambda_2)}$ is defined by (\ref{eqn:defK_A}). 
Here LHS is the moment generating function
of $\Xi^{(\Lambda_2)}_1(\Lambda_1)$ 
and RHS gives its Fredholm determinantal expression.
By replacing $\cW$ by $\widetilde{\cW}$ and letting $\alpha=-(1-z)$ in the proof of
Lemma \ref{thm:duality}, we obtain the equality,
\[
\Det_{L^2(S_1, \lambda_1)}
[I-(1-z) \cP_{\Lambda_1} \cK_{S_1}^{(\Lambda_2)} \cP_{\Lambda_1} ]
=
\Det_{L^2(S_2, \lambda_2)}
[I-(1-z) \cP_{\Lambda_2} \cK_{S_2}^{(\Lambda_1)} \cP_{\Lambda_2} ].
\]
Through (\ref{eqn:generating1})
and the similar equality for 
$\bE\left[ z^{\Xi_2^{(\Lambda_1)}(\Lambda_2)} \right]$, 
we obtain the corresponding
equivalence between the moment generating functions
of $\Xi_1^{(\Lambda_2)}(\Lambda_1)$ and $\Xi_2^{(\Lambda_1)}(\Lambda_2)$,
and hence the statement of the proposition is proved. \qed

Examples of duality relations will be given in 
Sections \ref{sec:duality_application1} and \ref{sec:duality_application2}.
Theorem \ref{thm:duality2} was used to analyze 
hyperuniformity \cite{Tor18}
of the Heisenberg family of DPPs in \cite{MKS21}.

\subsection{Orthonormal functions and correlation kernels}
\label{sec:orthogonal}

In addition to $L^2(S_{\ell}, \lambda_{\ell})$, $\ell=1, 2$, 
we introduce $L^2(\Gamma, \nu)$ as a parameter space
for functions in $L^2(S_{\ell}, \lambda_{\ell}), \ell=1,2$. 
Assume that there are two families of measurable
functions $\{\psi_1(x, \gamma) : x \in S_1, 
\gamma \in \Gamma \}$ and
$\{\psi_2(y, \gamma) : y \in S_2, 
\gamma \in \Gamma \}$ 
such that 
two bounded operators $\cU_{\ell} : L^2(S_{\ell}, \lambda_{\ell}) \to L^2(\Gamma, \nu)$
given by 
\[
(\cU_{\ell}f)(\gamma) := 
\int_{S_{\ell}} \overline{\psi_{\ell}(x, \gamma)} f(x)
\lambda_{\ell}(dx), \quad \ell=1,2, 
\]
are well-defined. 
Then, their adjoints $\cU_{\ell}^* : L^2(\Gamma, \nu) \to 
L^2(S_{\ell}, \lambda_{\ell}), \ell =1,2$ 
are given by 
\[
(\cU_{\ell}^* F)(\cdot) = 
\int_{\Gamma} \psi_{\ell}(\cdot, \gamma) F(\gamma)
\nu(d\gamma). 
\]
A typical example of $\cU_1$ is the Fourier transform,
i.e., $\psi_1(x,\gamma) = e^{i x \gamma}$. In this case, for
any $\gamma$, the function $\psi_1(\cdot,\gamma)$ is \textit{not} in $L^2(\R, dx)$.
Now we define $\cW : L^2(S_1, \lambda_1) \to L^2(S_2, \lambda_2)$ by 
$\cW = \cU_2^* \cU_1$, i.e., 
\begin{equation}
 (\cW f)(y) = \int_{\Gamma} \psi_2(y,\gamma) (\cU_1 f)(\gamma)
 \nu(d\gamma). 
\label{eqn:W_of}
\end{equation}
Let $I_{\Gamma}$ be an identity in $L^2(\Gamma, \nu)$.
We can see the following. 
\begin{lem} 
If $\cU_{\ell} \cU_{\ell}^* = I_{\Gamma}$ for $\ell=1,2$, then 
both $\cW$ and $\cW^*$ are partial isometries. 
\end{lem}
\noindent{\it Proof} \,
By the assumption, we see that 
\[
\cW \cW^* \cW = 
(\cU_2^* \cU_1) (\cU_1^* \cU_2) (\cU_2^* \cU_1) 
= \cU_2^* \cU_1 = \cW. 
\]
From Lemma~\ref{lem:A_Astar_partial_isometry2}, $\cW$ is a partial isometry. 
By symmetry, the assertion for $\cW^*$ also follows. 
\qed
\vskip 0.3cm
We note that 
$\cW^* \cW = \cU_1^* \cU_1$ and 
$\cW \cW^* = \cU_2^* \cU_2$. Hence,  
$\cU_{\ell}, \ell=1,2$ are also partial isometries.
In addition, 
$\cW^{\ast} \cW$ is a locally trace
class operator if and only if so is $\cU_1^{\ast} \cU_1$. 
Therefore, $\cW$ is of locally Hilbert--Schmidt class if and only if 
so is $\cU_1$. 

Now we rewrite the condition for 
$\cU_1$ to be of locally Hilbert--Schmidt class in terms of the function
$\psi_1(x, \gamma), x \in S_1, \gamma \in \Gamma$. 
\begin{lem}\label{lem:HScondition}
Let $\Psi_1(x) := \|\psi_1(x, \cdot)\|_{L^2(\Gamma, \nu)}$, $x \in S_1$
and assume that 
$\Psi_1 \in L^2_{\mathrm{loc}}(S_1, \la_1)$. 
Then, the operator $\cU_1$ is of locally Hilbert--Schmidt class. 
\end{lem}
\noindent{\it Proof} \quad
For a compact set $\La \subset S_1$, we see that 
\begin{align*}
|\cP_{\La} \cU_1^* \cU_1 \cP_{\La} f (x)|
&= \left|\mathbf{1}_{\La}(x) 
\int_{\Gamma} \nu(d\gamma) \psi_1(x,\gamma) 
\int_{S_1} \overline{\psi_1(x',\gamma)}
 \mathbf{1}_{\La}(x') f(x')\la_1(dx')\right|  \\
&\le \mathbf{1}_{\La}(x) \Psi_1(x) \int_{S_1} \mathbf{1}_{\La}(x') \Psi_1(x')
|f(x')| \la_1(dx') \\
&\le \cP_{\La} \Psi_1(x) \|\cP_{\La} \Psi_1\|_{L^2(S_1,\la_1)}
\|\cP_{\La} f\|_{L^2(S_1,\la_1)}. 
\end{align*}
By Fubini's theorem, we have 
\[
\cP_{\La} \cU_1^* \cU_1 \cP_{\La} f (x)
= 
\int_{S_1} \la_1(dx') f(x')
\left(\int_{\Gamma} \mathbf{1}_{\La}(x) \psi_1(x,\gamma) 
\overline{\mathbf{1}_{\La}(x') \psi_1(x',\gamma)}
\nu(d\gamma)\right)
\]
and hence
\[
\|\cU_1 \cP_{\La}\|^2_{\mathrm{HS}} = \tr (\cP_{\La} \cU_1^* \cU_1
\cP_{\La}) = 
\int_{S_1} \la_1(dx) \mathbf{1}_{\La}(x) \left(\int_{\Gamma}
|\psi_1(x,\gamma)|^2 \nu(d\gamma)\right) 
= \|\cP_{\La} \Psi_1\|^2_{L^2(S_1, \la_1)} < \infty. 
\]
This completes the proof. 
\qed
\vskip 0.3cm

Now we put the following.
\vskip 0.3cm
\noindent{\bf Assumption 3} \,
For $\ell=1,2$, 
\begin{description}
\item{(i)} \quad $\cU_{\ell} \cU_{\ell}^* = I_{\Gamma}$,
\item{(ii)} \quad $\Psi_{\ell} \in L^2_{\mathrm{loc}}(S_{\ell}, \la_{\ell})$,
where 
$\Psi_{\ell}(x):=\|\psi_{\ell}(x,\cdot)\|_{L^2(\Gamma, \nu)}$,
$x \in S_{\ell}$. 
\end{description}
\vskip 0.3cm
\noindent
Assumption 3(i) can be rephrased as the following 
orthonormality relations: 
\[
\langle 
\psi_{\ell}(\cdot, \gamma), \psi_{\ell}(\cdot, \gamma') \rangle_{L^2(S_{\ell}, \lambda_{\ell})}
\nu(d \gamma)
=\delta(\gamma-\gamma') d \gamma, 
\quad \gamma, \gamma' \in \Gamma,
\quad \ell=1,2.
\]
We often use these relations below. 
\vskip 0.3cm
The following is immediately obtained as a corollary of
Theorem \ref{thm:main1}.
\begin{cor}
\label{thm:main2}
Let $\cW = \cU_2^* \cU_1$ as in the above. 
We assume Assumption {\rm 3}. Then, there exists a unique pair of DPPs;
$(\Xi_1, K_{S_1}, \lambda_1(d x))$
on $S_1$ and
$(\Xi_2, K_{S_2}, \lambda_2(d y))$
on $S_2$.
Here the correlation kernels 
$K_{S_{\ell}}, \ell=1,2$ are given by 
\begin{align}
K_{S_1}(x, x')
&= \int_{\Gamma} \psi_1(x, \gamma) 
\overline{\psi_1(x', \gamma)} \nu(d \gamma)
=\langle \psi_1(x, \cdot), \psi_1(x', \cdot) \rangle_{L^2(\Gamma, \nu)},
\nonumber\\
K_{S_2}(y, y')
&= \int_{\Gamma} \psi_2(y, \gamma) 
\overline{\psi_2(y', \gamma)} \nu(d \gamma)
=\langle \psi_2(y, \cdot), \psi_2(y', \cdot) \rangle_{L^2(\Gamma, \nu)}.
\label{eqn:K_main2}
\end{align}
In particular, the densities of the DPPs are given by
$\rho_1(x)=K_{S_1}(x,x) = \Psi_1(x)^2, x \in S_1$ 
and $\rho_2(y)=K_{S_2}(y,y)=\Psi_2(y)^2, y \in S_2$ 
with respect to the background measures
$\lambda_1(dx)$ and $\lambda_2(dy)$, respectively.
\end{cor}
\vskip 0.3cm
\noindent{\bf Remark 1} \, 
Consider the symmetric case such that 
$L^2(S_1, \lambda_1)=L^2(S_2, \lambda_2)=: L^2(S, \lambda)$, 
$\psi_1=\psi_2 =: \psi$, 
$\nu=\lambda|_{\Gamma}$, $\Gamma \subseteq S$. 
In this case, $\cW = \cU^* \cU$ with 
\[
 (\cU f)(\gamma) = \int_S  \overline{\psi(x,\gamma)} f(x)
 \lambda(dx). 
\]
Then 
$K_{S_1}=K_{S_2}=W = : K$
is given by
\begin{equation}
K(x, x')=\int_{\Gamma} \psi(x, \gamma) \overline{\psi(x', \gamma)} \lambda(d \gamma).
\label{eqn:K_simple1}
\end{equation}
This is Hermitian; 
$\overline{K(x', x)}=K(x, x')$, and satisfies
the reproducing property 
\[
K(x, x')=\int_{S} K(x, \zeta) K(\zeta, x') \lambda(d \zeta).
\]
\vskip 0.3cm

Now we consider a simplified version of 
the preceding setting. 
Let $\Gamma \subseteq S_2$ and $\nu=\lambda_2|_{\Gamma}$.
We define 
$\cU_2 : L^2(S_2,\lambda_2) \to L^2(\Gamma, \nu)$ as the
restriction onto $\Gamma$, and then its adjoint $\cU_2^*$ is given by 
$(\cU_2^* F)(y) = F(y)$ for $y \in \Gamma$, and by $0$ for $y \in
S_2 \setminus \Gamma$. 
We write the extension $\tilde{F} = \cU^*_2 F$ for $F \in L^2(\Gamma,
\nu)$. 
It is obvious that $\cU_2 \cU_2^* =
I_{\Gamma}$ and hence $\cU_2$ is a partial isometry. 

For $\Gamma \subseteq S_2$, we assume that there is a family
of measurable functions 
$\{\psi_{1}(x, y) : x \in S_1, y \in \Gamma \}$ 
such that a bounded operator 
$\cU_1 : L^2(S_1,\lambda_1) \to L^2(\Gamma, \nu)$ given by 
\[
 (\cU_1 f)(\gamma) := \int_{S_1} \overline{\psi_1(x,\gamma)}
 f(x) \lambda_1(dx) \quad (\gamma \in \Gamma)
\]
is well-defined. 

\vskip 0.3cm
\noindent{\bf Assumption 3'} \,
\begin{description}
\item{(i)} \quad $\cU_1 \cU_1^* = I_{\Gamma}$,
\item{(ii)} \quad $\Psi_1 \in L^2_{\mathrm{loc}}(S_1, \la_1)$,
where 
$\Psi_1(x):=\|\psi_1(x,\cdot)\|_{L^2(\Gamma, \nu)}$,
$x \in S_1$. 
\end{description}
\vskip 0.3cm
\noindent
Assumption 3'(i) can be rephrased as the following 
orthonormality relation: 
\[
 \langle 
 \psi_{1}(\cdot, y), \psi_{1}(\cdot, y') \rangle_{L^2(S_{1}, \lambda_{1})}
 \lambda_2(d y)
 =\delta(y-y') d y, 
 \quad y, y' \in \Gamma.
 \]
Now we define $\cW : L^2(S_1, \lambda_1) \to L^2(S_2, \lambda_2)$ by 
$\cW = \cU_2^* \cU_1$ as before. In this case, we have  
\[
 (\cW f)(y) 
= \mathbf{1}_{\Gamma}(y) \int_{S_1}
\overline{\tilde{\psi}_1(x,y)}f(x) \lambda_1(dx),  
\]
and hence 
\begin{equation}
W(y, x)= \overline{\tilde{\psi}_1(x, y)} {\bf 1}_{\Gamma}(y).
\label{eqn:w_set2}
\end{equation}
It follows from Assumption 3' that $\cW$ is a partial isometry.  
Corollary \ref{thm:main2} is reduced to the following.
\begin{cor}
\label{thm:main3}
Let $\cW = \cU_2^* \cU_1$ as in the above. 
We assume Assumption {\rm 3'}. 
Then there exists a unique DPP, 
$(\Xi, K, \lambda_{1})$
on $S_{1}$ with
the correlation kernel
\begin{equation}
K_{S_1}(x, x^{\prime})
= \int_{\Gamma} \psi_1(x, y) \overline{\psi_1(x^{\prime}, y)} \lambda_2(d y)
=\langle \tilde{\psi}_1(x, \cdot), \tilde{\psi}_1(x^{\prime}, \cdot) \rangle_{L^2(\Gamma, \lambda_2)}.
\label{eqn:K_simple2}
\end{equation}
In particular, the density of the DPP is given by
$\rho_1(x)=K_{S_1}(x,x) = \Psi_1(x)^2, x \in S_1$ 
with respect to the background measures
$\lambda_1(dx)$.
\end{cor}
\noindent{\it Proof} \quad
The proof is the same as before. 
\qed
\vskip 0.3cm
\noindent{\bf Remark 2} \, 
The correlation kernel (\ref{eqn:K_simple2})
is the same as the correlation kernel (\ref{eqn:K_simple1}) 
shown in Remark 1 in the symmetric case, 
$L^2(S_1, \lambda_1)=L^2(S_2, \lambda_2)$
and $\psi_1=\psi_2$, of the pair of
DPPs given by Corollary \ref{thm:main2}.

\subsection{Weyl--Heisenberg ensembles of DPPs}
\label{sec:WH}

The family of DPPs given by Corollary \ref{thm:main3}
is a generalization of the class of DPPs
called the {\it Weyl--Heisenberg ensembles}
studied by Abreu {\it et al.} \cite{AGR16,APRT17,AGR19}. 
For $d \in \N$, let
\[
S_1=\C^d,
\quad
S_2=\Gamma= \R^d,
\]
with the Lebesgue measures
$\lambda_1(dx)=dx_{\rR} dx_{\rI}$, $\lambda_2(dy)=dy$, 
where $x=x_{\rR}+i x_{\rI}$ with $x_{\rR}, x_{\rI} \in \R^d$.
We consider the case that  
$\psi_1$ in the setting (\ref{eqn:w_set2}) of $W$ 
is given of the form
\begin{equation}
\psi_1(x, y)= \psi_1(x_{\rR}+ix_{\rI}, y) = G(y-x_{\rR}) e^{2 \pi i y \cdot x_{\rI}}
\quad \mbox{with $G \in L^2(\R^d, dx_{\rR})$}, 
\label{eqn:g_WH}
\end{equation}
where $y \cdot x_{\rI}$ denotes the inner product in $\R^d$. 
In this setting, 
since $\Psi_1(x) \equiv \|G\|_{L^2(\R^d,dx_{\rR})}$, $x \in \C^d$, 
Assumption 3'(ii) is
satisfied. Also, we have
\begin{align*}
\langle \psi_1(\cdot, y), \psi_1(\cdot, y') \rangle_{L^2(S_1, \lambda_1)}
&= \int_{\R^d} dx_{\rR} \, G(y-x_{\rR}) \overline{G(y'-x_{\rR})} 
\int_{\R^d} dx_{\rI} \, e^{2 \pi i (y-y') \cdot x_{\rI}} 
\end{align*}
Since 
$\int_{\R^d} e^{2 \pi i y \cdot x} dx =
\delta(y), y \in \R^d$, 
the above is equal to $\|G\|_{L^2(\R^d, dx_{\rR})}^2 \, \delta(y-y')$. 
Hence, the norm $\|G\|_{L^2(\R^d, dx_{\rR})}$ must be $1$ for Assumption 3'(i). 
Therefore, in this setting (\ref{eqn:g_WH}), 
Assumption 3' will be reduced to the following.

\vskip 0.3cm
\noindent
{\bf Assumption 4} \, 
The function $G$ in (\ref{eqn:g_WH}) has norm 1
in $L^2(\R^d, d x_{\rR})$.
\vskip 0.3cm

Under the setting (\ref{eqn:w_set2}) with $\Gamma=\R^d$ and (\ref{eqn:g_WH}), 
if Assumption 4 is satisfied, then the operator $\cW$
and the correlation kernel $K_{S_1}$ are written as
\begin{align}
(\cW_{\rm WH} f)(y)
&=\int_{\C^d} \overline{G(y-x_{\rR})} e^{-2 \pi i y \cdot x_{\rI}} f(x_{\rR}+i x_{\rI}) dx_{\rR} dx_{\rI},
\quad f \in L^2(\C^d, dx_{\rR} dx_{\rI}), 
\nonumber\\
(\cW^{\ast}_{\rm WH} g)(x)
&=\int_{\R^d} G(y-x_{\rR}) e^{2 \pi i y \cdot x_{\rI}} g(y) dy,
\quad g \in L^2(\R^d, dy), 
\nonumber\\
K_{\rm WH}(x, x')
&=\int_{\R^d} 
G(y-x_{\rR}) \overline{G(y-x'_{\rR})}
e^{2 \pi i y \cdot (x_{\rI}-x'_{\rI})} dy,
\label{eqn:K_WH}
\end{align}
for $(x, x')=(x_{\rR}+i x_{\rI}, x'_{\rR}+ i x'_{\rI} )
\in \C^d \times \C^d$. 
The second formula in (\ref{eqn:K_WH}) is
regarded as the {\it short-time Fourier transform}
of $g \in L^2(\R^d,dy)$ with respect to
a {\it window function} $G \in L^2(\R^d,dx_{\rR})$ \cite{Gro01}.
The formulas (\ref{eqn:K_WH}) 
define the Weyl--Heisenberg ensemble of DPP,
$(\Xi, K_{\rm WH}, d x_{\rR} d x_{\rI})$, studied
in \cite{AGR16,APRT17,AGR19}.

\begin{prop}
\label{thm:WH}
Under Assumption {\rm 4}, 
the Weyl--Heisenberg class of DPPs 
specified by the window function $G \in L^2(\R^d, d x_{\rR})$ 
is a special case 
of the family of DPPs given by Corollary {\rm \ref{thm:main3}},
in which $\Gamma=\R^d$,
$S_1=\C^d$, 
$\lambda_1(d x)=d x_{\rR} d x_{\rI}$, $\lambda_2(d y)= d y$, 
and $\psi_1$ in 
{\rm (\ref{eqn:w_set2})} is given of the form 
{\rm (\ref{eqn:g_WH})}.
\end{prop}

\SSC
{Examples in One-dimensional Spaces} \label{sec:examples_1d}
\subsection{Finite DPPs in $\R$ associated with 
classical orthonormal polynomials}
\label{sec:classical}
\subsubsection{Classical orthonormal polynomials and DPPs}
\label{sec:OP}
Let $S_1=S_2 = \R$. 
Assume that we have
two sets of orthonormal functions 
$\{\varphi_n \}_{n \in \N_0}$ 
and $\{\phi_n \}_{n \in \N_0}$
with respect to the measures $\lambda_1$ and $\lambda_2$,
respectively,
\begin{align}
\langle \varphi_n, \varphi_m \rangle_{L^2(\R, \lambda_1)}
&=
\int_{\R} \varphi_n(x) \overline{\varphi_m(x)} \lambda_1(d x) = \delta_{nm},
\nonumber\\
\langle \phi_n, \phi_m \rangle_{L^2(\R, \lambda_2)}
&=
\int_{\R} \phi_n(y) \overline{\phi_m(y)} \lambda_2(d y) = \delta_{nm},
\quad n, m \in \N_0.
\label{eqn:ON1}
\end{align}
Then for an arbitrary but fixed $N \in \N$, 
we set $\Gamma=\{0,1, \dots, N-1\} \subsetneq \N_0$,
$\psi_1(\cdot, \gamma)=\varphi_{\gamma}(\cdot)$,
$\psi_2(\cdot, \gamma)=\phi_{\gamma}(\cdot)$, $\gamma \in \Gamma$, 
and consider $\ell^2(\Gamma)$ as $L^2(\Gamma, \nu)$
in the setting of Section \ref{sec:orthogonal}. 
We see that 
$\int_{\R} \|\varphi_{\cdot}(x)\|^2_{\ell^2(\Gamma)} \lambda_1(dx)
=\sum_{n=0}^{N-1} \|\varphi_n \|^2_{L^2(\R, \lambda_1)}
=N$ 
and
$\int_{\R} \|\phi_{\cdot}(y)\|^2_{\ell^2(\Gamma)} \lambda_2(dy)
=\sum_{n=0}^{N-1} \| \phi_n \|^2_{L^2(\R, \lambda_2)}
=N$. 
Hence Assumption 3 is satisfied for any $N \in \N$.
Then the integral kernel for $\cW$ 
defined by (\ref{eqn:W_of}) is given by
\[
W(y, x)= \sum_{n=0}^{N-1} \overline{\varphi_n(x)} \phi_n(y).
\]
By Corollary \ref{thm:main2}, we have a pair of DPPs on $\R$,
$(\Xi_1, K_{\varphi}^{(N)}, \lambda_1(d x))$
and $(\Xi_2, K_{\phi}^{(N)}, \lambda_2(d y))$, 
where the correlation kernels are given by
\begin{equation}
K_{\varphi}^{(N)}(x, x')=\sum_{n=0}^{N-1} \varphi_n(x) \overline{\varphi_n(x')},
\quad
K_{\phi}^{(N)}(y, y')=\sum_{n=0}^{N-1} \phi_n(y) \overline{\phi_n(y')},
\label{eqn:K_classical}
\end{equation}
respectively.
Here $N$ gives the number of points for each DPPs.
If we can use the three-term relations 
in $\{\varphi_n\}_{n \in \N_0}$ or $\{\phi_n\}_{n \in \N_0}$, 
(\ref{eqn:K_classical}) can be written 
in the Christoffel--Darboux form
(see, for instance, Proposition 5.1.3 in \cite{For10}). 
As a matter of course, if we have three or more than three,
say $M$ distinct sets of orthonormal functions
satisfying Assumption 3 with a common $\Gamma$, then
by applying Corollary \ref{thm:main2} to every pair of them,
we will obtain $M$ distinct finite DPPs.
See examples given in 
Sections \ref{sec:Lie_group}, \ref{sec:affine}, \ref{sec:orthogonal_theta}, 
\ref{sec:rectangle_finite}, \ref{sec:strip_finite_infinite}.

Even if we have only one set of orthonormal functions,
for example, only the first one $\{\varphi_n\}_{n \in \N_0}$ in (\ref{eqn:ON1}),
we can obtain a DPP (labeled by the number of particles $N \in \N$) 
following Corollary \ref{thm:main3}.
In such a case, we set
\begin{equation}
W(n, x)=\overline{\varphi_n(x)} {\bf 1}_{\Gamma}(n)
\label{eqn:W_n_x1}
\end{equation}
with $\Gamma=\{0,1, \dots, N-1\}$
for (\ref{eqn:w_set2}). Then we have a DPP,
$(\Xi, K_{\varphi}^{(N)}, \lambda_1(d x))$.
See examples given in Sections
\ref{sec:sphere_finite}, \ref{sec:HeisenbergDPP}, 
and \ref{sec:DPP_S_d}.
\vskip 0.3cm
\noindent{\bf Remark 3} \, 
If $\Gamma$ is a finite set, $|\Gamma| < \infty$,
and the parameter space is given by $\ell^2(\Gamma)$,
Assumption 3(ii) (resp. Assumption 3'(ii))
is concluded from 3(i) (resp.~3'(i)) as shown below.
Since
$\Psi(x)^2 :=\|\varphi_{\cdot}(x) \|^2_{\ell^2(\Gamma)}
=\sum_{n \in \Gamma} |\varphi_n(x)|^2$, $x \in S$, 
we have
$\int_{S} \Psi(x)^2 \lambda(dx)
= \sum_{n \in \Gamma} \|\varphi_n \|^2_{L^2(S, \lambda)}$.
Then, if $\{\varphi_n\}_{n \in \Lambda}$ are normalized,
the above integral is equal to $|\Gamma| < \infty$.
This implies $\Psi \in L^2(S, \lambda)
\subset L^2_{\mathrm{loc}}(S, \lambda)$.
See finite DPPs given in 
Sections \ref{sec:Lie_group}, \ref{sec:affine}, \ref{sec:sphere_finite}, 
\ref{sec:rectangle_finite}, \ref{sec:strip_finite_infinite}, 
\ref{sec:DPP_S_d}.
\vskip 0.3cm

Now we give classical examples of DPPs
associated with real-valued orthonormal polynomials.
Let $\lambda_{\rN(m, \sigma^2)}(dx)$ denote the
{\it normal distribution}, 
\[
\lambda_{\rN(m, \sigma^2)}(dx)
:= \frac{1}{\sqrt{2 \pi} \sigma} 
e^{-(x-m)^2/(2\sigma^2)} d x,
\quad m \in \R, \quad \sigma >0, 
\]
and 
$\lambda_{\Gamma(a, b)}(d y)$ do the
{\it Gamma distribution}, 
\[
\lambda_{\Gamma(a, b)}(d y)
:= 
\frac{b^a}{\Gamma(a)} y^{a-1} e^{-b y} {\bf 1}_{[0, \infty)}(y) dy,
\quad a >0, \quad b >0, 
\]
with the Gamma function
$\Gamma(z) := \int_0^{\infty} u^{z-1} e^{-u} du, \Re z > 0$.
We set 
\begin{align}
\lambda_1(dx) &= \lambda_{\rN(0,1/2)}(d x)
= \frac{1}{\sqrt{\pi}} e^{-x^2} dx, 
\nonumber\\
\varphi_n(x) &= \frac{1}{\sqrt{2^n n!}} H_n(x), \quad n \in \N_0, 
\label{eqn:HermiteZ1}
\end{align}
and
\begin{align}
\lambda_2(d y) &=\lambda_{\Gamma(\nu+1, 1)}(d y)
= \frac{1}{\Gamma(\nu+1)} y^{\nu} e^{-y} {\bf 1}_{[0, \infty)}(y)dy, 
\nonumber\\
\phi_n(y) &= \phi_n^{(\nu)}(y)=
\sqrt{ \frac{\Gamma(n+1) \Gamma(\nu+1)}{\Gamma(n+\nu+1)}} \, L^{(\nu)}_n(y),
\quad n \in \N_0, 
\label{eqn:LaguerreZ1}
\end{align}
with parameter $\nu \in (-1, \infty)$.
Here 
$\{H_n(x)\}_{n \in \N_0}$ are the {\it Hermite polynomials}, 
\begin{align}
H_n(x) & : =(-1)^n e^{x^2} \frac{d^n}{dx^n} e^{-x^2}
\nonumber\\
&= n! \sum_{k=0}^{[n/2]} 
\frac{(-1)^k (2x)^{n-2k}}{k! (n-2k)!}, \quad n \in \N_0,
\label{eqn:Hermite1}
\end{align}
where $[a]$ denotes the largest integer not greater than 
$a \in \R$, and
$\{L^{(\nu)}_n(x)\}_{n \in \N_0}$ are
the {\it Laguerre polynomials}, 
\begin{align}
L^{(\nu)}_n(x) & : = \frac{1}{n!} x^{-\nu} e^{x} \frac{d^n}{dx^n} \Big(x^{n+\nu} e^{-x} \Big)
\nonumber\\
&= \sum_{k=0}^{n}
\frac{(\nu+k+1)_{n-k}}{(n-k)! k!} (-x)^k,
\quad n \in \N_0, \quad \nu \in (-1, \infty), 
\label{eqn:Laguerre1}
\end{align}
where $(\alpha)_n := \alpha(\alpha+1) \cdots (\alpha+n-1)
= \Gamma(\alpha+n)/\Gamma(\alpha)$, $n \in \N$, 
$(\alpha)_0 := 1$. 
The correlation kernels (\ref{eqn:K_classical}) are written 
in the {\it Christoffel--Darboux form} as,
\[
K^{(N)}_{\varphi}(x, x^{\prime})=K_{\rm Hermite}^{(N)}(x, x')
=\sqrt{\frac{N}{2}}
\frac{\varphi_N(x) \varphi_{N-1}(x')-\varphi_N(x') \varphi_{N-1}(x)}
{x-x^{\prime}},
\quad x, x^{\prime} \in \R,
\]
and
\begin{align*}
K^{(N)}_{\phi}(y, y')
&=K_{\rm Laguerre}^{(\nu, N)}(y, y')
\nonumber\\
&=- \sqrt{N(N+\nu)}
\frac{\phi^{(\nu)}_N(y) \phi^{(\nu)}_{N-1}(y')
-\phi^{(\nu)}_N(y') \phi^{(\nu)}_{N-1}(y)}
{y-y'},
\quad y, y' \in [0, \infty).
\end{align*}
When $x=x'$ or $y=y'$, we make sense of the above formulas
by using L'H\^opital's rule.
The former is called the {\it Hermite kernel}
and the latter is the {\it Laguerre kernel}. 

By definition, 
for a finite DPP $(\Xi, K, \lambda(dx))$ with $N$ points in $S$, 
the probability density with respect to
$\lambda^{\otimes N}(dx_1 \cdots dx_N)$ is given by
$\rho^{N}(x_1, \dots, x_N)
=\det_{1 \leq j, k \leq N}[K(x_j, x_k)], 
\x=(x_1, \dots, x_N) \in S^N$.
Using the {\it Vandermonde determinantal formula}, 
$\det_{1 \leq j, k \leq N}(z_k^{j-1})
=\prod_{1 \leq j < k \leq N}(z_k-z_j)$,
which is also given as the type $\AN$ of
Weyl denominator formula (\ref{eqn:Weyl_denominator})
in Appendix \ref{sec:Weyl_denominators}, 
we can verify that 
the probability densities
of the DPPs
$(\Xi, K^{(N)}_{\rm Hermite}, \lambda_{\rN(0,1/2)}(dx))$
and 
$(\Xi, K^{(N)}_{\rm Laguerre}, \lambda_{\Gamma(\nu+1, 1)}(dy))$
with respect to 
the Lebesgue measures
$d \x =\prod_{j=1}^N d x_j$ and 
$d \y =\prod_{j=1}^N d y_j$ are given
as
\begin{align}
\bp^{(N)}_{\rm Hermite}(\x) 
&= \frac{1}{Z^{(N)}_{\rm Hermite}}
\prod_{1 \leq j < k \leq N} (x_k-x_j)^2
\prod_{\ell=1}^N e^{-x_{\ell}^2},
\quad 
\x \in \R^N,
\nonumber\\
\bp^{(\nu, N)}_{\rm Laguerre}(\y) 
&= \frac{1}{Z^{(\nu, N)}_{\rm Laguerre}}
\prod_{1 \leq j < k \leq N} (y_k-y_j)^2
\prod_{\ell=1}^N y_{\ell}^{\nu} e^{-y_{\ell}},
\quad \nu > -1, \quad \y \in [0, \infty)^N,
\label{eqn:P_Laguerre}
\end{align}
with the normalization constants 
$Z^{(N)}_{\rm Hermite}$ and
$Z^{(\nu, N)}_{\rm Laguerre}$. 

The DPP $(\Xi, K^{(N)}_{\rm Hermite}, \lambda_{\rN(0,1/2)}(dy))$
describes
the eigenvalue distribution of 
$N \times N$ Hermitian random 
matrices in the {\it Gaussian unitary ensemble} (GUE).
When $\nu \in \N_0$, 
the DPP 
$(\Xi, K^{(N)}_{\rm Laguerre}$, $\lambda_{\Gamma(\nu+1, 1)}(dx))$
describes the distribution of the nonnegative square roots
of eigenvalues of $M^{\dagger} M$,
where $M$ is $(N+\nu) \times N$  complex
random matrix in the 
{\it chiral Gaussian ensemble} (chGUE)
and $M^{\dagger}$ is its Hermitian conjugate.
The probability density (\ref{eqn:P_Laguerre}) 
can be extended to any $\nu \in (-1, \infty)$ and it
is called the {\it complex Laguerre ensemble}
or the {\it complex Wishart ensemble}.
Many other examples of one-dimensional DPPs 
are given as eigenvalue ensembles
of Hermitian random matrices in the literature of
random matrix theory (see, for instance, \cite{Meh04,For10,KT04}).

\subsubsection{Duality relations between 
DPPs in continuous and discrete spaces}
\label{sec:duality_application1}
We consider the simplified setting (\ref{eqn:W_n_x1}) of $W$
with $\Gamma=\N_0$.
If we set
$\Lambda_1=[r, \infty) \subset S_1=\R, r \in \R$ and
$\Lambda_2=\{0,1, \dots, N-1\} \subset S_2=\Gamma=\N_0, N \in \N$
in (\ref{eqn:K_A}), we obtain 
\begin{align}
K_{\R}^{\{0,1, \dots, N-1 \}}(x, x')
&= \sum_{n=0}^{N-1} \varphi_n(x) \overline{\varphi_n(x')},
\quad x, x' \in \R,
\nonumber\\
K^{[r, \infty)}_{\N_0}(n, n')
&= \int_r^{\infty} \overline{\varphi_n(x)} \varphi_{n'}(x) \lambda_1(dx),
\quad n, n' \in \N_0.
\label{eqn:Kr_inf}
\end{align}
When $\lambda_1(dx)$ and $\{\varphi_n\}_{n \in \N_0}$
are given by (\ref{eqn:HermiteZ1}) or by (\ref{eqn:LaguerreZ1}),
the kernels (\ref{eqn:Kr_inf}) are given by
\begin{align*}
& K_{{\rm DHermite}^+(r)}(n, n')
= ( \pi 2^{n+n'} n! n'!)^{-1/2}
\int_r^{\infty} H_n(x) H_{n'}(x) e^{-x^2} dx
\nonumber\\
& \qquad = - ( \pi n! n'! 2^{n+n'+2} )^{-1/2} 
e^{-r^2} 
\frac{H_{n+1}(r) H_{n'}(r) - H_n(r) H_{n'+1}(r)}{n-n'},
\end{align*}
and, provided $r > 0$, 
\begin{align*}
& K_{{\rm DLaguerre}^+(r, \nu+1)}(n, n')
= \left( \frac{n ! n' !}{\Gamma(n+\nu+1) \Gamma(n'+\nu+1)} \right)^{1/2}
\int_r^{\infty} L^{(\nu)}_n(x) L^{(\nu)}_{n'}(x) x^{\nu} e^{-x} dx
\nonumber\\
& \qquad = \left( \frac{n ! n' !}{\Gamma(n+\nu+1) \Gamma(n'+\nu+1)} \right)^{1/2}
r^{\nu+1} e^{-r}
\frac{L^{(\nu+1)}_{n-1}(r) L^{(\nu)}_{n'}(r)
- L^{(\nu)}_{n}(r) L^{(\nu+1)}_{n'-1}(r)}{n-n'},
\end{align*}
with the convention that $L^{(\nu)}_{-1}(r)=0$, 
respectively (see Propositions 3.3 and 3.4 in \cite{BO17}).
Borodin and Olshanski called
the correlation kernels $K_{{\rm DHermite}^+(r)}$ and
$K_{{\rm DLaguerre}^+(r, \nu+1)}$ the {\it discrete Hermite kernel}
and the {\it discrete Laguerre kernel}, respectively \cite{BO17}.
Theorem~\ref{thm:duality2} gives 
\begin{equation}
\bP(\Xi_1^{\{0,1, \dots, N-1\}}([r, \infty))=m)
=\bP(\Xi_2^{[r, \infty)}(\{0,1, \dots, N-1\})=m),
\quad \forall m \in \N_0,
\label{eqn:duality_eq3}
\end{equation}
where LHS denotes the probability that
the number of points in the interval $[r, \infty)$
is $m$ for the $N$-point continuous DPP on $\R$ such as
$(\Xi_1, K^{(N)}_{\rm Hermite}, \lambda_{N(0, 1/2)}(dx))$
or $(\Xi_1, K^{(N)}_{\rm Laguerre}, \lambda_{\Gamma(\nu+1, 1)}(dx))$, 
$\nu \in (-1, \infty)$, 
while RHS does the probability that
the number of points in $\{0,1,\dots, N-1\}$ is $m$ for
the discrete DPP on $\N_0$ such as
$(\Xi_2, K_{{\rm DHermite}^{+}(r)})$
or $(\Xi_2, K_{{\rm DLaguerre}^{+}(r, \nu+1)})$, $\nu \in (-1, \infty)$. 
The {\it duality between continuous and discrete ensembles}
of Borodin and Olshanski (Theorem 3.7 in \cite{BO17}) is
a special case with $m=0$ of the equality (\ref{eqn:duality_eq3}).

\subsection{Finite DPPs in intervals
related with classical root systems}
\label{sec:Lie_group}

Let $N \in \N$ and consider the four types of
{\it classical root systems} denoted by $\AN, \BN, \CN$, and $\DN$
(see Appendix \ref{sec:Weyl_denominators}).
We set
$S^{\AN}=\S^1=[0, 2 \pi)$, the unit circle, 
with a uniform measure 
$\lambda^{\AN}(dx)=\lambda_{[0, 2 \pi)}(dx) :=dx/(2 \pi)$, and
$S^{\RN}=[0, \pi]$, the upper half-circle, with 
$\lambda^{\RN}(dx)=\lambda_{[0, \pi]}(dx) :=dx/\pi$ for
$\RN=\BN, \CN, \DN$. 

For a fixed $N \in \N$, we introduce the four sets of
functions $\{\varphi^{\RN}_n\}_{n=1}^N$ on $S^{\RN}$ defined as
\begin{align*}
\varphi^{\RN}_n(x) := \begin{cases}
\displaystyle{
e^{-i(\cN^{\AN}-2 J^{\AN}(n))x/2}
}, 
& \RN=\AN,
\cr
\displaystyle{
\sin \big[(\cN^{\RN}-2J^{\RN}(n)) x/2 \big]
}, 
& \RN=\BN, \CN, 
\cr
\displaystyle{
\cos \big[(\cN^{\DN}-2J^{\DN}(n)) x/2 \big]
}, 
& \RN=\DN,
\end{cases}
\end{align*}
where
\begin{equation}
\cN^{\RN} := \begin{cases}
N, \quad & \RN = \AN, \\
2N-1, \quad & \RN = \BN, \\
2(N+1), \quad & \RN =\CN, \\
2(N-1), \quad & \RN = \DN.
\end{cases}
\label{eqn:N_R_classic}
\end{equation}
and
\begin{equation}
J^{\RN}(n) := \begin{cases}
n-1/2, \quad & \RN = \AN, 
\\
n-1, \quad & \RN=\BN, \DN,
\\
n, \quad & \RN=\CN.
\end{cases}
\label{eqn:J_R_classic}
\end{equation}
It is easy to verify that they satisfy the following
orthonormality relations,
\begin{align*}
\langle \varphi^{\AN}_n, \varphi^{\AN}_m \rangle_{L^2(\S^1, \lambda_{[0, 2 \pi)})}
&=\delta_{n m}, 
\nonumber\\
\langle \varphi^{\RN}_n, \varphi^{\RN}_m \rangle_{L^2([0, \pi], \lambda_{[0, \pi]})}
&=\delta_{n m}, 
\quad \RN=\BN, \CN, \DN, \quad
\mbox{if $n, m \in \{1, \dots, N\}$}.
\end{align*}
We put $\Gamma = \{1, \dots, N\}, N \in \N$ and
$L^2(\Gamma, \nu) = \ell^2(\Gamma)$.
By the argument given in Remark 3 in Section \ref{sec:OP}, 
Assumption 3 is verified, and hence
Corollary \ref{thm:main2} gives the four types of DPPs; 
$(\Xi, K^{\AN}, \lambda_{[0, 2 \pi)}(dx))$ on $\S^1$,
and 
$(\Xi, K^{\RN}, \lambda_{[0, \pi]}(dx))$ on $[0, \pi]$,
$\RN=\BN, \CN, \DN$, with the correlation kernels,
\begin{align*}
K^{\RN}(x, x^{\prime})
&= \sum_{n=1}^N \varphi^{\RN}_n(x) \overline{\varphi^{\RN}_n(x^{\prime})}
\nonumber\\
&= \begin{cases}
\displaystyle{
\frac{\sin\{N(x-x^{\prime})/2\}}{\sin\{(x-x^{\prime})/2\}}
},
& \RN=\AN,
\cr
\displaystyle{ \frac{1}{2} \left[
\frac{\sin\{N(x-x^{\prime})\}}{\sin\{(x-x^{\prime})/2\}}
- \frac{\sin\{N(x+x^{\prime})\}}{\sin\{(x+x^{\prime})/2\}}
\right]},
& \RN=\BN, 
\cr
& \cr
\displaystyle{\frac{1}{2} \left[
\frac{\sin\{(2N+1)(x-x^{\prime})/2\}}{\sin\{(x-x^{\prime})/2\}}
- \frac{\sin\{(2N+1)(x+x^{\prime})/2\}}{\sin\{(x+x^{\prime})/2\}}
\right]},
& \RN=\CN, 
\cr
& \cr
\displaystyle{\frac{1}{2} \left[
\frac{\sin\{(2N-1)(x-x^{\prime})/2\}}{\sin\{(x-x^{\prime})/2\}}
+ \frac{\sin\{(2N-1)(x+x^{\prime})/2\}}{\sin\{(x+x^{\prime})/2\}}
\right]},
& \RN=\DN.
\cr
\end{cases}
\end{align*}

By Lemma \ref{thm:Weyl_trigonometric}
in Appendix \ref{sec:Weyl_denominators}, 
the probability densities for these finite DPPs
with respect to 
the Lebesgue measures,
$d \x =\prod_{j=1}^N d x_j$ are given
as
\begin{align*}
\bp^{\AN}(\x) 
&= \frac{1}{Z^{\AN}} \prod_{1 \leq j < k \leq N}
\sin^2 \frac{x_k-x_j}{2}, 
\quad \x \in [0, 2\pi)^N,
\nonumber\\
\bp^{\BN}(\x) 
&= \frac{1}{Z^{\BN}} 
\prod_{1 \leq j < k \leq N} 
\left( \sin^2 \frac{x_k-x_j}{2} \sin^2 \frac{x_k+x_j}{2} \right)
\prod_{\ell=1}^N \sin^2 \frac{x_{\ell}}{2},
\quad \x \in [0, \pi]^N,
\nonumber\\
\bp^{\CN}(\x) 
&= \frac{1}{Z^{\CN}} 
\prod_{1 \leq j < k \leq N} 
\left( \sin^2 \frac{x_k-x_j}{2} \sin^2 \frac{x_k+x_j}{2} \right)
\prod_{\ell=1}^N \sin^2 x_{\ell},
\quad \x \in [0, \pi]^N,
\nonumber\\
\bp^{\DN}(\x) 
&= \frac{1}{Z^{\DN}} 
\prod_{1 \leq j < k \leq N} 
\left( \sin^2 \frac{x_k-x_j}{2} \sin^2 \frac{x_k+x_j}{2} \right),
\quad \x \in [0, \pi]^N,
\end{align*}
with the normalization constants $Z^{\RN}$.

The DPP, $(\Xi, K^{\AN}, \lambda_{[0, 2 \pi)}(dx))$ is known as
the {\it circular unitary ensemble} (CUE)
in random matrix theory (see Section 11.8 in \cite{Meh04}).
These four types of DPPs, 
$(\Xi, K^{\AN}, \lambda_{[0, 2 \pi)}(d x))$,
$(\Xi, K^{\RN}, \lambda_{[0, \pi]}(d x))$, 
$\RN=\BN, \CN, \DN$
are realized as the eigenvalue distributions of
random matrices in the {\it classical groups}, 
$\U(N)$, $\SO(2N+1)$, $\Sp(N)$, and $\SO(2N)$, respectively.
(See Section 2.3 c) in \cite{Sos00}
and Section 5.5 in \cite{For10}.)

\subsection{Finite DPPs in intervals 
related with affine root systems} 
\label{sec:affine}

We define the following four types of functions;
\begin{align}
\Theta^{A}(\sigma, z, \tau) 
&:= e^{2 \pi i \sigma z} 
\vartheta_2 (\sigma \tau + z; \tau),
\nonumber\\
\Theta^{B}(\sigma, z, \tau) 
&:= e^{2 \pi i \sigma z} \vartheta_1 (\sigma \tau + z; \tau) 
- e^{-2 \pi i \sigma z} \vartheta_1 (\sigma \tau - z ; \tau), 
\nonumber\\
\Theta^{C}(\sigma, z, \tau) 
&:= e^{2 \pi i \sigma z} \vartheta_2 (\sigma \tau + z; \tau)
- e^{-2 \pi i \sigma z} \vartheta_2 (\sigma \tau - z; \tau),
\nonumber\\
\Theta^{D}(\sigma, z, \tau)  
&:= e^{2 \pi i \sigma z} \vartheta_2 (\sigma \tau + z; \tau)
+ e^{-2 \pi i \sigma z} \vartheta_2 (\sigma \tau - z; \tau),
\label{eqn:RS}
\end{align}
for $\sigma \in \R, z \in \C, \tau \in \H := \{z \in \C: \Im z >0\}$,
where $\vartheta_{\mu}(v; \tau), \mu=1,2$ 
are the {\it Jacobi theta functions}.
See Appendix \ref{sec:Jacobi}
for definitions and the basic properties of
the Jacobi theta functions.

Here we consider the seven types
of {\it irreducible reduced affine root systems}
$\RN=\AN$, $\BN$, $\BNv$, $\CN$, 
$\CNv$, $\BCN$, $\DN$, $N \in \N$ \cite{Mac72,RS06}
(see Appendix \ref{sec:Macdonald_denominators}).
We put
$S^{\AN}=\S^1=[0, 2\pi)$ with 
$\lambda^{\AN}(dx)=\lambda_{[0, 2 \pi)}(dx)$, and
$S^{\RN}=[0, \pi]$ with $\lambda^{\RN}(dx)=\lambda_{[0, \pi]}(dx)$ 
for $\RN=\BN$, $\BNv$, $\CN$, 
$\CNv$, $\BCN$, $\DN$.
We assume that $\tau \in \H$ is pure imaginary, that is,
\[
\tau = i \Im \tau \in i (0, \infty).
\]
For a fixed $N \in \N$, we define the seven sets of functions
$\{\varphi^{\RN}_n(x; \tau)\}_{n=1}^N$ on $S^{\RN}$ as
\[
\varphi_n^{\RN}(x; \tau) := \frac{1}{\sqrt{m^{\RN}_n(\tau)}}
\Theta^{\sharp(\RN)}
\left( \frac{J^{\RN}(n)}{\cN^{\RN}}, 
\cN^{\RN} \frac{x}{2 \pi}, \tau \right), 
\]
where 
\begin{equation}
\sharp(\RN)
:= \begin{cases}
A, \quad & \mbox{if $\RN=\AN$},
\cr
B, \quad & \mbox{if $\RN=\BN, \BNv$},
\cr
C, \quad & \mbox{if $\RN=\CN, \CNv, \BCN$},
\cr
D, \quad & \mbox{if $\RN=\DN$},
\end{cases}
\label{eqn:sharp}
\end{equation}
\begin{equation}
\cN^{\RN} := \begin{cases}
N, \quad & \RN = \AN, \\
2N-1, \quad & \RN = \BN, \\
2N, \quad & \RN = \BNv, \CNv, \\
2(N+1), \quad & \RN =\CN, \\
2N+1, \quad & \RN = \BCN, \\
2(N-1), \quad & \RN = \DN,
\end{cases}
\label{eqn:N_R}
\end{equation}
\begin{equation}
J^{\RN}(n) := \begin{cases}
n-1/2, \quad & \RN = \AN, \CNv,
\\
n-1, \quad & \RN=\BN, \BNv, \DN,
\\
n, \quad & \RN=\CN, \BCN,
\end{cases}
\label{eqn:J_R}
\end{equation}
and we set
\begin{align*}
m^{\AN}_n(\tau) &:=
\vartheta_2(2 J^{\AN}(n) \tau/\cN^{\AN}; 2 \tau),
\quad n \in \{1, \dots, N \},
\nonumber\\
m^{\RN}_n(\tau) &:=
2 \vartheta_2(2 J^{\RN}(n) \tau/\cN^{\RN} ; 2 \tau),
\quad n \in \{1, \dots, N \},
\quad \mbox{for $\RN=\CN, \CNv, \BCN$},
\nonumber\\
m^{\RN}_n(\tau) &:= \begin{cases}
4 \vartheta_2 (0; 2 \tau),
& n=1, \\
2 \vartheta_2 (2 J^{\RN}(j) \tau/\cN^{\RN} ; 2 \tau),
& n \in \{2, 3, \dots, N\},
\end{cases}
\quad \mbox{for $\RN=\BN, \BNv$},
\nonumber\\
m^{\DN}_n(\tau) &:= \begin{cases}
4 \vartheta_2 (0; 2 \tau),
& n=1,
\\
2 \vartheta_2 (2 J^{\DN}(j) \tau/\cN^{\DN}; 2 \tau),
& n \in \{2, 3, \dots, N-1\},
\\
4 \vartheta_2 (2 (N-1) \tau/\cN^{\DN}; 2 \tau),
& n=N.
\end{cases}
\end{align*}

For $N \in \N$, 
the following orthonormality relations can be proved 
as a special case of Lemma 2.1 in \cite{Kat19a},
\begin{align*}
\langle \varphi^{\AN}_n(\cdot, \tau), 
\varphi^{\AN}_m(\cdot, \tau) \rangle_{L^2(\S^1, \lambda_{[0, 2 \pi)})}
&=\delta_{n m},
\nonumber\\
\langle \varphi^{\RN}_n(\cdot, \tau), 
\varphi^{\RN}_m(\cdot, \tau) \rangle_{L^2([0, \pi], \lambda_{[0, \pi]})}
&=\delta_{n m}, \quad
\RN=\BN, \BNv, \CN, \CNv, \BCN, \DN, 
\nonumber\\
& \qquad
n, m \in \Gamma :=\{1, \dots, N\}.
\end{align*}
By the argument given in Remark 3 in Section \ref{sec:OP}, 
Assumption 3 is verified, and hence
Corollary \ref{thm:main2} gives
the seven types of DPPs,
$(\Xi, K^{\RN}_{\tau}, \lambda^{\RN}(dx))$
with the correlation kernels,
\[
K^{\RN}_{\tau}(x, x^{\prime})
= \sum_{n=1}^N \varphi^{\RN}_n(x; \tau) \overline{\varphi^{\RN}_n(x'; \tau)},
\quad \RN=\AN, \BN, \BNv, \CN, \CNv, \BCN, \RN.
\]
Thanks to the {\it Macdonald denominator formula}
proved by Rosengren and Schlosser \cite{RS06}
(see (3.1) with (3.2) in \cite{Kat19a} in the present notations), 
the probability densities for these finite DPPs
with respect to 
the Lebesgue measure,
$d \x =\prod_{j=1}^N d x_j$ are given
as follows, 
\begin{align*}
\bp^{\AN}_{\tau}(\x) 
&= \begin{cases}
\displaystyle{
\frac{1}{Z^{\AN}(\tau)}
\left|
\vartheta_0 \left( \sum_{j=1}^N \frac{x_j}{2\pi}; \frac{\tau}{\cN^{\AN}} \right)
W^{\AN}\left( \frac{\x}{2 \pi}; \frac{\tau}{\cN^{\AN}} \right) 
\right|^2
},
& \mbox{if $N$ is even}, \cr
\displaystyle{
\frac{1}{Z^{\AN}(\tau)}
\left|
\vartheta_3 \left( \sum_{j=1}^N \frac{x_j}{2\pi}; \frac{\tau}{\cN^{\AN}} \right)
W^{\AN}\left( \frac{\x}{2 \pi}; \frac{\tau}{\cN^{\AN}} \right) 
\right|^2
},
& \mbox{if $N$ is odd}, \cr
& \x \in [0, 2 \pi)^N, 
\end{cases}
\nonumber\\
\bp^{\RN}_{\tau}(\x)
&=\frac{1}{Z^{\RN}(\tau)} \left|
W^{\RN} \left( \frac{\x}{2 \pi}; \frac{\tau}{\cN^{\RN}} \right) 
\right|^2,
\,
\x \in [0, \pi]^N, \, \RN=\BN, \BNv, \CN, \CNv, \BCN, \DN, 
\end{align*}
where $W^{\RN}$ are the Macdonald denominators
given by (\ref{eqn:Macdonald_denominators})
in Appendix \ref{sec:Macdonald_denominators}
and $Z^{\RN}(\tau)$ are the normalization constants.
By the properties (\ref{eqn:quasi_periodic1}) and (\ref{eqn:theta_values})
of the Jacobi theta functions, 
it is easy to verify the following,
\begin{align*}
& \cS_{2\pi/N} (\Xi, K_{\tau}^{\AN}, \lambda_{[0, 2 \pi)})
\law (\Xi, K_{\tau}^{\AN}, \lambda_{[0, 2 \pi)}),
\nonumber\\
& \rho_{\tau}^{\RN}(0)=0, \quad
\RN=\BN, \CNv, \BCN,
\nonumber\\
& \rho_{\tau}^{\RN}(0)=\rho_{\tau}^{\RN}(\pi)=0,
\quad \RN=\BNv, \CN.
\end{align*}

In \cite{Kat19a}, it was proved that 
these seven types of DPPs
are realized as the particle configurations 
at the middle time $t=t_{\ast}/2$ 
of the 
{\it noncolliding Brownian bridges}
in time duration $[0, t_{\ast}]$, 
provided $t_{\ast}=4 \pi \Im \tau > 0$, whose
initial configurations at $t=0$ and 
final configurations at $t=t_{\ast}$ are fixed 
to be specially chosen configurations 
depending on the types $\RN, N \in \N$.

As $\Im \tau \to \infty$, the temporal inhomogeneity 
in such noncolliding Brownian bridges vanishes.
Associated with such limit transitions, 
the following degeneracies are observed 
in the weak convergence of DPPs 
from the seven types of affine root systems
to the four types of classical root systems, 
\begin{align}
(\Xi, K_{\tau}^{\AN}, \lambda_{[0, 2 \pi)}(d x))
&\weakT (\Xi, K^{\AN}, \lambda_{[0, 2 \pi)}(d x)), 
\nonumber\\
\left.
\begin{array}{l}
\displaystyle{(\Xi, K_{\tau}^{\BN}, \lambda_{[0, \pi]}(d x))} \cr
\displaystyle{(\Xi, K_{\tau}^{\CNv}, \lambda_{[0, \pi]}(d x))} \cr
\displaystyle{(\Xi, K_{\tau}^{\BCN}, \lambda_{[0, \pi]}(d x))}
\end{array}
\right\}
&\weakT (\Xi, K^{\BN}, \lambda_{[0, \pi]}(d x)), 
\nonumber\\
\left.
\begin{array}{l}
\displaystyle{(\Xi, K_{\tau}^{\BNv}, \lambda_{[0, \pi]}(d x))} \cr
\displaystyle{(\Xi, K_{\tau}^{\CN}, \lambda_{[0, \pi]}(d x))}
\end{array}
\right\}
&\weakT (\Xi, K^{\CN}, \lambda_{[0, \pi]}(d x)),
\nonumber\\
(\Xi, K_{\tau}^{\DN}, \lambda_{[0, \pi]}(d x))
&\weakT (\Xi, K^{\DN}, \lambda_{[0, \pi]}(d x)), 
\label{eqn:T_inf_1}
\end{align}
where the DPPs, $(\Xi, K^{\AN}, \lambda_{[0, 2\pi)})$
and $(\Xi, K^{\RN}, \lambda_{[0, \pi]})$, 
$\RN=\BN, \CN, \DN$
were given in Section \ref{sec:Lie_group}.

\subsection{Infinite DPPs in $\R$ 
associated with classical orthonormal functions} 
\label{sec:classical_functions}

Here we give examples of infinite DPPs obtained
by Corollary \ref{thm:main3}.

\begin{description}
\item{(i)} \quad DPP with the {\it sinc kernel} : 
We set $S_1=\R$, $\lambda_1(dx)=dx$, 
$\Gamma=(-1, 1)$, $\nu(dy)=\lambda_2(d y)=dy$, and put
\[
\psi_1(x, y)=\frac{1}{\sqrt{2 \pi}} e^{ i x y} .
\]
We see that 
$\Psi_1(x)^2 \equiv 1/\pi, x \in \R$ and thus Assumption 3'(ii) is satisfied.
The correlation kernel $K_{S_1}$ is given by
\[
K_{\rm sinc}(x, x^{\prime})
= \frac{1}{2 \pi}\int_{-1}^1 e^{i y(x-x^{\prime})} dy
= \frac{\sin(x-x^{\prime}) }{ \pi (x-x^{\prime}) },
\quad x, x^{\prime} \in \R.
\]

\item{(ii)} \quad DPP with the {\it Airy kernel} :
We set $S_1=\R$, $\lambda_1(dx)=dx$, 
$\Gamma=[0, \infty)$, $\nu(dy)=\lambda_2(d y)=dy$, 
and put
\[
\psi_1(x,y)=\Ai(x+y),
\]
where $\Ai(x)$ denotes the {\it Airy function} \cite{NIST10}
\[
\Ai(x) := 
\frac{1}{\pi} \int_0^{\infty} \cos \left( \frac{k^3}{3}+kx \right) dk.
\]
We see that 
$\Psi_1(x)^2=-x \Ai(x)^2 + \Ai'(x)^2, x \in \R$
and thus Assumption 3(ii) is satisfied. 
The correlation kernel $K_{S_1}$ is given by
\[
K_{\rm Airy}(x, x^{\prime})
= \int_0^{\infty} \Ai(x+y) \Ai(x^{\prime}+y) dy
= \frac{\Ai(x) \Ai^{\prime}(x^{\prime}) - \Ai(x^{\prime}) \Ai^{\prime}(x)}
{x-x^{\prime}},
\quad x, x^{\prime} \in \R,
\]
where $\Ai^{\prime}(x) := d \Ai(x)/dx$.

\item{(iii)} \quad DPP with the {\it Bessel kernel} :
We set
$S_1=[0,\infty)$, $\lambda_1(dx)= dx$,
$\Gamma=[0, 1]$,
$\nu(dy)=\lambda_2(d y)=dy$.
With parameter $\nu > -1$ we put
\[
\psi_1(x,y) =\sqrt{x y} J_{\nu}(x y),
\]
where $J_{\nu}$ is the 
{\it Bessel function of the first kind} defined by
\begin{equation}
J_{\nu}(x) := \sum_{n=0}^{\infty} 
\frac{(-1)^n}{n! \Gamma(\nu+n+1)}
\left(\frac{x}{2} \right)^{2n+\nu},
\quad x \in \C \setminus (-\infty, 0).
\label{eqn:Bessel_function}
\end{equation}
We see that 
$\Psi_1(x)^2= x \{J_{\nu}(x)^2 - J_{\nu-1}(x)J_{\nu+1}(x)\}/2$,
$x \in [0, \infty)$ 
and thus Assumption 3'(ii) is satisfied.
The correlation kernel $K_{S_1}$ is given by
\begin{align}
K_{\rm Bessel}^{(\nu)}(x, x')
&= \int_0^1 \sqrt{x y} J_{\nu}(x y) \sqrt{x' y} J_{\nu}(x' y) dy 
\nonumber\\
&= \frac{\sqrt{x x'}}{x^2-(x')^2}
 \Big\{ J_{\nu}(x) x' J_{\nu}^{\prime}(x')
 - x J_{\nu}^{\prime}(x) J_{\nu}(x') \Big\},
 \quad x, x' \in [0, \infty),
 \label{eqn:Bessel_kernel}
\end{align}
where $J_{\nu}^{\prime}(x) := d J_{\nu}(x)/dx$.
\end{description}

These three kinds of infinite DPPs,
$(\Xi, K_{\rm sinc}, d x)$, $(\Xi, K_{\rm Airy}, d x)$,
and $(\Xi, K_{\rm Bessel}^{(\nu)}, {\bf 1}_{[0, \infty)}(x) d x)$,
are obtained as the scaling limits
of the finite DPPs,
$(\Xi, K^{(N)}_{\rm Hermite}, \lambda_{\rN(0, 1/2)}(d x))$ 
and $(\Xi, K^{(\nu, N)}_{\rm Laguerre}$, $\lambda_{\Gamma(\nu+1, 1)}(d x))$, 
given in Section \ref{sec:classical}
as follows.

\begin{description}
\item{(i)} {\it Bulk scaling limit},
\[
\sqrt{2N} \circ
(\Xi, K^{(N)}_{\rm Hermite}, \lambda_{\rN(0, 1/2)}(d x))
\weak (\Xi, K_{\rm sinc}, d x). 
\]

\item{(ii)} {\it Soft-edge scaling limit},
\[
\sqrt{2}N^{1/6} \circ \cS_{\sqrt{2N}} \,
(\Xi,  K^{(N)}_{\rm Hermite}, \lambda_{\rN(0, 1/2)}(d x))
\weak (\Xi, K_{\rm Airy}, d x). 
\]

\item{(iii)} {\it Hard-edge scaling limit},
for $\nu > -1$, 
\[
4N \circ \Big(
(\Xi, K^{(\nu, N)}_{\rm Laguerre}, \lambda_{\Gamma(\nu+1,1)}(d x))^{\langle 1/2 \rangle}
\Big)
\weak (\Xi, K^{(\nu)}_{\rm Bessel}, {\bf 1}_{[0, \infty)}(x) d x).
\]
\end{description}

\noindent
See, for instance, \cite{Meh04,For10,AGZ10,Kat15_Springer}, 
for more details.

The DPPs with the sinc kernel and 
the Bessel kernel with the special values of
parameter $\nu$ can be obtained as the bulk
scaling limits of the DPPs,
$(\Xi, K^{\RN}, \lambda^{\RN}(d x))$,
$\RN=\AN, \BN, \CN, \DN$
given in Section \ref{sec:Lie_group} as
\begin{align}
\frac{N}{2} \circ (\Xi, K^{\AN}, \lambda_{[0, 2 \pi)}(d x))
&\weak (\Xi, K_{\rm sinc}, d x),
\nonumber\\
\left.
\begin{array}{l}
N \circ (\Xi, K^{\BN}, \lambda_{[0, \pi]}(d x))
\cr
N \circ (\Xi, K^{\CN}, \lambda_{[0, \pi]}(d x))
\end{array}
\right\}
&\weak (\Xi, K_{\rm Bessel}^{(1/2)}, {\bf 1}_{[0, \infty)}(x) d x),
\nonumber\\
N \circ (\Xi, K^{\DN}, \lambda_{[0, \pi]}(d x))
&\weak (\Xi, K_{\rm Bessel}^{(-1/2)}, {\bf 1}_{[0, \infty)}(x) d x),
\label{eqn:Lie_bulk}
\end{align}
where 
\begin{align*}
K_{\rm Bessel}^{(1/2)}(x, x^{\prime})
&= \frac{\sin (x-x^{\prime})}{\pi (x-x^{\prime})} 
- \frac{\sin(x+x^{\prime})}{\pi (x+x^{\prime})}, 
\quad x, x^{\prime} \in [0, \infty), 
\nonumber\\
K_{\rm Bessel}^{(-1/2)}(x, x^{\prime})
&= \frac{\sin(x-x^{\prime})}{\pi (x-x^{\prime})} 
+ \frac{\sin(x+x^{\prime})}{\pi (x+x^{\prime})}, 
\quad x, x^{\prime} \in [0, \infty).
\end{align*}
Since $J_{1/2}(x)=\sqrt{2/(\pi x)} \sin x$ and
$J_{-1/2}(x)=\sqrt{2/(\pi x)} \cos x$, 
the above correlation kernels 
are readily obtained from (\ref{eqn:Bessel_kernel})
by setting $\nu=1/2$ and $-1/2$, respectively.

\subsection{Infinite DPPs in $\R$ 
associated with orthonormal theta functions} 
\label{sec:orthogonal_theta}

Let $S =\R$ with
$\lambda^{A}(dx)=dx$,
and
$\lambda^{R}(dx)={\bf 1}_{[0, \infty)}(x) dx$ 
for the types $R=B, C, D$. 
Here we assume that $\tau \in \H$ is pure imaginary. 
We put
\begin{align*}
\psi^A(x, \gamma; \tau) 
&:= \frac{\Theta^{A}( \gamma, x/\pi, \tau)}{\sqrt{\pi \vartheta_2(2 \tau \gamma; 2 \tau)}},
\nonumber\\
\psi^R(x, \gamma; \tau) 
&:= \frac{\Theta^{R}(\gamma/2, x/\pi, \tau)}{\sqrt{2 \pi \vartheta_2(\tau \gamma; 2 \tau)}},
\quad R=B, C, D,
\end{align*}
and $\Gamma =(0,1)$ with
$\nu(d \gamma)=d \gamma$, 
where $\Theta^{R}, R=A, B, C, D$ are given by (\ref{eqn:RS}).

Using the equalities
\begin{align*}
& \frac{1}{\pi} 
\int_{\R} e^{ 2 i \{(\gamma-\gamma')+(n-m)\}x} dx
= \delta_{nm} \delta(\gamma-\gamma'),
\nonumber\\
& \frac{1}{2\pi} \int_{\R} e^{i \{(\gamma+\gamma')+2(n+m-1)\}x} dx
= 0,
\quad \mbox{for $n, m \in \Z, \quad \gamma, \gamma' \in \Gamma$},
\end{align*}
we can show the orthonormality relations \cite{Kat19a};
\begin{align*}
\langle \psi^A(\cdot, \gamma; \tau),
\psi^A(\cdot, \gamma'; \tau) \rangle_{L^2(\R, d x)}
&=\delta(\gamma-\gamma'),
\nonumber\\
\langle \psi^R(\cdot, \gamma; \tau),
\psi^R(\cdot, \gamma'; \tau) \rangle_{L^2(\R, {\bf 1}_{[0, \infty)}(x) d x)}
&=\delta(\gamma-\gamma'), \quad R=A, B, C, D,
\quad \gamma, \gamma' \in \Gamma.
\end{align*}
We can also evaluate the upper bounds for
$\Psi^R(x; \tau)^2 
:= \|\psi^R(x, \cdot; \tau) \|_{L^2(\Gamma, \nu)}^2$, $x \in \R$, 
$R=A, B, C, D$,
and confirm that Assumption 3(ii) is also satisfied.
Hence by Corollary \ref{thm:main2}, 
we obtain the
four types of infinite DPPs 
$(\Xi, K^A_{\tau}, d x)$,
$(\Xi, K^R_{\tau}, {\bf 1}_{[0, \infty)}(x) d x)$,
$R=B,C,D$.
The correlation kernels are written as follows,
\begin{align}
K^A_{\tau}(x, x^{\prime})
&= \frac{1}{\pi} \int_0^{1} 
e^{2 i  (x-x^{\prime}) \gamma} 
\frac{\vartheta_2 (x/\pi+\tau \gamma; \tau)
\vartheta_2 (x^{\prime}/\pi-\tau \gamma; \tau)}
{\vartheta_2 (2 \tau \gamma; 2 \tau)} d \gamma,
\quad x, x^{\prime} \in \R,
\nonumber\\
K^B_{\tau}(x, x^{\prime})
&= \frac{1}{2 \pi} \left[
\int_{-1}^{1} 
e^{i (x-x^{\prime}) \gamma} 
\frac{\vartheta_1 (x/\pi + \tau \gamma/2; \tau)
\vartheta_1 (x^{\prime}/\pi-\tau \gamma/2; \tau)}
{\vartheta_2 (\tau \gamma; 2 \tau)} d\gamma
\right.
\nonumber\\
& \qquad \left. 
+ \int_{-1}^{1} 
e^{i (x+x^{\prime}) \gamma} 
\frac{\vartheta_1 (x/\pi+\tau \gamma/2; \tau)
\vartheta_1 (x^{\prime}/\pi+\tau \gamma/2) ; \tau)}
{\vartheta_2 (\tau \gamma; 2 \tau)} d\gamma \right],
\quad x, x^{\prime} \in [0, \infty), 
\nonumber\\
K^C_{\tau}(x, x^{\prime})
&= \frac{1}{2 \pi} \left[
\int_{-1}^{1} 
e^{i (x-x^{\prime}) \gamma} 
\frac{\vartheta_2 (x/\pi + \tau \gamma/2; \tau)
\vartheta_2 (x^{\prime}/\pi-\tau \gamma/2; \tau)}
{\vartheta_2 (\tau \gamma; 2 \tau)} d\gamma
\right.
\nonumber\\
& \qquad \left. 
- \int_{-1}^{1} 
e^{i (x+x^{\prime}) \gamma} 
\frac{\vartheta_2 (x/\pi+ \tau \gamma/2; \tau)
\vartheta_2 (x^{\prime}/\pi+ \tau \gamma/2) ; \tau)}
{\vartheta_2 (\tau \gamma; 2 \tau)} d\gamma \right],
\quad x, x^{\prime} \in [0, \infty), 
\nonumber\\
K^D_{\tau}(x, x^{\prime})
&= \frac{1}{2 \pi} \left[
\int_{-1}^{1} 
e^{i (x-x^{\prime}) \gamma} 
\frac{\vartheta_2 (x/\pi + \tau \gamma/2; \tau)
\vartheta_2 (x^{\prime}/\pi- \tau \gamma/2; \tau)}
{\vartheta_2 (\tau \gamma; 2 \tau)} d\gamma
\right.
\nonumber\\
& \qquad \left. 
+ \int_{-1}^{1} 
e^{i (x+x^{\prime}) \gamma} 
\frac{\vartheta_2 (x/\pi+ \tau \gamma/2; \tau)
\vartheta_2 (x^{\prime}/\pi+ \tau \gamma/2) ; \tau)}
{\vartheta_2 (\tau \gamma; 2 \tau)} d\gamma \right],
\quad x, x^{\prime} \in [0, \infty).
\label{eqn:inf_affine_K}
\end{align}
If we change the integral variables appropriately,
the above become the correlation kernels
$\cK_{t_{\ast}/2}^{R}$, $R=A, B, C, D$, 
with $t_{\ast}=4 \pi \Im \tau$, given in
Lemma 3.5 in \cite{Kat19a}.
Using the quasi-periodicity (\ref{eqn:quasi_periodic1})
of the theta functions, we can show that
the DPP $(\Xi, K^A_{\tau}, dx)$ has a periodicity of $\pi$; 
$\cS_{\pi} K^A_{\tau}(x, x')=K^A_{\tau}(x, x')$, $x, x' \in \R$.
By the symmetry (\ref{eqn:even_odd}) of the theta functions, 
we see that $\rho^{R}_{\tau}(0)=K^{R}_{\tau}(0,0)=0$,
$R=B, C$.

The infinite DPPs associated with the
above correlation kernels (\ref{eqn:inf_affine_K}) are
obtained as the bulk scaling limits
of the finite DPPs in intervals 
studied in Section \ref{sec:affine} \cite{Kat19a}; 
\begin{align}
\frac{N}{2} \circ (\Xi, K_{\tau}^{\AN}, \lambda_{[0, 2 \pi)}(d x))
& \weak (\Xi, K_{\tau}^A, d x), 
\nonumber\\
\left.
\begin{array}{l}
\displaystyle{
N \circ (\Xi, K_{\tau}^{\BN}, \lambda_{[0, \pi]}(d x))}
\cr
\displaystyle{
N \circ (\Xi, K_{\tau}^{\BNv}, \lambda_{[0, \pi]}(d x))}
\end{array}
\right\}
& \weak (\Xi, K_{\tau}^B, {\bf 1}_{[0, \infty)}(x) d x), 
\nonumber\\
\left.
\begin{array}{l}
\displaystyle{
N \circ (\Xi, K_{\tau}^{\CN}, \lambda_{[0, \pi]}(d x))}
\cr
\displaystyle{
N \circ (\Xi, K_{\tau}^{\CNv}, \lambda_{[0, \pi]}(d x))}
\cr
\displaystyle{
N \circ (\Xi, K_{\tau}^{\BCN}, \lambda_{[0, \pi]}(d x))}
\end{array}
\right\}
& \weak (\Xi, K_{\tau}^C, {\bf 1}_{[0, \infty)}(x) d x),
\nonumber\\
N \circ (\Xi, K_{\tau}^{\DN}, \lambda_{[0, \pi]}(d x))
& \weak (\Xi, K_{\tau}^D, {\bf 1}_{[0, \infty)}(x) d x).
\label{eqn:affine_bulk}
\end{align}

If we take the limit $\Im \tau \to \infty$ in (\ref{eqn:inf_affine_K}),
we obtain the following three infinite DPPs,
\begin{align}
(\Xi, K_{\tau}^A, d x)
&\weakT (\Xi, K_{\rm sinc}, d x),
\nonumber\\
\left.
\begin{array}{l}
(\Xi, K_{\tau}^{B}, {\bf 1}_{[0, \infty)}(x) d x)
\cr
(\Xi, K_{\tau}^{C}, {\bf 1}_{[0, \infty)}(x)d x)
\end{array}
\right\}
&\weakT (\Xi, K_{\rm Bessel}^{(1/2)}, {\bf 1}_{[0, \infty)}(x) d x),
\nonumber\\
(\Xi, K_{\tau}^{D}, {\bf 1}_{[0, \infty)}(x)d x)
&\weakT (\Xi, K_{\rm Bessel}^{(-1/2)}, {\bf 1}_{[0, \infty)}(x) d x),
\label{eqn:T_inf_2}
\end{align}
which are the same as the limiting DPPs
given by (\ref{eqn:Lie_bulk}).

\vskip 0.3cm
\noindent{\bf Remark 4} \, 
The results (\ref{eqn:T_inf_1}), (\ref{eqn:Lie_bulk}),
(\ref{eqn:affine_bulk}), and (\ref{eqn:T_inf_2})
imply that,
in the limit transitions from the finite DPPs $(\Xi, K_{\tau}^{\RN}, \lambda^{\RN})$,
$\RN=\AN$, $\BN$, $\BNv$, $\CN$, $\CNv$, $\BCN$, $\DN$, 
to the infinite DPPs $(\Xi, K_{\rm sinc}, d x)$, 
$(\Xi, K^{(1/2)}_{\rm Bessel}, {\bf 1}_{[0, \infty)}(x) d x)$,
$(\Xi, K^{(-1/2)}_{\rm Bessel}, {\bf 1}_{[0, \infty)}(x) d x)$,
the scaling limits associated with $N \to \infty$
and the limit $\Im \tau \to \infty$ are commutable.

\SSC
{Examples in Two-dimensional Spaces} \label{sec:examples_2d}

\subsection{Infinite DPPs on $\C$ : Ginibre and Ginibre-type DPPs}
\label{sec:Ginibre_ACD}
\subsubsection{Three types of Ginibre DPPs}
\label{sec:three_Ginibre}
Let 
$\lambda_{\rN(m, \sigma^2; \C)}(d x)$ denote
the {\it complex normal distribution}, 
\begin{align*}
\lambda_{\rN(m, \sigma; \C)}(d x)
&:= \frac{1}{\pi \sigma^2} e^{-|x-m|^2/\sigma^2} d x
\nonumber\\
&= \frac{1}{\pi \sigma^2} 
e^{-(x_{\rR}-m_{\rR})^2/\sigma^2 -(x_{\rI}-m_{\rI})^2/\sigma^2} 
d x_{\rR}  d x_{\rI}, 
\end{align*}
$m \in \C, m_{\rR} :=\Re m, m_{\rI} :=\Im m, \sigma>0$.
We set
$S=\C$,
\begin{align*}
\lambda(dx) &=\lambda_{\rN(0,1;\C)}(dx)
= \frac{1}{\pi} e^{-|x|^2} dx
\nonumber\\
&=\lambda_{\rN(0,1/2)}(dx_{\rR}) \lambda_{\rN(0,1/2)}(dx_{\rI}),
\end{align*}
and 
\begin{align*}
\psi^{A}(x, \gamma) & :=
e^{-(x_{\rR}^2-x_{\rI}^2)/2+2 x \gamma}, 
\nonumber\\
\psi^{C}(x, \gamma) &:=
\sqrt{2} \sinh(2 x \gamma)
e^{-(x_{\rR}^2-x_{\rI}^2)/2},
\nonumber\\
\psi^{D}(x, \gamma) &:= 
\sqrt{2} \cosh(2 x \gamma)
e^{-(x_{\rR}^2-x_{\rI}^2)/2}.
\end{align*}
It is easy to confirm that
\begin{align*}
\frac{1}{\pi} 
\int_{\R} \psi^A(x, \gamma) \overline{\psi^A(x, \gamma')} e^{-x_{\rI}^2} d x_{\rI}
&=e^{-(x_{\rR}^2-4 x_{\rR} \gamma)} \delta(\gamma-\gamma'),
\cr
\frac{1}{\pi}
\int_{\R} \psi^R(x, \gamma) \overline{\psi^R(x, \gamma')} e^{- x_{\rI}^2} d x_{\rI}
&= e^{- x_{\rR}^2} 
\cosh(4 x_{\rR} \gamma)
\times 
\begin{cases}
\delta(\gamma-\gamma')-\delta(\gamma+\gamma'), & R=C, \cr
\delta(\gamma-\gamma')+\delta(\gamma+\gamma'), & R=D.
\end{cases}
\end{align*}
Therefore, we have
\begin{align*}
& \langle \psi^A(\cdot, \gamma), \psi^A(\cdot, \gamma') \rangle_{L^2(\C, \lambda_{\rN(0, 1; \C)})}
\nu(d \gamma)
=\delta(\gamma-\gamma') d \gamma, \quad
\gamma, \gamma' \in \Gamma^A :=\R,
\cr
& \langle \psi^R(\cdot, \gamma), \psi^R(\cdot, \gamma') \rangle_{L^2(\C, \lambda_{\rN(0,1; \C)})}
\nu(d \gamma)
=\delta(\gamma-\gamma') d \gamma, \quad
\gamma, \gamma' \in \Gamma^R :=(0, \infty), \quad R=C, D,
\end{align*}
with 
$\nu(d \gamma)=\lambda_{\rN(0, 1/4)}(d \gamma)
=\sqrt{2/\pi} e^{-2 \gamma^2} d \gamma$. 
We also see that
$\Psi^A(x)^2 :=\| \psi^A(x, \cdot) \|_{L^2(\Gamma^A, \nu)}^2
=e^{|x|^2}$,
$\Psi^C(x)^2 :=\| \psi^C(x, \cdot) \|_{L^2(\Gamma^C, \nu)}^2
=\sinh |x|^2$,
and 
$\Psi^D(x)^2 :=\| \psi^D(x, \cdot) \|_{L^2(\Gamma^D, \nu)}^2
=\cosh |x|^2$,
$x \in \C$.
Thus Assumption 3 is satisfied and we can
apply Corollary \ref{thm:main2}. 
The kernels (\ref{eqn:K_main2}) of obtained DPPs are given as
\begin{align*}
K^A(x, x^{\prime})
&= \sqrt{\frac{2}{\pi}} e^{-\{(x_{\rR}^2-x_{\rI}^2)+({x'_{\rR}}^2-{x'_{\rI}}^2) \}/2}
\int_{-\infty}^{\infty} 
e^{- 2 \{\gamma^2-(x+\overline{x'}) \gamma\}} d\gamma,
\nonumber\\
K^C(x, x^{\prime})
&= 2 \sqrt{\frac{2}{\pi}} e^{-\{(x_{\rR}^2-x_{\rI}^2)+({x'_{\rR}}^2-{x'_{\rI}}^2) \}/2}
\int_0^{\infty} 
e^{-2 \gamma^2}
\sinh(2 x \gamma) \sinh(2 \overline{x^{\prime}} \gamma) d\gamma, 
\nonumber\\
K^D(x, x^{\prime})
&= 2 \sqrt{\frac{2}{\pi}} e^{-\{(x_{\rR}^2-x_{\rI}^2)+({x'_{\rR}}^2-{x'_{\rI}}^2) \}/2}
\int_0^{\infty} 
e^{-2 \gamma^2}
\cosh(2 x \gamma) \cosh(2 \overline{x^{\prime}} \gamma) d\gamma. 
\end{align*}
The integrals are performed and we obtain 
$
K^R(x, x^{\prime})
= e^{i x_{\rR} x_{\rI}} K^{R}_{{\rm Ginibre}}(x, x^{\prime})
e^{-i x_{\rR}^{\prime} x_{\rI}^{\prime}}
$,
$R=A, C, D$, 
with
\begin{align}
K^{A}_{{\rm Ginibre}}(x, x^{\prime})
&= e^{x \overline{x^{\prime}}},
\label{eqn:K_Ginibre_A}
\\
K^{C}_{{\rm Ginibre}}(x, x^{\prime})
&= \sinh(x \overline{x^{\prime}}),  
\label{eqn:K_Ginibre_C}
\\
K^{D}_{{\rm Ginibre}}(x, x^{\prime})
&= \cosh(x \overline{x^{\prime}}), \quad x, x^{\prime} \in \C. 
\label{eqn:K_Ginibre_D}
\end{align}
Due to the gauge invariance of DPP mentioned in Section \ref{sec:basic_DPP}, 
the obtained three types of infinite DPPs on $\C$ are written as
$(\Xi, K_{\rm Ginibre}^R, \lambda_{\rN(0,1; \C)}(dx))$, $R=A, C, D$.
The DPP, $(\Xi$, $K_{\rm Ginibre}^A$, $\lambda_{\rN(0,1; \C)}(dx))$ 
with (\ref{eqn:K_Ginibre_A})
describes the eigenvalue distribution of the
Gaussian random complex matrix in the bulk scaling limit,
which is called the {\it complex Ginibre ensemble}
\cite{Gin65,Meh04,HKPV06,HKPV09,For10,Shirai15}.
This density profile is uniform 
with the Lebesgue measure $dx$ on $\C$ as
\[
\rho_{{\rm Ginibre}}(x) dx
=K^{A}_{{\rm Ginibre}}(x, x) \lambda_{\rN(0, 1; \C)}(d x)
=\frac{1}{\pi} dx_{\rR} dx_{\rI},
\quad x \in \C.
\]
On the other hands, the Ginibre DPPs
of types $C$ and $D$ 
with the correlation kernels (\ref{eqn:K_Ginibre_C})
and (\ref{eqn:K_Ginibre_D}) 
are rotationally symmetric around the origin, 
but non-uniform on $\C$.
The density profiles 
with the Lebesgue measure $dx$ on $\C$ are given by
\begin{align*}
\rho^{C}_{{\rm Ginibre}}(x) dx
&=
K^{C}_{{\rm Ginibre}}(x, x) \lambda_{\rN(0, 1; \C)}(d x)
=\frac{1}{2 \pi} (1-e^{-2 |x|^2}) dx_{\rR} dx_{\rI}, \quad x \in \C, 
\nonumber\\
\rho^{D}_{{\rm Ginibre}}(x) dx
&=K^{D}_{{\rm Ginibre}}(x, x) \lambda_{\rN(0, 1; \C)}(d x)
=\frac{1}{2 \pi} (1+e^{-2 |x|^2})dx_{\rR} dx_{\rI}, \quad x \in \C.
\end{align*}
They were first obtained in \cite{Kat19b} by taking
the limit $W \to \infty$ keeping the density of points of the infinite DPPs 
in the strip on $\C$, $\{z \in \C : 0 \leq \Im z \leq W\}$. 
See (\ref{eqn:cylinder_Ginibre}) in 
Section \ref{sec:strip_finite_infinite}, 
which represents the corresponding limit transitions. 
See also Remarks 5 and 8 in \cite{Kat19b}
in which the present Ginibre DPPs of types $C$ and $D$
are discussed as new examples of the 
{\it Mittag--Leffler fields} studied by \cite{AK13,AKM19,AKS18}.

\subsubsection{Ginibre and Ginibre-type DPPs 
as examples of Weyl--Heisenberg ensembles}
\label{sec:example_WH}
Let $d=1$ and consider the following window function, 
\begin{equation}
G(x_{\rR})= 2^{1/4} e^{-\pi x_{\rR}^2}, \quad x_{\rR} \in \R.
\label{eqn:g_A2}
\end{equation}
It is obvious that Assumption 4
is satisfied, $\|G\|^2_{L^2(\R, d x_{\rR})}=1$. 
In this case (\ref{eqn:K_WH}) becomes \cite{AGR16,APRT17,AGR19}
\[
K_{\rm WH}(x, x^{\prime})
= \frac{e^{\pi i x_{\rR} x_{\rI}}}{e^{\pi i x_{\rR}^{\prime} x_{\rI}^{\prime}}}
K_{{\rm Ginibre}}^A(\sqrt{\pi} x, \sqrt{\pi} x^{\prime}) e^{- \pi(|x|^2+|x'|^2)/2}. 
\]

By taking into account the direct decomposition
\[
L^2(\R)=L^2_{\rm odd}(\R) \oplus L^2_{\rm even}(\R),
\]
we have
\[
\cW^{\ast}_{\rm WH}(L^2(\R))
=\cW^{\ast}_{\rm WH}(L^2_{\rm odd}(\R)) \oplus
\cW^{\ast}_{\rm WH}(L^2_{\rm even}(\R)).
\]
When $G(-y)=G(y)$, we have
\begin{align*}
& \cW^{\ast}_{\rm WH}(L^2_{\rm odd}(\R))
\subset \{ F \in L^2(\C) : F(-x) =-F(x),  x \in \C \}
=: L^2_{\rm odd}(\C),
\nonumber\\
& \cW^{\ast}_{\rm WH}(L^2_{\rm even}(\R))
\subset \{ F \in L^2(\C) : F(-x) =F(x),  x \in \C \}
=: L^2_{\rm even}(\C).
\end{align*}

We consider the restriction of operator 
\[
\cW_{\rm WH} \Big|_{L^2_{\rm odd}(\C)} :
L^2_{\rm odd}(\C) \to L^2_{\rm odd}(\R).
\]
and its adjoint
\[
\cW^{\ast}_{\rm WH} \Big|_{L^2_{\rm odd}(\R)} :
L^2_{\rm odd}(\R) \to L^2_{\rm odd}(\C).
\]
Then, the kernel of the operator
$\cW_{\rm WH}^* \Big|_{L^2_{\rm odd}(\R)}
\cW_{\rm WH} \Big|_{L^2_{\rm odd}(\C)}$
is given by
\[
K_{\rm WH}^{\rm odd}(x, x^{\prime})
= \frac{1}{2} (K_{\rm WH}(x, x^{\prime})-K_{\rm WH}(x, -x^{\prime})).
\]
Similarly, we have the kernel, 
\[
K_{\rm WH}^{\rm even}(x, x^{\prime})
= \frac{1}{2} (K_{\rm WH}(x, x^{\prime})+K_{\rm WH}(x, -x^{\prime})).
\]
When the window function $G$ is given by (\ref{eqn:g_A2}), 
we obtain
\begin{align*}
K_{\rm WH}^{\rm odd}(x, x^{\prime})
&= \frac{e^{\pi i x_{\rR} x_{\rI}}}{e^{\pi i x_{\rR}^{\prime} x_{\rI}^{\prime}}}
K^{C}_{{\rm Ginibre}}(\sqrt{\pi} x, \sqrt{\pi} x^{\prime})e^{- \pi(|x|^2+|x'|^2)/2},
\nonumber\\
K_{\rm WH}^{\rm even}(x, x^{\prime})
&= \frac{e^{\pi i x_{\rR} x_{\rI}}}{e^{\pi i x_{\rR}^{\prime} x_{\rI}^{\prime}}}
K^{D}_{{\rm Ginibre}}(\sqrt{\pi} x, \sqrt{\pi} x^{\prime})e^{-\pi (|x|^2+|x'|^2)/2},
\quad x, x^{\prime} \in \C, 
\end{align*}
where $K^{C}_{{\rm Ginibre}}$ and $K^{D}_{{\rm Ginibre}}$
are
given by (\ref{eqn:K_Ginibre_C}) and (\ref{eqn:K_Ginibre_D}),
respectively.

The Ginibre DPP of type A is extended to 
{\it Ginibre-type DPPs} indexed by $q \in \N_0$,  
which are introduced in \cite{Shirai15} and 
also known as the infinite pure {\it polyanalytic ensembles} (cf. \cite{AGR19}). 
Each Ginibre-type DPP with index 
$q \in \N_0$ is associated with the correlation kernel 
\begin{equation}
 K^{(q)}_{\mathrm{Ginibre\text{-}type}}(x,x') := 
L_q^{(0)}(|x-x'|^2) K_{\mathrm{Ginibre}}^A(x, x'), 
\quad x, x' \in \C, 
\label{eqn:K_G_type}
\end{equation}
where $L_q^{(0)}$ is the $q$-th Laguerre polynomial (\ref{eqn:Laguerre1})
with parameter $\nu=0$ and $K_{\mathrm{Ginibre}}^A$ is defined by
(\ref{eqn:K_Ginibre_A}). 
This DPP can be viewed as 
the Weyl--Heisenberg ensemble,
$(\Xi, K_{\rm WH}^{h_q}, dx_{\rR} dx_{\rI})$, 
with the window function $G(x)=h_q(x)$, $x \in \R$, which is defined 
using the $q$-th Hermite polynomial (\ref{eqn:Hermite1}) as
\[
h_q(x) := \frac{2^{-q/2+1/4}}{\sqrt{q!}} e^{-\pi x^2}
H_q(\sqrt{2 \pi} x), \quad x \in \R, \quad q \in \N_0.
\]
Indeed, we see that 
\[
K_{\rm WH}^{h_q}(x, x^{\prime})
= \frac{e^{\pi i x_{\rR} x_{\rI}}}{e^{\pi i x_{\rR}^{\prime} x_{\rI}^{\prime}}}
K^{(q)}_{\mathrm{Ginibre\text{-}type}}(\sqrt{\pi} x, \sqrt{\pi} x') 
e^{- \pi(|x|^2+|x'|^2)/2},
\quad x, x' \in \C, 
\quad q \in \N_0.
\]
See \cite{AGR19} for more details about the Weyl--Heisenberg
aspect of finite polyanalytic ensembles. 
Other examples of the Weyl--Heisenberg ensembles
are given in \cite{AGR16,APRT17,AGR19}. 

\subsubsection{Representations of Ginibre and Ginibre-type 
kernels in the Bargmann--Fock space and the eigenspaces of Landau levels}
\label{sec:Fock}
We consider an application of Corollary \ref{thm:main3}.
Let
$S_1=\C$
and $S_2=\N_0$ with 
$\lambda_1(dx)=\lambda_{\rN(0, 1; \C)}(d x)$.
We put
\begin{equation}
\varphi_n(x) := \frac{x^n}{\sqrt{n!}},
\quad n \in \N_0.
\label{eqn:BF_phi}
\end{equation}
Note that $\{\varphi_n(x)\}_{n \in \N_0}$ forms 
a complete orthonormal system 
of the {\it Bargmann--Fock space}, 
which is the space of square-integrable analytic functions
on $\C$ with respect to the complex Gaussian
measure; 
\[
\langle \varphi_n, \varphi_m \rangle_{L^2(\C, \lambda_{\rN(0,1;\C)})}
=\delta_{n m}, \quad n, m \in \N_0.
\]
We assume that $\Gamma=S_2=\N_0$.
We can see that 
$\|\varphi_{\cdot}(x) \|_{\ell^2(\Gamma)}
=\sum_{n \in \N_0} |x|^{2n}/n!=e^{|x|^2}$, $x \in \C$.
Hence Assumption 3' is satisfied.
By Corollary \ref{thm:main3}, 
we obtain the DPP on $\C$ in which
the correlation kernel 
with respect to $\lambda_{\rN(0, 1: \C)}$ is given by
\begin{align*}
K_{\rm BF}(x, x^{\prime}) 
&= \sum_{n \in \N_0} \varphi_n(x) \overline{\varphi_n(x^{\prime})}
=\sum_{n=0}^{\infty} \frac{(x \overline{x^{\prime}})^n}{n!}
\nonumber\\
&= e^{x \overline{x^{\prime}}},
\quad x, x^{\prime} \in \C.
\end{align*}
This is the reproducing kernel in the Bargmann--Fock space
and obtained DPP is identified with 
$(\Xi, K^A_{\rm Ginibre}, \lambda_{\rN(0,1; \C)}(dx))$.
See \cite{Shirai15,BQ17,AGR19}.

If we set $\Gamma=2 \N_0+1 =\{1,3,5,\dots\}$ or 
$\Gamma=2 \N_0=\{0,2,4,\dots\}$, 
we will obtain the DPPs with the following kernels 
\begin{align*}
K_{\rm BF}^{\rm odd}(x, x^{\prime})
&= \sum_{k=0}^{\infty} \frac{(x \overline{x^{\prime}})^{2k+1}}{(2k+1)!}
=\sinh(x \overline{x^{\prime}}), 
\nonumber\\
K_{\rm BF}^{\rm even}(x, x^{\prime})
&= \sum_{k=0}^{\infty} \frac{(x \overline{x^{\prime}})^{2k}}{(2k)!}
=\cosh(x \overline{x^{\prime}}),
\quad x, x^{\prime} \in \C. 
\end{align*}
The obtained DPPs are identified with
$(\Xi$, $K^C_{\rm Ginibre}$, $\lambda_{\rN(0,1; \C)}(dx))$ 
and $(\Xi$, $K^D_{\rm Ginibre}$, $\lambda_{\rN(0,1; \C)}(dx))$,
respectively.

The correlation kernel of Ginibre-type DPP
(\ref{eqn:K_G_type}) admits the similar representation 
in terms of the {\it complex Hermite polynomials} defined by 
\[
 H_{p,q}(\zeta, \zetabar) 
 := (-1)^{p+q} e^{\zeta \zetabar} 
 \frac{\partial^p}{\partial \zetabar^p} \frac{\partial^q}{\partial \zeta^q} 
e^{-\zeta \zetabar}, \quad \zeta \in \C, 
\quad p, q \in \N_0, 
\]
which were introduced by It\^{o} \cite{Ito53}. 
We note that their generating function is given by 
\[
\sum_{p=0}^{\infty} \sum_{q=0}^{\infty} 
H_{p,q}(\zeta,\zetabar) \frac{s^p t^q}{p!q!} 
= \exp(\zeta s + \zetabar t - st)
\]
and the set $\{ H_{p,q}(\zeta, \zetabar)/\sqrt{p!q!} : p, q
\in \N_0\}$ forms a complete orthonormal system of 
$L^2(\C, \lambda_{\rN(0, 1; \C)}(d\zeta))$. 
Let $S_1=\C$ and $S_2=\N_0$ with 
$\lambda_1(dx)=\lambda_{\rN(0, 1; \C)}(d x)$, 
and for fixed $q \in \N_0$, define
\[
\varphi^{(q)}_n(x) := \frac{1}{\sqrt{n!q!}} H_{n,q}(x, \overline{x}),
\quad x \in \C, 
\quad n \in \N_0.
\]
Then $\{\varphi^{(q)}_n(x)\}_{n \in \N_0}$ forms 
a complete orthonormal system of the eigenspace corresponding
to the $q$-th {\it Landau level}, which coincides with 
the Bargmann--Fock space when $q=0$. 
Since the following formula is known 
\[
L_q^{(0)}(|\zeta-\eta|^2) e^{\zeta \etabar}
= \sum_{p=0}^{\infty} \frac{1}{p!q!}
H_{p,q}(\zeta,\zetabar)
\overline{H_{p,q}(\eta,\etabar)}, 
\quad \zeta, \eta \in \C, 
\quad q \in \N_0, 
\]
we obtain the following expansion formula for (\ref{eqn:K_G_type}), 
\[
K^{(q)}_{\mathrm{Ginibre\text{-}type}}(x,x') 
= \sum_{n=0}^{\infty}  \varphi^{(q)}_n(x) \overline{\varphi^{(q)}_n(x')},
\quad x, x' \in \C, 
\quad q \in \N_0. 
\]
The obtained DPPs are identified with
$(\Xi$, $K^{(q)}_{\mathrm{Ginibre\text{-}type}}, 
\lambda_{\rN(0,1; \C)}(dx))$, $q \in \N_0$ 
constructed as Weyl--Heisenberg ensembles 
in Section \ref{sec:example_WH}. 

\subsubsection{Application of duality relations}
\label{sec:duality_application2}
We consider the simplified setting (\ref{eqn:W_n_x1}) of $W$
with (\ref{eqn:BF_phi}) and 
$\Gamma=\N_0$.  
If we set $\lambda_1(dx)=\lambda_{\rN(0,1;\C)}(dx)$, 
$\Lambda_1$ be a disk (i.e., two-dimensional 
ball) $\B^2_r$ with radius $r \in (0, \infty)$
centered at the origin in $S_1=\C \simeq \R^2$ and
$\Lambda_2=S_2=\N_0$ in (\ref{eqn:K_A}), we obtain
\begin{align*}
K_{\C}^{(\N_0)}(x, x')
&= \sum_{n=0}^{\infty} \varphi_n(x) \overline{\varphi_n(x')}
=e^{x \overline{x'}} 
\nonumber\\
&= K_{\rm Ginibre}^{A}(x, x'),
\quad x, x' \in \C,
\end{align*}
where $K_{\rm Ginibre}^A$ denotes the 
correlation kernel of the Ginibre DPP 
of type $A$, and
\begin{align*}
K^{(\B^2_r)}_{\N_0}(n, n')
&= \int_{\B^2_r} \overline{\varphi_n(x)} \varphi_{n'}(x) \lambda_{\rN(0,1;\C)}(dx)
= \frac{1}{\pi \sqrt{n! n'!}} \int_0^r ds \, e^{-s^2} s^{n+n'+1} 
\int_0^{2 \pi} d \theta \, e^{i \theta (n'-n)}
\nonumber\\
&= 2 \delta_{n n'} \int_0^r \frac{s^{2n+1} e^{-s^2} }{n!} ds
=\delta_{n n'} \int_0^{r^2} \lambda_{\Gamma(n+1, 1)} (du),
\quad n, n' \in \N_0.
\end{align*}
Define
\[
\lambda_n(r) :=
\int_0^{r^2} \frac{u^n e^{-u}}{n!} du
= \sum_{k=n+1}^{\infty} \frac{r^{2k} e^{-r^2}}{k!},
\quad n \in \N_0, \quad r \in (0, \infty),
\]
where the second equality is due to Eq.(4.1) in \cite{Shirai15}.
That is, if we write the Gamma distribution with parameters
$(a, b)$ as $\Gamma(a, b)$ (see Section \ref{sec:OP})
and the Poisson distribution with parameter $c$ as
${\rm Po}(c)$,
\[
\lambda_n(r) := \bP(R_n \leq r^2) = \bP(Y_{r^2} \geq n+1),
\]
provided $R_n \sim \Gamma(n+1, 1)$ and
$Y_{r^2} \sim {\rm Po}(r^2)$.
Then DPP 
$(\Xi_2^{(\B^2_r)}, K^{(\B^2_r)}_{\N_0})$ on $\N_0$ 
is the product measure 
$\bigotimes_{n \in \N_0} \mu_{\lambda_n(r)}^{\rm Bernoulli}$
under the natural identification between 
$\{0,1\}^{\N_0}$ and the power set of $\N_0$, 
where $\mu_p^{\rm Bernoulli}$ denotes the Bernoulli measure of probability $p \in [0, 1]$.
Theorem~\ref{thm:duality2} gives the duality relation
\[
\bP(\Xi_{\rm Ginibre}^A(\B^2_r)=m)
=
\bP(\Xi_2^{(\B^2_r)}(\N_0)=m), 
\quad
\forall m \in \N_0,
\]
where we have identified the DPP, 
$(\Xi_1^{(\N_0)}, K_{\C}^{(\N_0)}, \lambda_1(dx))$
with the Ginibre DPP of type A, 
$(\Xi_{\rm Ginibre}^A$, $K_{\rm Ginibre}^A, \lambda_{\rN(0, 1; \C)})$.
If we introduce a series of random variables
$X_n^{(r)} \in \{0, 1\}, n \in \N_0$,
which are mutually independent and 
$X_n^{(r)} \sim \mu_{\lambda_n(r)}^{\rm Bernoulli}$, $n \in \N_0$, 
then the above implies the equivalence in probability law
\[
\Xi_{\rm Ginibre}^A(\B^2_r) \law 
\Xi_2^{(\B^2_r)}(\N_0) 
\law \sum_{n \in \N_0} X_n^{(r)},
\quad r \in (0, \infty).
\]
Similarly, we have the following equalities
by the results in Section \ref{sec:Fock} and 
Theorem \ref{thm:duality2},
\[
\Xi_{\rm Ginibre}^C(\B^2_r) 
\law \sum_{n \in 2\N_0+1} X_n^{(r)},
\quad
\Xi_{\rm Ginibre}^D(\B^2_r) 
\law \sum_{n \in 2\N_0} X_n^{(r)},
\quad
r \in (0, \infty).
\]

The argument above is valid for general 
{\it radially symmetric DPPs}
associated with radially symmetric finite measure $\lambda_1(dx) =
p(|x|) dx$ on $\C$. 
Let $\varphi_n(x) = a_n x^n, n \in \N_0$ be an
orthonormal system in $L^2(\C, \lambda_1)$ where 
$a_n > 0, n \in \N_0$ are the normalization constants, 
and we set 
\begin{align*}
K_{\C}^{(\N_0)}(x, x')
&= \sum_{n=0}^{\infty} \varphi_n(x) \overline{\varphi_n(x')}
= \sum_{n=0}^{\infty} a_n^2 (x \overline{x'})^n 
\quad x, x' \in \C,
\nonumber\\
K^{(\B^2_r)}_{\N_0}(n, n')
&= \int_{\B^2_r} \overline{\varphi_n(x)} \varphi_{n'}(x)
 \lambda_1(dx)
=\delta_{n n'} \lambda_n(r) \quad n, n' \in \N_0, 
\end{align*}
where 
\[
\la_n(r) := \frac{1}{Z_n} \int_0^{r^2} u^n p(\sqrt{u}) du 
\]
with $Z_n = \int_0^{\infty} u^n p(\sqrt{u})du$. 
Then DPP 
$(\Xi_1^{(\N_0)}, K_{\C}^{(\N_0)}, p(|x|)dx)$ on $\C$
is radially symmetric 
and DPP 
$(\Xi_2^{(\B^2_r)}, K^{(\B^2_r)}_{\N_0})$ on $\N_0$
is again identified
with the product measure 
$\bigotimes_{n \in \N_0} \mu_{\lambda_n(r)}^{\rm
Bernoulli}$.  
For example, if $p(s) = \pi^{-1} e^{-s^2}$ and $a_n = 1/\sqrt{n!}$,
then 
$(\Xi_1^{(\N_0)}, K_{\C}^{(\N_0)}, p(|x|) dx)$
is the Ginibre DPP of type A. 
The function $\la_n(r)$ is considered as a probability distribution
function on $[0,\infty)$ and hence there exist independent random
variables $R_n, n \in \N_0$ such that 
\[
 \la_n(r) = \bP(R_n \le r^2). 
\]
If we define $X_n^{(r)} = \mathbf{1}_{\{R_n \le r^2\}}$ for
each $n \in
\N_0$, then Theorem \ref{thm:duality2} gives the duality relation
\[
\Xi_1^{(\N_0)}(\B^2_r) \law \Xi_2^{(\B^2_r)}(\N_0)
\law \sum_{n \in \N_0} X_n^{(r)},
\quad r \in (0, \infty).
\]
Indeed, $\{X_n^{(r)}, n \in \N_0\}$ 
are mutually independent $\{0,1\}$-valued random variables
whose laws are given by 
$\{\mu_{\lambda_n(r)}^{\rm Bernoulli}$, $n \in\N_0\}$. 
If we take a set $\La_2 \subset \N_0$, then DPP 
$(\Xi_1^{(\La_2)}, K_{\C}^{(\Lambda_2)}, p(|x|) dx)$
satisfies 
\[
\Xi_1^{(\La_2)}(\B^2_r) \law \Xi_2^{(\B^2_r)}(\Lambda_2)
\law \sum_{n \in \Lambda_2} X_n^{(r)},
\quad r \in (0, \infty).
\]
We note that if we write 
$\Xi_1^{(\N_0)} = \sum_{j} \delta_{X_j}$,
then 
$\sum_j \delta_{|X_j|^2}$
is equal to
$\sum_{n \in \N_0} \delta_{R_n}$ in law, which was 
discussed in Theorem 4.7.1 in \cite{HKPV09} by constructing
$\{R_n \}_{n \in \N_0}$ in
terms of {\it size-biased sampling}.  

\subsection{Finite DPPs on sphere $\S^2$} 
\label{sec:sphere_finite}
Let
$\S^2 := \{x \in \R^3 : \|x\|_{\R^3}=1\}$ be the
two-dimensional unit sphere centered at the
origin in the three-dimensional Euclidean space $\R^3$,
where $\| \cdot \|_{\R^3}$ denotes the
Euclidean distance in $\R^3$.
We will use the following coordinates for $x=(x^{(1)}, x^{(2)}, x^{(3)})$ on $\S^2$,
\begin{equation}
x^{(1)}=\sin \theta \cos \varphi, \quad
x^{(2)}=\sin \theta \sin \varphi, \quad
x^{(3)}=\cos \theta, \quad \theta \in [0, \pi], \quad \varphi \in [0, 2 \pi). 
\label{eqn:polar1}
\end{equation}
We consider the case that
$S_1=\S^2$ and $S_2=\N_0$, in which
we assume that $\lambda_1(d x)$
is given by 
the {\it Lebesgue surface area measure} 
$d \sigma_2(x)$
on $\S^2$ such that
\[
\lambda_1(dx)=d \sigma_2(x)=d \sigma_2(\theta, \varphi)
:=\sin \theta d \theta d \varphi,
\quad
\lambda_1(\S^2)=\sigma_2(\S^2)=4 \pi.
\]

For $n \in \{0, 1, \dots, N-1\}, N \in \N$, put
\begin{equation}
\varphi_n^{\S^2}(x)=\varphi_n^{\S^2}(\theta, \varphi)
:= \frac{1}{\sqrt{h_n}} e^{- i n \varphi} \sin^n \frac{\theta}{2} \cos^{N-1-n} \frac{\theta}{2},
\quad \theta \in [0, \pi], \quad \varphi \in [0, 2 \pi),
\label{eqn:ON_S2}
\end{equation}
with
\[
h_n=h_n^{(N)} := \frac{4 \pi}{N} {\binom{N-1}{n}}^{-1}.
\]
It is easy to confirm the following orthonormality relations on $\S^2$,
\[
\langle \varphi_n^{\S^2}(\cdot), \varphi_m^{\S^2}(\cdot) 
\rangle_{L^2(\S^2; d \sigma_2)}
=\int_0^{\pi} d \theta \int_0^{2 \pi} d \varphi \,
\varphi_n^{\S^2}(\theta, \varphi)
\overline{\varphi_m^{\S^2}(\theta, \varphi)} d\sigma_2(\theta, \varphi)
=\delta_{n m}, \quad
n, m \in \N_0.
\]
We set $\psi_1(\cdot, n)=\varphi_n^{\S^2}(\cdot)$, 
$n \in \Gamma :=\{0,1, \dots, N-1\}$, $N \in \N_0$, 
By the argument given in Remark 3 in Section \ref{sec:OP},
we see Assumption 3' is satisfied.
Then Corollary \ref{thm:main3} gives
the DPP with $N$ points on $\S^2$, $(\Xi, K_{\S^2}^{(N)}, d \sigma_2(x))$, 
whose correlation kernel is given by
\begin{align}
K^{(N)}_{\S^2}(x, x')
&= K^{(N)}_{\S^2}((\theta, \varphi), (\theta', \varphi'))
\nonumber\\
&= \frac{N}{4 \pi} \sum_{n=0}^{N-1}
\binom{N-1}{n} \left(e^{-i(\varphi-\varphi')} \sin \frac{\theta}{2} \sin \frac{\theta'}{2} \right)^n
\left(\cos \frac{\theta}{2} \cos \frac{\theta'}{2} \right)^{N-1-n}
\nonumber\\
&= \frac{N}{4 \pi}
\left(e^{-i(\varphi-\varphi')} \sin \frac{\theta}{2} \sin \frac{\theta'}{2} 
+\cos \frac{\theta}{2} \cos \frac{\theta'}{2} \right)^{N-1}.
\label{eqn:K_S2}
\end{align}
The density of points with respect to $d \sigma_2(x)$ is given by
\[
\rho(x)=K^{(N)}_{\S^2}(x, x)= \frac{N}{4 \pi}
=\mbox{constant},
\quad x \in \S^2.
\]
 
For two points 
$x=(\theta, \varphi)$ and $x'=(\theta', \varphi')$ on $\S^2$,
\begin{align*}
\|x-x'\|_{\R^3}^2 &=  (\sin \theta \cos \varphi- \sin \theta' \cos \varphi')^2
+(\sin \theta \sin \varphi-\sin \theta' \sin \varphi')^2
+(\cos \theta- \cos \theta')^2
\nonumber\\
&= |\Phi(x-x')|^2,
\end{align*}
with
\begin{align*}
\Phi(x-x') & := 
2 \left[ \sin \frac{\theta- \theta'}{2} \cos \frac{\varphi-\varphi'}{2}
- i \sin \frac{\theta+\theta'}{2} \sin \frac{\varphi-\varphi'}{2} \right]
\nonumber\\
&= 2  \cos \frac{\theta}{2} \cos \frac{\theta'}{2}
e^{i(\varphi+\varphi')/2}
\left[ e^{-i \varphi} \tan \frac{\theta}{2}
- e^{-i \varphi'} \tan \frac{\theta'}{2} \right].
\end{align*}
Then we can show that the 
probability density of this DPP with respect to 
$d \sigma_2(\x) =\prod_{j=1}^N d \sigma_2(x_j)$ is given
as
\[
\bp^{(N)}_{\S^2}(\x)
= \frac{1}{Z^{(N)}_{\S^2}}
\prod_{1 \leq j < k \leq N} \|x_k-x_j\|_{\R^3}^2,
\]
with
\[
Z^{(N)}_{\S^2}
=\frac{2^{N(N+1)} \pi^N}{(N!)^{N-1}}
\left( \prod_{j=1}^N (j-1)! \right)^2.
\]
Since $\|x-x^{\prime}\|_{\R^3}^2=2-2 x \cdot x^{\prime}$
for $x, x^{\prime} \in \S^2$, we have the equality
\[
\frac{1}{2}(1+x \cdot x^{\prime})
= \left|
e^{-i(\varphi-\varphi')} \sin \frac{\theta}{2} \sin \frac{\theta'}{2}
+\cos \frac{\theta}{2} \cos \frac{\theta'}{2}
\right|^2.
\]
Hence the absolute value of (\ref{eqn:K_S2}) 
is written as
\[
\left| K^{(N)}_{\S^2}(x, x') \right|
=\frac{N}{4 \pi} 
\left( \frac{1+x \cdot x^{\prime}}{2} \right)^{(N-1)/2},
\]
and hence the two-point correlation function
(\ref{eqn:rho2}) with respect to
$d \sigma_2(x)$ is given by
\[
\rho^2(x, x^{\prime})
=\left( \frac{N}{4 \pi} \right)^2
\left[ 1-\left( \frac{1+x \cdot x^{\prime}}{2} \right)^{N-1} \right],
\quad x, x^{\prime} \in \S^2.
\]
The system $(\Xi, K^{(N)}_{\S^2}, d \sigma_2(x))$ is
uniform and isotropic on $\S^2$,
which is called the {\it spherical ensemble}
\cite{Kri09,For10,AZ15,BE18,BE19}.

\vskip 0.3cm
\noindent{\bf Remark 5} \, 
Let $G_1$ and $G_2$ be $N \times N$ independent
random matrices, whose entries are
i.i.d. following $\rN(0,1; \C)$.
Krishnapur \cite{Kri09} studied the statistical ensemble
of the eigenvalues $\z=(z_1, \dots, z_N)$ on $\C$
of $G_1^{-1} G_2$ and proved that
it gives the DPP $(\sum_{j} \delta_{Z_j}, K^{(N)}_{G_1^{-1} G_2}, \lambda(dz))$
with
\[
K^{(N)}_{G_1^{-1} G_2}(z, z') = (1+z \overline{z'})^{N-1},
\quad
\lambda(dz)=\frac{N}{\pi} \frac{dz}{(1+|z|^2)^{N+1}},
\]
which implies that the probability density of $\z$
with respect to the Lebesgue measure $d \z=\prod_{j=1}^N dz_j$
on $\C$ is given by
\[
\bp^{(N)}_{G_1^{-1} G_2}(\z)
=\frac{1}{Z_{G_1^{-1} G_2}^{(N)}}
\prod_{1 \leq j < k \leq N} |z_k-z_j|^2
\prod_{\ell=1}^N \frac{1}{(1+|z_{\ell}|^2)^{N+1}}
\]
with a normalization constant $Z_{G_1^{-1} G_2}^{(N)}$.
Krishnapur claimed that if we consider the stereographic projection
from $\S^2$ to $\widehat{\C} := \C \cup \{\infty\}$
which makes an equatorial plane of $\S^2$,
then the DPP, 
$(\sum_j \delta_{Z_j}, K_{G_1^{-1} G_2}^{(N)}, \lambda(dz))$
is realized as the image of the DPP in the spherical ensemble, 
$(\Xi$, $K_{\S^2}^{(N)}$, $d \sigma_2(x))$ \cite{Kri09}.
Actually if we consider the stereographic projection
such that the north pole of $\S^2$ ($\theta=0$) is mapped
to the origin of $\widehat{\C}$ and the south pole of $\S^2$
($\theta=\pi$) is to $\infty$,
the image of 
$x=(\sin \theta \cos \varphi, \sin \theta \sin \varphi, \cos \theta)
\in \S^2$ is given by
\[
z= e^{i \varphi} \tan \frac{\theta}{2} \in \widehat{\C},
\quad \theta \in [0, \pi], \quad \varphi \in [0, 2 \pi).
\]
We see that
\[
\|x-x'\|_{\R^3}^2= \frac{4}{(1+|z|^2)(1+|z'|^2)} |z-z'|^2
\]
and
\[
d \sigma_2(x)=\frac{4}{(1+|z|^2)^2} dz.
\]
Hence we can verify the statement of Krishnapur \cite{Kri09}.
\vskip 0.3cm

The equivalent system with the spherical ensemble of DPP
was studied by Caillol \cite{Cai81}
as a {\it two-dimensional one-component plasma model}
in physics. 
It is interesting to see that he used the
{\it Cayley--Klein parameters} defined by
\[
\alpha := e^{i \varphi/2} \cos \frac{\theta}{2},
\quad
\beta := -i e^{-i \varphi/2} \sin \frac{\theta}{2},
\quad \varphi \in [0, 2 \pi), \quad \theta \in [0, \pi].
\]
The orthonormal functions (\ref{eqn:ON_S2}) can be
identified with the following up to irrelevant factors,
\[
\widetilde{\varphi}_n^{\S^2}(\alpha, \beta)
=\frac{1}{\sqrt{h_n}} \alpha^{N-1-n} \beta^n, 
\quad n \in \{0,1, \dots, N-1 \}.
\]
If we define
\[
\langle (\alpha, \beta), (\alpha', \beta') \rangle_{\rm CK}
: =\alpha \overline{\alpha'} + \beta \overline{\beta'},
\]
the correlation kernel (\ref{eqn:K_S2}) is written as
\[
K_{\S^2}^{(N)}(x, x')
= K_{\S^2}^{(N)}((\alpha, \beta), (\alpha', \beta'))
=\frac{N}{4 \pi} \left( \langle (\alpha, \beta), (\alpha', \beta') \rangle_{\rm CK} \right)^{N-1}.
\]

Following the claim given in \cite{Cai81}
(see also Section 15.6.2 in \cite{For10}), 
we consider the vicinity of the north pole,
$x_{\rm np}=(0,0,1) \in \R^3$,
that is $\theta \sim 0$.
We put
\[
\theta=\frac{2r}{\sqrt{N}}, \quad
\theta^{\prime}= \frac{2 r'}{\sqrt{N}},
\]
and take the limit $N \to \infty$ keeping
$r$ and $r'$ be constants.
Then in (\ref{eqn:K_S2}), we see that
\begin{align*}
& \sin \frac{\theta}{2} \sin \frac{\theta'}{2} 
\sim \frac{1}{4} \theta \theta^{\prime}
= \frac{r r^{\prime}}{N},
\nonumber\\
& \cos \frac{\theta}{2} \cos \frac{\theta'}{2} \sim 
1-\frac{\theta^2+{\theta'}^2}{8}
=1- \frac{r^2+{r'}^2}{2N}.
\end{align*}
We set $r e^{i \varphi}=z, r' e^{i \varphi'}=z' \in \C$
with $r dr d \varphi= dz$.
Then the kernel given by (\ref{eqn:K_S2}) 
multiplied by $d \sigma_2$ 
has the following limit,
\begin{align*}
& \lim_{N \to \infty}
K^{(N)}_{\S^2}((\theta, \varphi), (\theta', \varphi')) d \sigma_2(\theta, \varphi)
\Big|_{\theta=2 |z|/\sqrt{N}, \theta'=2|z'|/\sqrt{N}}
\nonumber\\
& \quad
= \lim_{N \to \infty}
\frac{N}{4 \pi} \left( 1+ \frac{1}{N} \left\{z \overline{z^{\prime}}
- \frac{|z|^2+|z'|^2}{2} \right\} \right)^N \frac{4}{N} dz
\nonumber\\
& \quad
= \frac{1}{\pi} e^{z \overline{z^{\prime}} - (|z|^2+|z'|^2)/2} dz.
\end{align*}
Since the spherical ensemble is uniform and isotropic on $\S^2$,
we obtain the same limiting DPP in the vicinity 
of any point on $\S^2$.
This implies the following limit theorem \cite{Cai81,KS1}.
\begin{prop}
\label{thm:S2_Ginibre}
The following weak convergence is established,
\[
\frac{\sqrt{N}}{2}
\circ \Big(\Xi, K_{\S^2}^{(N)}, d \sigma_2(x) \Big)
\weak \Big(\Xi, K_{\rm Ginibre}^A, \lambda_{\rN(0, 1; \C)}(d x) \Big),
\]
where the limit point process is the
Ginibre DPP of type A given in Section \ref{sec:Ginibre_ACD}.
\end{prop}

\subsection{Finite DPPs on torus $\T^2$} 
\label{sec:rectangle_finite}

We will consider the finite DPPs on a surface of torus 
with double periodicity of 
$2 \omega_1 := 2 \pi$ and $2 \omega_3 :=2 \tau \pi$, 
where we assume that
$\tau=i \Im \tau \in i(0, \infty)$.
The surface of such a torus 
$\T^2=\T^2(2\pi, 2 \tau \pi) := \S^1(2 \pi) \times \S^1(2 \pi \Im \tau)$
can be identified with a rectangular domain in $\C$,
\[
D_{(2 \pi, 2 \tau \pi)} :=
\{z \in \C : 0 \leq \Re z < 2 \pi, 0 \leq \Im z < 2 \pi \Im \tau \} \subset \C
\quad \mbox{with double periodicity of $(2 \pi, 2 \tau \pi)$}. 
\]
So we first consider the systems on $D_{(2 \pi, 2 \tau \pi)}$.

Let $S=\C$ with 
$\lambda(d x) =
{\bf 1}_{D_{(2\pi, 2 \tau \pi)}}(x) dx_{\rR} dx_{\rI}$.
For $N \in \N$, put
\[
\varphi^{\RN, (2 \pi, 2 \tau \pi)}_n(x)
:= 
\frac{e^{-\cN^{\RN} i x_{\rI}^2/(4 \tau \pi)}}{\sqrt{h^{\RN}_n(\tau)}} 
\Theta^{\sharp(\RN)} \left( \frac{J^{\RN}(n)}{\cN^{\RN}}, 
\cN^{\RN} \frac{x}{2 \pi}, \cN^{\RN} \tau \right), 
\quad
n \in \{1, \dots, N\}.
\]
where
$\Theta^{R}$, 
$\sharp(\RN)$, $\cN^{\RN}$, and
$J^{\RN}(n)$ are given by 
(\ref{eqn:RS}), (\ref{eqn:sharp}), 
(\ref{eqn:N_R}), and (\ref{eqn:J_R}), respectively, 
and
\begin{align*}
h^{\AN}_n(\tau)
&:= 
4 \pi^2 \sqrt{\frac{\Im \tau}{2 \cN^{\AN}}}
e^{-2 \tau \pi i {J^{\AN}(n)}^{2}/ \cN^{\AN}},
\quad n \in \{1, \dots, N \},
\nonumber\\
h^{\RN}_n(\tau)
&:= 
8 \pi^2 \sqrt{\frac{\Im \tau}{2 \cN^{\RN}}}
e^{-2 \tau \pi i {J^{\RN}(n)}^{2}/\cN^{\RN}},
\quad n \in \{1, \dots, N\},
\quad \mbox{for $\RN=\CN, \CNv, \BCN$},
\nonumber\\
h^{\RN}_n(\tau)
&:=  \begin{cases}
\displaystyle{16 \pi^2 \sqrt{\frac{\Im \tau}{2 \cN^{\RN} }} },
\quad &n=1,
\cr
& \cr
\displaystyle{8 \pi^2 \sqrt{\frac{\Im \tau}{2 \cN^{\RN}}}
e^{- 2 \tau \pi i {J^{\RN}(n)}^{2}/\cN^{\RN}}},
\quad &n \in \{2,3, \dots, N\}, \quad
\end{cases}
\quad \mbox{for $\RN=\BN, \BNv$},
\nonumber\\
h^{\DN}_n(\tau)
&:=  \begin{cases}
\displaystyle{16 \pi^2 \sqrt{\frac{\Im \tau}{2 \cN^{\RN}}}
e^{- 2 \tau \pi i {J^{\DN}(n)}^{2}/\cN^{\DN}}},
\quad &n \in \{1, N \}, 
\cr
& \cr
\displaystyle{8 \pi^2 \sqrt{\frac{\Im \tau}{2 \cN^{\RN}}}
e^{- 2 \tau \pi i {J^{\DN}(n)}^{2}/\cN^{\DN}}},
\quad &n \in \{2,3, \dots, N-1 \}.
\end{cases}
\end{align*}
The following orthonormality relations were proved in \cite{Kat19b},
\[
\langle \varphi^{\RN, (2\pi, 2 \tau \pi)}_n,
 \varphi^{\RN, (2\pi, 2 \tau \pi )}_m
\rangle_{L^2(\C, {\bf 1}_{D_{(2\pi, 2 \tau \pi)}}(x) d x)} = \delta_{nm},
\quad n, m \in \Gamma :=\{1, \dots, N\}, 
\]
$\RN=\AN, \BN, \BNv, \CN, \CNv, \BCN, \DN$. 
By the argument given in Remark 3 in Section \ref{sec:OP},
we see Assumption 3 is satisfied.
Then Corollary \ref{thm:main2} gives
the seven types of DPPs with the correlation kernels,
\begin{equation}
K^{\RN, (2 \pi, 2 \tau \pi)}(x, x^{\prime})
= \sum_{n=1}^N \varphi^{\RN, (2 \pi, 2 \tau \pi )}_n(x) 
\overline{\varphi^{\RN, (2 \pi, 2 \tau \pi)}_n(x^{\prime})},
\label{eqn:K_torus}
\end{equation}
with respect to the measure
$\lambda(d x)={\bf 1}_{D_{(2\pi, 2 \tau \pi)}} d x$
on $\C$ for
$\RN=\AN, \BN, \BNv, \CN, \CNv, \BCN, \DN$. 

Using the quasi-periodicity of the Jacobi theta functions 
(\ref{eqn:quasi_periodic1}) and (\ref{eqn:quasi_periodic2}),
we can show that the correlation kernels are quasi-double-periodic as \cite{Kat19b},
\begin{align*}
K^{\RN, (2 \pi, 2 \tau \pi)}(x+2 \pi, x^{\prime}) 
&= K^{\RN, (2 \pi, 2 \tau \pi)}(x, x^{\prime}+2 \pi)
\nonumber\\
&=\begin{cases}
(-1)^{\cN^{\AN}} K^{\RN, (2 \pi, 2 \tau \pi )}(x, x^{\prime}),
\quad & \RN=\AN,
\cr
- K^{\RN, (2 \pi, 2 \tau \pi)}(x, x^{\prime}),
\quad & \RN=\BN, \CNv, \BCN,
\cr
K^{\RN, (2 \pi, 2 \tau \pi)}(x, x^{\prime}),
\quad & \RN=\BNv, \CN, \DN,
\end{cases}
\nonumber\\
K^{\RN, (2 \pi, 2 \tau \pi)}(x+ 2 \tau \pi, x^{\prime})
&= \begin{cases}
e^{- \cN^{\RN} i x_{\rR}} K^{\RN, (2 \pi, 2 \tau \pi)}(x, x^{\prime}),
\quad & \RN=\AN, \CN, \CNv, \BCN, \DN,
\cr
-e^{- \cN^{\RN} i x_{\rR}} K^{\RN, (2 \pi, 2 \tau \pi)}(x, x^{\prime}),
\quad & \RN=\BN, \BNv,
\end{cases}
\nonumber\\
K^{\RN, (2 \pi, 2 \tau \pi)}(x, x^{\prime}+2 \tau \pi)
&= \begin{cases}
e^{\cN^{\RN} i x_{\rR}^{\prime}} K^{\RN, (2 \pi, 2 \tau \pi)}(x, x^{\prime}),
\quad & \RN=\AN, \CN, \CNv, \BCN, \DN,
\cr
-e^{\cN^{\RN} i x_{\rR}^{\prime}} K^{\RN, (2 \pi, 2 \tau \pi)}(x, x^{\prime}),
\quad & \RN=\BN, \BNv.
\end{cases}
\end{align*}
The above implies the following double periodicity 
(up to an irrelevant gauge transformation), 
\begin{align*}
\cS_{2 \pi} K^{\RN, (2 \pi, 2 \tau \pi)}(x, x')
&= \frac{e^{\cN^{\RN} i x_{\rR}}}{e^{\cN^{\RN} i x_{\rR}^{\prime}}}
\cS_{2 \tau \pi} K^{\RN, (2 \pi, 2 \tau \pi)}(x, x' )
\nonumber\\
&=K^{\RN, (2 \pi, 2 \tau \pi)}(x, x'),
\quad x, x' \in D_{(2\pi, 2 \tau \pi)}.
\end{align*}
In other words, we have obtained the seven types
of DPPs with a finite number of points $N$ 
on a surface of torus $\T^2(2 \pi, 2 \tau \pi)$.
Hence here we write them as 
$\left(\Xi, K^{\RN}_{\T^2(2 \pi, 2 \tau \pi)}, d x \right)$,
$\RN=\AN$, $\BN$, $\BNv$, $\CN$, $\CNv$, $\BCN$, $\DN$. 
Using the {\it Macdonald denominator formula}
given by Rosengren and Schlosser \cite{RS06}
(see (2.6) in \cite{Kat19b} in the present notations), 
the probability densities for these finite DPPs
with respect to 
the Lebesgue measures,
$d \x =\prod_{j=1}^N d x_j$ are given
as follows; 
\begin{align}
\bp^{\AN}_{\T^2(2 \pi, 2 \tau \pi)}(\x) 
&= \frac{1}{Z^{\AN}_{\T^2(2 \pi, 2 \tau \pi)}}
\exp\left(-\frac{\cN^{\AN}}{2 \pi \Im \tau} \sum_{j=1}^N (x_j)_{\rI}^2 \right)
\nonumber\\
&
\times
\begin{cases}
\displaystyle{
\left|
\vartheta_0 \left( \sum_{k=1}^N \frac{x_k}{2\pi}; \tau \right)
W^{\AN}\left( \frac{\x}{2 \pi}; \tau \right) 
\right|^2
},
& \mbox{if $N$ is even}, \cr
\displaystyle{
\left|
\vartheta_3 \left( \sum_{k=1}^N \frac{x_k}{2\pi}; \tau \right)
W^{\AN}\left( \frac{\x}{2 \pi}; \tau \right) 
\right|^2
},
& \mbox{if $N$ is odd}, 
\end{cases}
\nonumber\\
\bp^{\RN}_{\T^2(2 \pi, 2 \tau \pi)}(x)
&=\frac{1}{Z^{\RN}_{\T^2(2 \pi, 2 \tau \pi)}} 
\exp\left(-\frac{\cN^{\RN}}{2 \pi \Im \tau} \sum_{j=1}^N (x_j)_{\rI}^2 \right)
\left|
W^{\RN} \left( \frac{\x}{2 \pi}; \tau \right) 
\right|^2,
\nonumber\\
& \qquad \RN=\BN, \BNv, \CN, \CNv, \BCN, \DN, 
\label{eqn:p_torus}
\end{align}
for $\x \in (\T^2(2 \pi, 2 \tau \pi))^N$,
where $W^{\RN}$ are the Macdonald denominators
given by (\ref{eqn:Macdonald_denominators})
in Appendix \ref{sec:Macdonald_denominators}
and $Z^{\RN}_{\T^2(2 \pi, 2 \tau \pi)}$ are normalization constants \cite{Kat19b}.

We can prove the following symmetry properties
for the present DPPs on $\T^2(2\pi, 2 \tau \pi)$.
\begin{prop}
\label{thm:symmetry_torus}
\begin{description}
\item{\rm (i)} \,
The finite DPPs
$\left(\Xi, K^{\RN}_{\T^2(2 \pi, 2 \tau \pi)}, d x \right)$
with 
$\tau=i \Im \tau \in i(0, \infty)$
have the following shift invariance, 
\begin{align*}
& \cS_{2 \pi/N} (\Xi, K_{\T^2(2 \pi, 2 \tau \pi)}^{\AN}, dx)
\law (\Xi, K_{\T^2(2 \pi, 2 \tau \pi)}^{\AN}, dx),
\nonumber\\
& \cS_{2 \tau \pi/N} (\Xi, K_{\T^2(2 \pi, 2 \tau \pi)}^{\AN}, dx)
\law (\Xi, K_{\T^2(2 \pi, 2 \tau \pi)}^{\AN}, dx),
\nonumber\\
& \cS_{\pi} (\Xi, K_{\T^2(2 \pi, 2 \tau \pi)}^{\RN}, dx)
\law (\Xi, K_{\T^2(2 \pi, 2 \tau \pi)}^{\RN}, dx),
\quad \RN=\BNv, \CN, \DN,
\nonumber\\
& \cS_{\tau \pi} (\Xi, K_{\T^2(2 \pi, 2 \tau \pi)}^{\RN}, dx)
\law (\Xi, K_{\T^2(2 \pi, 2 \tau \pi)}^{\RN}, dx),
\quad \RN=\CN, \CNv, \BCN, \DN.
\end{align*}
\item{\rm (ii)} \,
The densities of points 
$\rho_{\T^2(2 \pi, 2 \tau \pi)}^{\RN}(x)$
given by $K_{\T^2(2 \pi, 2 \tau \pi)}^{\RN}(x, x)$
have the following zeros,
\begin{align*}
&\rho^{\BN}_{\T^2(2 \pi, 2 \tau \pi)}(0)=0, 
\nonumber\\
&\rho^{\BNv}_{\T^2(2 \pi, 2 \tau \pi)}(0)= \rho^{\BNv}_{\T^2(2 \pi, 2 \tau \pi)}(\pi)=0, 
\nonumber\\
&\rho^{\RN}_{\T^2(2 \pi, 2 \tau \pi)}(0)=\rho^{\RN}_{\T^2(2 \pi, 2 \tau \pi)}(\tau \pi)=0, \quad
\RN=\CNv, \BCN,
\nonumber\\
&\rho^{\CN}_{\T^2(2 \pi, 2 \tau \pi)}(0)
= \rho^{\CN}_{\T^2(2 \pi, 2 \tau \pi)}(\pi)
=\rho^{\CN}_{\T^2(2 \pi, 2 \tau \pi)}(\tau \pi)=0.
\end{align*}
\end{description}
\end{prop}
\noindent{\it Proof} \quad
(i) It is easy to verify the statements if we use the
formulas (\ref{eqn:p_torus}) with (\ref{eqn:Macdonald_denominators})
for the probability densities.
Use the formulas 
(\ref{eqn:quasi_periodic1})--(\ref{eqn:half_tau}),
which show the change of values of $\vartheta_{\mu}(v; \tau)$,
$\mu=0,1,3$, due to the shift of variable
$v \to v+1, v \to v+\tau$ and $v \to v+\tau/2$,
respectively.
For the shift $\cS_{\tau \pi}$, note the fact that
\[
\cS_{\tau \pi} \exp \left(
-\frac{\cN^{\RN}}{2 \pi \Im \tau} \sum_{j=1}^N (x_j)_{\rm I}^2 \right)
=\exp \left(
-\frac{\cN^{\RN}}{2 \pi \Im \tau} \sum_{j=1}^N (x_j)_{\rm I}^2 \right)
\prod_{\ell=1}^N e^{2 \pi i \cN^{\RN}(i x_{\ell}/2 \pi + \tau/4)}.
\]
As a matter of course, the statements can be proved
also by showing the shift invariance of
the correlation kernels (\ref{eqn:K_torus}) 
up to irrelevant gauge transformations. 
(ii) By the properties (\ref{eqn:theta_values}) 
of the Jacobi theta functions,
the zeros of densities $\rho_{\T^2(2 \pi, 2 \tau \pi)}^{\RN}(x)$
are determined as above. 
Then the proof is complete. \qed
\vskip 0.3cm

We note that the 
periods $2 \pi/N \in (0, \infty)$
and $2 \tau \pi/N \in i (0, \infty)$
of $(\Xi, K_{\T^2(2 \pi, 2 \tau \pi)}^{\AN}, dx)$
shown by Proposition \ref{thm:symmetry_torus} (i) 
become zeros as $N \to \infty$.
Hence, as the $N \to \infty$ limit of
$(\Xi, K_{\T^2(2 \pi, 2 \tau \pi)}^{\AN}, dx)$,
it is expected to obtain a uniform system
of infinite number of points on $\C$.
Actually we can prove the following limit
theorems.
\begin{prop}
\label{thm:torus_Ginibre}
The following weak convergence is established,
\begin{align*}
\frac{1}{2} \sqrt{\frac{N}{\pi \Im \tau}}
\circ \Big(\Xi, K^{\AN}_{\T^2(2 \pi, 2 \tau \pi)}, d x \Big)
&\weak \Big(\Xi, K_{\rm Ginibre}^A, \lambda_{\rN(0, 1; \C)}(d x) \Big),
\nonumber\\
\sqrt{\frac{N}{2 \pi \Im \tau}}
\circ \Big(\Xi, K^{\RN}_{\T^2(2 \pi, 2 \tau \pi)}, d x \Big)
&\weak \Big(\Xi, K_{\rm Ginibre}^C, \lambda_{\rN(0, 1; \C)}(d x) \Big),
\quad \RN=\BN, \BNv, \CN, \CNv, \BCN,
\nonumber\\
\sqrt{\frac{N}{2 \pi \Im \tau}}
\circ \Big(\Xi, K^{\DN}_{\T^2(2 \pi, 2 \tau \pi)}, d x \Big)
&\weak \Big(\Xi, K_{\rm Ginibre}^D, \lambda_{\rN(0, 1; \C)}(d x) \Big),
\end{align*}
where the limit point processes are the three types of 
Ginibre DPPs given in Section \ref{sec:Ginibre_ACD}.
\end{prop}
\noindent{\it Proof} \,
By (\ref{eqn:theta_asym}), we see that
\begin{align*}
\varphi_n^{\AN, (2 \pi, 2 \tau \pi)} 
\left( 2 \sqrt{\frac{\pi \Im \tau}{N}} x \right)
&\sim
N^{1/4} \frac{1}{2 \pi} \left( \frac{2}{\Im \tau} \right)^{1/4}
e^{-N \pi \Im \tau/4 -i \sqrt{N \pi \Im \tau} x -x_{\rI}^2} 
\nonumber\\
&\times
\exp \left[ - \pi \Im \tau \left( \frac{n-1/2}{\sqrt{N}} \right)^2
+(2 i \sqrt{\pi \Im \tau} x + \sqrt{N} \pi \Im \tau) \frac{n-1/2}{\sqrt{N}} \right],
\end{align*}
as $N \to \infty$ with $(n-1/2)/\sqrt{N}={\rm O}(1) > 0$.
Then
\begin{align*}
& K^{\AN}_{\T^2(2 \pi, 2 \tau \pi)}
\left( 2 \sqrt{\frac{\pi \Im \tau}{N}} x,
2 \sqrt{\frac{\pi \Im \tau}{N}} x^{\prime} 
\right) 
\left( 2 \sqrt{\frac{\pi \Im \tau}{N}} \right)^2 
\nonumber\\
& \quad
\sim \frac{1}{\pi} \sqrt{2 \Im \tau}
e^{-N \pi \Im \tau/2-i \sqrt{N \pi \Im \tau}(x-\overline{x^{\prime}})} 
e^{-(x_{\rI}^2+{x_{\rI}^{\prime}}^2)}
\int_0^{\infty} e^{-2 \pi \Im \tau u^2
+\{2 i \sqrt{\pi \Im \tau}(x-\overline{x^{\prime}})
+ 2 \sqrt{N} \pi \Im \tau \} u} du
\nonumber\\
& \quad
= \frac{1}{\pi} \sqrt{ 2 \Im \tau}
e^{-(x_{\rI}^2+{x_{\rI}^{\prime}}^2)-(x-\overline{x^{\prime}})^2/2}
\int_{-\sqrt{N}/2 -i(x-\overline{x^{\prime}})/\{2\sqrt{\pi \Im \tau} \}}^{\infty}
e^{-2 \pi \Im \tau v^2} dv
\nonumber\\
& \quad 
\rightarrow
\frac{1}{\pi} e^{-(x_{\rI}^2+{x_{\rI}^{\prime}}^2)-(x-\overline{x^{\prime}})^2/2}
\quad \mbox{as $N \to \infty$}. 
\end{align*}
It implies that
\[
\lim_{N \to \infty} 
K^{\AN}_{\T^2(2 \pi, 2 \tau \pi)}
\left( 2 \sqrt{\frac{\pi \Im \tau}{N}} x,
2 \sqrt{\frac{\pi \Im \tau}{N}} x^{\prime} 
\right) \left( 2 \sqrt{\frac{\pi \Im \tau}{N}} \right)^2 
=\frac{e^{-i x_{\rR} x_{\rI}}}{e^{-i x_{\rR}^{\prime} x_{\rI}^{\prime}}}
K^{A}_{\rm Ginibre}(x, x^{\prime}) 
\frac{1}{\pi} e^{-(|x|^2+|x^{\prime}|^2)/2}.
\]
Similarly, we can show that
\[
\lim_{N \to \infty} 
K^{\RN}_{\T^2(2 \pi, 2 \tau \pi)}
\left( \sqrt{\frac{2\pi \Im \tau}{N}} x,
\sqrt{\frac{2 \pi \Im \tau}{N}} x^{\prime} 
\right) \left( \sqrt{\frac{2 \pi \Im \tau}{N}} \right)^2 
=\frac{e^{-i x_{\rR} x_{\rI}}}{e^{-i x_{\rR}^{\prime} x_{\rI}^{\prime}}}
K^{C}_{\rm Ginibre}(x, x^{\prime}) 
\frac{1}{\pi} e^{-(|x|^2+|x^{\prime}|^2)/2},
\]
for $\RN=\BN, \BNv, \CN, \CNv, \BCN$,
and
\[
\lim_{N \to \infty} 
K^{\DN}_{\T^2(2 \pi, 2 \tau \pi)}
\left( \sqrt{\frac{2\pi \Im \tau}{N}} x,
\sqrt{\frac{2 \pi \Im \tau}{N}} x^{\prime} 
\right) \left( \sqrt{\frac{2 \pi \Im \tau}{N}} \right)^2 
=\frac{e^{-i x_{\rR} x_{\rI}}}{e^{-i x_{\rR}^{\prime} x_{\rI}^{\prime}}}
K^{D}_{\rm Ginibre}(x, x^{\prime}) 
\frac{1}{\pi} e^{-(|x|^2+|x^{\prime}|^2)/2}.
\]
Then the statement will be proved. \qed

\subsection{Finite and infinite DPPs on cylinder $\R \times \S^1$} 
\label{sec:strip_finite_infinite}

Here we consider the finite DPPs on a surface of
cylinder with infinite length having periodicity of $2 \pi \alpha$
in the circumference direction, $\alpha \in (0, \infty)$, 
which we write here as $\R \times \S^1(2 \pi \alpha)$.
The surface of $\R \times \S^1(2 \pi \alpha)$ can be identified with a strip 
with width $2 \pi \alpha$ in $\C$,
\[
D_{2 \pi \alpha} := 
\{z \in \C: 0 \leq \Im z < 2 \pi \alpha \} \subset \C
\quad \mbox{with periodicity of $2 \pi i \alpha$}. 
\]
So we first consider the systems on $D_{2 \pi \alpha}$.
Let $S=\C$. 
For $N \in \N$, we set
\[
\lambda(d x) 
=\lambda_{\rN(0, 1/4)}(d x_{\rR}) \lambda_{[0, 2 \pi \alpha)}(d x_{\rI})
=\frac{1}{\sqrt{2} \pi^{3/2} \alpha} e^{-2 x_{\rR}^2}
{\bf 1}_{D_{2 \pi \alpha}}(x) dx_{\rR} dx_{\rI}.
\]
Define
\begin{align*}
\varphi^{\AN, 2 \pi \alpha}_n(x) 
&:=
e^{-[ (\cN^{\AN}-2 J^{\AN}(n))^2/(16 \alpha^2)
+ (\cN^{\AN}-2 J^{\AN}(n) )x/(2 \alpha)]},
\quad n \in \{1, \dots, N\}, 
\nonumber\\
\varphi^{\RN, 2 \pi \alpha}_n(x)
&:= \sqrt{2}
e^{-(\cN^{\RN}-2J^{\RN}(n))^2/(16 \alpha^2)}
\sinh \left[ (\cN^{\RN}-2 J^{\RN}(n)) \frac{x}{2 \alpha} \right], 
\quad n \in \{1,\dots, N\}, 
\nonumber\\
& \hskip 10cm
\mbox{for $\RN=\BN, \CN$}, 
\nonumber\\
\varphi^{\DN, 2 \pi \alpha}_n(x)
&:= \begin{cases}
\displaystyle{
\sqrt{2} e^{-(\cN^{\DN}-2J^{\DN}(n))^2/(16 \alpha^2)}
\cosh \left[ (\cN^{\DN}-2 J^{\DN}(n)) \frac{x}{2 \alpha} \right]
},
& n \in \{1, \dots, N-1\},
\cr
1,
& n=N,
\cr
\end{cases}
\end{align*}
where $\cN^{\RN}$ and $J^{\RN}(n)$, $\RN=\AN, \BN, \CN, \DN$
are given by (\ref{eqn:N_R_classic}) and (\ref{eqn:J_R_classic}), respectively.
They have periodicity or quasi-periodicity of $2 \pi i \alpha$,
\begin{equation}
\varphi^{\RN, 2 \pi \alpha}_n(x+2 \pi i \alpha)=
\begin{cases}
(-1)^{N+1} \varphi^{\RN, 2 \pi \alpha}_n(x), & \RN=\AN,
\cr
- \varphi^{\RN, 2 \pi \alpha}_n(x), & \RN=\BN,
\cr
\varphi^{\RN, 2 \pi \alpha}_n(x), & \RN=\CN, \DN.
\end{cases}
\label{eqn:psi_iW}
\end{equation}

It is easy to verify the following orthonormality relations;
for $\RN=\AN, \BN, \CN, \DN$, 
\begin{align*}
\langle \varphi^{\RN, 2 \pi \alpha}_n,
\varphi^{\RN, w \pi \alpha}_m \rangle_{L^2(\C, \lambda_{\rN(0, 1/4)}(d x_{\rR}) \lambda_{[0, 2 \pi \alpha)}(d x_{\rI}))}
= \delta_{nm}, 
\quad
n, m \in \Gamma:=\{1, \dots, N\}.
\end{align*}
By the argument given in Remark 3 in Section \ref{sec:OP},
we see Assumption 3' is satisfied.
Then Corollary \ref{thm:main2}
gives the following four types of DPPs with the 
correlation kernels,
\[
K^{\RN, 2 \pi \alpha}(x, x^{\prime})
= \sum_{n=1}^N \varphi^{\RN, 2 \pi \alpha}_n(x) 
\overline{\varphi^{\RN, 2 \pi \alpha}_n(x^{\prime})},
\quad \RN=\AN, \BN, \CN, \DN.
\]
From (\ref{eqn:psi_iW}), we can see the
periodicity or quasi-periodicity of $2 \pi i \alpha$ 
in the correlation kernels, 
\[
K^{\RN, 2 \pi \alpha}(x+2 \pi i \alpha, x^{\prime})=
K^{\RN, 2 \pi \alpha}(x, x^{\prime}+ 2 \pi i \alpha)
=\begin{cases}
(-1)^{N+1} K^{\RN, 2 \pi \alpha}(x,x^{\prime}), & \RN=\AN,
\cr
- K^{\RN, 2 \pi \alpha}(x,x^{\prime}), & \RN=\BN,
\cr
K^{\RN, 2 \pi \alpha}(x,x^{\prime}), & \RN=\CN, \DN,
\end{cases}
\]
which implies
\[
\cS_{2 \pi i \alpha} K^{\RN, 2 \pi \alpha}(x, x')
=K^{\RN, 2 \pi \alpha}(x, x'),
\quad x, x' \in D_{2 \pi \alpha}.
\]
That is, we have obtained the four types
of DPPs with $N$ points
on a surface of cylinder, $\R \times \S^1(2 \pi \alpha)$.
Hence here we write them as 
$\left(\Xi, K^{\RN}_{\R \times \S^1(2 \pi \alpha)}, 
\lambda_{\rN(0, 1/4)}(d x_{\rR}) \lambda_{[0, 2 \pi \alpha)}(d x_{\rI}) \right)$,
$\RN=\AN$, $\BN$, $\CN$, $\DN$.
(See \cite{CJ88} for related systems in two dimensions.)

By Lemma \ref{thm:Weyl_trigonometric}
in Appendix \ref{sec:Weyl_denominators}, 
the probability densities for these finite DPPs
with respect to 
the Lebesgue measures,
$d \x =\prod_{j=1}^N d x_j$ are given
as follows; for $\x \in (\R \times \S^1(2 \pi \alpha))^N$, 
$\alpha \in (0, \infty)$, 
\begin{align}
\bp^{\AN}_{\alpha}(\x) 
&= \frac{1}{Z^{\AN}_{\alpha}} \prod_{1 \leq j < k \leq N}
\sinh^2 \frac{x_k-x_j}{2 \alpha}
\prod_{\ell=1}^N e^{-2(x_{\ell})_{\rR}^2}
\nonumber\\
\bp^{\BN}_{\alpha}(\x) 
&= \frac{1}{Z^{\BN}_{\alpha}}
\prod_{1 \leq j < k \leq N} 
\left( \sinh^2 \frac{x_k-x_j}{2 \alpha} \sinh^2 \frac{x_k+x_j}{2 \alpha} \right)
\prod_{\ell=1}^N \left( \sinh^2 \frac{x_{\ell}}{2 \alpha} e^{-2(x_{\ell})_{\rR}^2} \right), 
\nonumber\\
\bp^{\CN}_{\alpha}(\x) 
&= \frac{1}{Z^{\CN}_{\alpha}}
\prod_{1 \leq j < k \leq N} 
\left( \sinh^2 \frac{x_k-x_j}{2 \alpha} \sinh^2 \frac{x_k+x_j}{2 \alpha} \right)
\prod_{\ell=1}^N \left( \sinh^2 \frac{x_{\ell}}{\alpha} e^{-2(x_{\ell})_{\rR}^2} \right),
\nonumber\\
\bp^{\DN}_{\alpha}(\x) 
&= \frac{1}{Z^{\DN}_{\alpha}} 
\prod_{1 \leq j < k \leq N} 
\left( \sinh^2 \frac{x_k-x_j}{2 \alpha} \sinh^2 \frac{x_k+x_j}{2 \alpha} \right)
\prod_{\ell=1}^N e^{-2(x_{\ell})_{\rR}^2},
\label{eqn:p_cylinder}
\end{align}
with normalization constants $Z^{\RN}_{\alpha}$.

If we use the formulas (\ref{eqn:p_cylinder}),
it is easy to verify the following symmetry properties.
\begin{prop}
\label{thm:symmetry_cylinder}
\begin{description}
\item{\rm (i)} \,
The infinite DPPs 
$\left(\Xi, K^{\RN}_{\R \times \S^1(2 \pi \alpha)}, 
\lambda_{\rN(0, 1/4)}(d x_{\rR}) \lambda_{[0, 2 \pi \alpha)}(d x_{\rI}) \right)$ with
$\alpha \in (0, \infty)$
have the following shift invariance, 
\begin{align*}
&\cS_{i \theta} 
\left(\Xi, K^{\AN}_{\R \times \S^1(2 \pi \alpha)}, 
\lambda_{\rN(0, 1/4)}(d x_{\rR}) \lambda_{[0, 2 \pi \alpha)}(d x_{\rI}) \right)
\nonumber\\
& \qquad 
\law
\left(\Xi, K^{\AN}_{\R \times \S^1(2 \pi \alpha)}, 
\lambda_{\rN(0, 1/4)}(d x_{\rR}) \lambda_{[0, 2 \pi \alpha)}(d x_{\rI}) \right), 
\quad \forall \theta \in [0, 2 \pi \alpha),
\nonumber\\
&\cS_{i \pi \alpha} 
\left(\Xi, K^{\RN}_{\R \times \S^1(2 \pi \alpha)}, 
\lambda_{\rN(0, 1/4)}(d x_{\rR}) \lambda_{[0, 2 \pi \alpha)}(d x_{\rI}) \right)
\nonumber\\
& \qquad 
\law
\left(\Xi, K^{\RN}_{\R \times \S^1(2 \pi \alpha)}, 
\lambda_{\rN(0, 1/4)}(d x_{\rR}) \lambda_{[0, 2 \pi \alpha)}(d x_{\rI}) \right), 
\quad \RN=\CN, \DN.
\end{align*}
\item{\rm (ii)} \,
The densities of points $\rho_{\R \times \S^1(2 \pi \alpha)}^{\RN}(x)$
given by $K_{\R \times \S^1(2 \pi \alpha)}^{\RN}(x, x)$
have the following zeros,
\[
\rho_{\R \times \S^1(2 \pi \alpha)}^{\BN}(0)=0, \quad
\rho_{\R \times \S^1(2 \pi \alpha)}^{\CN}(0)=
\rho_{\R \times \S^1(2 \pi \alpha)}^{\CN}(i \pi \alpha)=0.
\]
\end{description}
\end{prop}
\vskip 0.3cm

Using the Jacobi theta functions (\ref{eqn:theta}), 
the limits of the correlation kernels in $N \to \infty$ can be
expressed as follows,
\begin{align}
\lim_{\ell \to \infty} K^{A_{2 \ell-1}}_{\R \times \S^1(2 \pi \alpha)}(x, x')
&= \vartheta_2 \left( \frac{i(x+\overline{x^{\prime}})}{2 \pi \alpha};
\frac{i}{2 \pi \alpha^2} \right)
=: K^{A, {\rm even}}_{\R \times \S^1(2 \pi \alpha)}(x, x'),
\nonumber\\
\lim_{\ell \to \infty} K^{A_{2 \ell}}_{\R \times \S^1(2 \pi \alpha)}(x, x')
&= \vartheta_3 \left( \frac{i(x+\overline{x^{\prime}})}{2 \pi \alpha};
\frac{i}{2 \pi \alpha^2} \right)
=: K^{A, {\rm odd}}_{\R \times \S^1(2 \pi \alpha)}(x, x'),
\nonumber\\
\lim_{N \to \infty} K^{\BN}_{\R \times \S^1(2 \pi \alpha)}(x, x')
&= \frac{1}{2} \left\{
\vartheta_2 \left( \frac{i(x+\overline{x^{\prime}})}{2 \pi \alpha}; 
\frac{i}{2 \pi \alpha^2} \right)
-
\vartheta_2 \left( \frac{i(x-\overline{x^{\prime}})}{2 \pi \alpha}; 
\frac{i}{2 \pi \alpha^2} \right) 
\right\}, 
\nonumber\\
&=: K^{B}_{\R \times \S^1(2 \pi \alpha)}(x, x')
\nonumber\\
\lim_{N \to \infty} K^{\CN}_{\R \times \S^1(2 \pi \alpha)}(x, x')
&= \frac{1}{2} \left\{
\vartheta_3 \left( \frac{i(x+\overline{x^{\prime}})}{2 \pi \alpha}; 
\frac{i}{2 \pi \alpha^2} \right)
-
\vartheta_3 \left( \frac{i(x-\overline{x^{\prime}})}{2 \pi \alpha}; 
\frac{i}{2 \pi \alpha^2} \right) 
\right\}
\nonumber\\
&=: K^{C}_{\R \times \S^1(2 \pi \alpha)}(x, x')
\nonumber\\
\lim_{N \to \infty} K^{\DN}_{\R \times \S^1(2 \pi \alpha)}(x, x')
&= \frac{1}{2} \left\{
\vartheta_3 \left( \frac{i(x+\overline{x^{\prime}})}{2 \pi \alpha}; 
\frac{i}{2 \pi \alpha^2} \right)
+
\vartheta_3 \left( \frac{i(x-\overline{x^{\prime}})}{2 \pi \alpha}; 
\frac{i}{2 \pi \alpha^2} \right) 
\right\}
\nonumber\\
&=: K^{D}_{\R \times \S^1(2 \pi \alpha)}(x, x'), 
\label{eqn:kernels_cylinder1}
\end{align}
$x, x' \in \R \times \S^1(2 \pi \alpha)$.

\begin{prop}
\label{thm:infiniteDPP_cylinder}
The five limit kernels 
{\rm (\ref{eqn:kernels_cylinder1})}
define the five kinds of infinite DPPs 
on the cylinder,
$(\Xi, K^{A, \sharp}_{\R \times \S^1(2 \pi \alpha)}, \lambda_{\rN(0, 1/4)}(d x_{\rR}) \lambda_{[0, 2 \pi \alpha)}(d x_{\rI}))$,
$\sharp=$even, odd, and
$(\Xi$, $K^{R}_{\R \times \S^1(2 \pi \alpha)}$, 
$\lambda_{\rN(0, 1/4)}(d x_{\rR})$ $\lambda_{[0, 2 \pi \alpha)}(d x_{\rI}))$,
$R=B, C, D$.
\end{prop}

The particle densities at 
$x \in D_{2 \pi \alpha} \simeq \R \times \S^1(2 \pi \alpha)$
are obtained from the limit kernels (\ref{eqn:kernels_cylinder1}) 
by setting $x'=x$. 
Since $x+\overline{x}=2 x_{\rR}$, $x-\overline{x}=2 i x_{\rI}$,
the definitions of the Jacobi theta functions (\ref{eqn:theta})
with parity (\ref{eqn:even_odd})
give explicit expressions for them; for instance,
\begin{align*}
\rho_{\R \times \S^1(2 \pi \alpha)}^B(x)
&= \frac{1}{2} 
 \left\{
\vartheta_2 \left( \frac{i x_{\rR}}{\pi \alpha}; 
\frac{i}{2 \pi \alpha^2} \right)
-
\vartheta_2 \left( \frac{ x_{\rI}}{2 \pi \alpha}; 
\frac{i}{2 \pi \alpha^2} \right) 
\right\}
\nonumber\\
&= \sum_{n=1}^{\infty}
e^{-(n-1/2)^2/2 \alpha^2}
\Big[ \cosh\{(2n-1) x_{\rR}/\alpha\} - \cos \{(2n-1) x_{\rI}/\alpha\} \Big],
\nonumber\\
\rho_{\R \times \S^1(2 \pi \alpha)}^C(x)
&= \frac{1}{2} 
 \left\{
\vartheta_3 \left( \frac{i x_{\rR}}{\pi \alpha}; 
\frac{i}{2 \pi \alpha^2} \right)
-
\vartheta_3 \left( \frac{ x_{\rI}}{2 \pi \alpha}; 
\frac{i}{2 \pi \alpha^2} \right) 
\right\}
\nonumber\\
&= \sum_{n=1}^{\infty}
e^{-n^2/2 \alpha^2}
\Big[ \cosh(2n x_{\rR}/\alpha) - \cos (2n x_{\rI}/\alpha) \Big].
\end{align*}
They show that the obtained particle densities are
indeed nonnegative, and that
\[
\rho_{\R \times \S^1(2 \pi \alpha)}^B(0)=0, \quad
\rho_{\R \times \S^1(2 \pi \alpha)}^C(0)
=\rho_{\R \times \S^1(2 \pi \alpha)}^C( i \pi \alpha)=0.
\]

\vskip 0.3cm
\noindent{\bf Remark 6} \, 
In \cite{Kat19b}, the infinite DPPs on a strip in $\C$ 
were introduced 
by taking an anisotropic scaling limit
associated with $N \to \infty$ 
of the doubly periodic DPPs
(Theorem 3.4 in \cite{Kat19b}).
There the limiting correlation kernels
are expressed by the integrals
of products of Jacobi's theta functions.
We have found that, 
if we correctly perform Jacobi's imaginary transformations 
(\ref{eqn:Jacobi_imaginary}) of the integrands, 
the integrals can be calculated 
and the results are identified with the correlation kernels
simply given by (\ref{eqn:kernels_cylinder1}). 
\vskip 0.3cm

Using the quasi-periodicity of the Jacobi theta functions (\ref{eqn:quasi_periodic2}), 
we can show that, for $\sharp={\rm even, odd}$,
\[
\cS_{1/(2 \alpha)} K^{A, \sharp}_{\R \times \S^1(2 \pi \alpha)}(x, x^{\prime})
\sqrt{\frac{2}{\pi}} e^{-(x_{\rR}^2+{x_{\rR}^{\prime}}^2)} 
\frac{1}{2 \pi \alpha}
= \frac{e^{i x_{\rI}/\alpha}}{e^{i x_{\rI}'/\alpha}}
K^{A, \sharp}_{\R \times \S^1(2 \pi \alpha)}(x, x^{\prime})
\sqrt{\frac{2}{\pi}} e^{-(x_{\rR}^2+{x_{\rR}^{\prime}}^2)} 
\frac{1}{2 \pi \alpha}.
\]
By the gauge invariance, this implies the 
shift invariance, 
\begin{align}
&\cS_{1/(2 \alpha)}
\left(\Xi, K^{A, \sharp}_{\R \times \S^1(2 \pi \alpha)}, 
\lambda_{\rN(0, 1/4)}(d x_{\rR}) \lambda_{[0, 2 \pi \alpha)}(d x_{\rI}) \right)
\nonumber\\
& \qquad 
\law \left(\Xi, K^{A, \sharp}_{\R \times \S^1(2 \pi \alpha)}, 
\lambda_{\rN(0, 1/4)}(d x_{\rR}) \lambda_{[0, 2 \pi \alpha)}(d x_{\rI}) \right),
\label{eqn:shift_typeA_cylinder}
\end{align}
$\sharp = \rm{even, odd}$.
Moreover, by the properties (\ref{eqn:half_tau}) of the Jacobi theta functions,
we obtain the equality
\[
\cS_{1/(4 \alpha)} K^{A, {\rm even}}_{\R \times \S^1(2 \pi \alpha)}(x, x^{\prime})
\sqrt{\frac{2}{\pi}} e^{-(x_{\rR}^2+{x_{\rR}^{\prime}}^2)} 
\frac{1}{2 \pi \alpha}
= \frac{e^{i x_{\rI}/(2 \alpha)}}{e^{i x_{\rI}'/(2 \alpha)}}
K^{A, {\rm odd}}_{\R \times \S^1(2 \pi \alpha)}(x, x^{\prime})
\sqrt{\frac{2}{\pi}} e^{-(x_{\rR}^2+{x_{\rR}^{\prime}}^2)} 
\frac{1}{2 \pi \alpha}.
\]
Hence again by the gauge invariance, 
\[
\cS_{1/(4 \alpha)} 
\left(\Xi, K^{A, {\rm even}}_{\R \times \S^1(2 \pi \alpha)}, 
\lambda_{\rN(0, 1/4)}(d x_{\rR}) \lambda_{[0, 2 \pi \alpha)}(d x_{\rI}) \right)
\law
\left(\Xi, K^{A, {\rm odd}}_{\R \times \S^1(2 \pi \alpha)}, 
\lambda_{\rN(0, 1/4)}(d x_{\rR}) \lambda_{[0, 2 \pi \alpha)}(d x_{\rI}) \right),
\]
that is, the even-limit and the odd-limit of type $A$
are equivalent up to the shift by $1/(4 \alpha)$ in the real-axis direction.

We note that the 
period $1/(2 \alpha) \in (0, \infty)$
of $\left(\Xi, K^{A, \sharp}_{\R \times \S^1(2 \pi \alpha)}, 
\lambda_{\rN(0, 1/4)}(d x_{\rR}) \lambda_{[0, 2 \pi \alpha)}(d x_{\rI}) \right)$,
$\sharp=$ even, odd, shown by (\ref{eqn:shift_typeA_cylinder}) 
becomes zero as $\alpha \to \infty$.
Hence, as the $\alpha \to \infty$ limits of these DPPs, 
a uniform system of infinite number of points on $\C$ will be obtained.
In order to see such limit transitions, first
we perform Jacobi's imaginary transformations
(\ref{eqn:Jacobi_imaginary}). 

\begin{lem}
\label{thm:infiniteDPP_cylinder_Jacobi_Imaginary}
The following equalities hold,
\begin{align*}
& K^{A, {\rm even}}_{\R \times \S^1(2 \pi \alpha)}(x, x^{\prime})
\sqrt{\frac{2}{\pi}} e^{-(x_{\rR}^2+{x_{\rR}^{\prime}}^2)} 
\frac{1}{2 \pi \alpha}
= \frac{e^{i x_{\rR} x_{\rI}}}{e^{i x^{\prime}_{\rR} x^{\prime}_{\rI}}}
e^{x \overline{x^{\prime}}} 
\frac{1}{\pi} e^{-(|x|^2+|x^{\prime}|^2)/2}
\vartheta_0((x+\overline{x^{\prime}}) \alpha; 2 \pi i \alpha^2),
\nonumber\\
& K^{A, {\rm odd}}_{\R \times \S^1(2 \pi \alpha)}(x, x^{\prime})
\sqrt{\frac{2}{\pi}} e^{-(x_{\rR}^2+{x_{\rR}^{\prime}}^2)} 
\frac{1}{2 \pi \alpha}
= \frac{e^{i x_{\rR} x_{\rI}}}{e^{i x^{\prime}_{\rR} x^{\prime}_{\rI}}}
e^{x \overline{x^{\prime}}} 
\frac{1}{\pi} e^{-(|x|^2+|x^{\prime}|^2)/2}
\vartheta_3((x+\overline{x^{\prime}}) \alpha; 2 \pi i \alpha^2),
\nonumber\\
& K^B_{\R \times \S^1(2 \pi \alpha)}(x, x^{\prime})
\sqrt{\frac{2}{\pi}} e^{-(x_{\rR}^2+{x_{\rR}^{\prime}}^2)} 
\frac{1}{2 \pi \alpha}
\nonumber\\
& \quad 
= \frac{e^{i x_{\rR} x_{\rI}}}{e^{i x^{\prime}_{\rR} x^{\prime}_{\rI}}}
\Big\{ e^{x \overline{x^{\prime}}} 
\vartheta_0( (x+\overline{x^{\prime}}) \alpha; 2 \pi i \alpha^2)
-
e^{- x \overline{x^{\prime}}} 
\vartheta_0( (x-\overline{x^{\prime}}) \alpha; 2 \pi i \alpha^2) \Big\}
\frac{1}{2 \pi} e^{-(|x|^2+|x^{\prime}|^2)/2},
\nonumber\\
& K^C_{\R \times \S^1(2 \pi \alpha)}(x, x^{\prime})
\sqrt{\frac{2}{\pi}} e^{-(x_{\rR}^2+{x_{\rR}^{\prime}}^2)} 
\frac{1}{2 \pi \alpha}
\nonumber\\
& \quad 
= \frac{e^{i x_{\rR} x_{\rI}}}{e^{i x^{\prime}_{\rR} x^{\prime}_{\rI}}}
\Big\{ e^{x \overline{x^{\prime}}} 
\vartheta_3( (x+\overline{x^{\prime}}) \alpha; 2 \pi i \alpha^2)
-
e^{- x \overline{x^{\prime}}} 
\vartheta_3( (x-\overline{x^{\prime}}) \alpha; 2 \pi i \alpha^2) \Big\}
\frac{1}{2 \pi} e^{-(|x|^2+|x^{\prime}|^2)/2},
\nonumber\\
& K^D_{\R \times \S^1(2 \pi \alpha)}(x, x^{\prime})
\sqrt{\frac{2}{\pi}} e^{-(x_{\rR}^2+{x_{\rR}^{\prime}}^2)} 
\frac{1}{2 \pi \alpha}
\nonumber\\
& \quad 
= \frac{e^{i x_{\rR} x_{\rI}}}{e^{i x^{\prime}_{\rR} x^{\prime}_{\rI}}}
\Big\{ e^{x \overline{x^{\prime}}} 
\vartheta_3( (x+\overline{x^{\prime}}) \alpha; 2 \pi i \alpha^2)
+
e^{- x \overline{x^{\prime}}} 
\vartheta_3( (x-\overline{x^{\prime}}) \alpha; 2 \pi i \alpha^2) \Big\}
\frac{1}{2 \pi} e^{-(|x|^2+|x^{\prime}|^2)/2}.
\end{align*}
\end{lem}
By the asymptotics of the Jacobi theta functions (\ref{eqn:theta_asym}),
the following limit transitions are immediately concluded
from the expressions in Lemma \ref{thm:infiniteDPP_cylinder_Jacobi_Imaginary}.
\begin{prop}
\label{thm:alpha_inf}
The following limit transitions 
from the four types of infinite DPPs on
$\R \times \S^1(2 \pi \alpha)$
to the three types of Ginibre DPPs on $\C$ are
established,
\begin{align}
\left(\Xi, K^{A, \sharp}_{\R \times \S^1(2 \pi \alpha)}, 
\lambda_{\rN(0, 1/4)}(d x_{\rR}) \lambda_{[0, 2 \pi \alpha)}(d x_{\rI}) \right)
& \weaka
\Big(\Xi, K^A_{{\rm Ginibre}}, \lambda_{\rN(0, 1; \C)}(d x) \Big), 
\quad \sharp={\rm even, odd}, 
\nonumber\\
\left.
\begin{array}{l}
\displaystyle{
\left(\Xi, K^{B}_{\R \times \S^1(2 \pi \alpha)}, 
\lambda_{\rN(0, 1/4)}(d x_{\rR}) \lambda_{[0, 2 \pi \alpha)}(d x_{\rI}) \right)
} \cr
\displaystyle{
\left(\Xi, K^{C}_{\R \times \S^1(2 \pi \alpha)}, 
\lambda_{\rN(0, 1/4)}(d x_{\rR}) \lambda_{[0, 2 \pi \alpha)}(d x_{\rI}) \right)
}
\end{array}
\right\}
&\weaka
\Big(\Xi, K^C_{{\rm Ginibre}}, \lambda_{\rN(0, 1; \C)}(d x) \Big),
\nonumber\\
\left(\Xi, K^{D}_{\R \times \S^1(2 \pi \alpha)}, 
\lambda_{\rN(0, 1/4)}(d x_{\rR}) \lambda_{[0, 2 \pi \alpha)}(d x_{\rI}) \right)
&\weaka
\Big(\Xi, K^D_{{\rm Ginibre}}, \lambda_{\rN(0, 1; \C)}(d x) \Big).
\label{eqn:cylinder_Ginibre}
\end{align}
\end{prop}

\SSC
{Examples in Spaces with Arbitrary Dimensions} 
\label{sec:example_d_dim}
\subsection{Heisenberg family of infinite DPPs on $\C^{d}$}
\label{sec:HeisenbergDPP}

The Ginibre DPP of type A on $\C$ given in Section \ref{sec:Ginibre_ACD}
can be generalized to the DPPs on $\C^d$ for $d \geq 2$.
This generalization was done by \cite{AGR16,APRT17,AGR19}
as the Weyl--Heisenberg ensembles of DPP,
but here we derive the DPPs on $\C^d$, $d \in \N$,
following Corollary \ref{thm:main3} given in
Section \ref{sec:orthogonal}. 

Let $S_1=\C^d$, $S_2=\Gamma=\R^d$,  
\begin{align*}
\lambda_1(d x) & = \prod_{a=1}^d \lambda_{\rN(0, 1; \C)}(d x^{(a)})
= \frac{1}{\pi^d} e^{-|x|^2} = \frac{1}{\pi^d} e^{-(|x_{\rR}|^2+|x_{\rI}|^2)}
\nonumber\\
&=: \lambda_{\rN(0, 1; \C^d)}(d x),
\nonumber\\
\lambda_2(d y) &= \nu(dy) 
= \prod_{a=1}^d \lambda_{\rN(0, 1/4)}(d y^{(a)})
=\left( \frac{2}{\pi} \right)^{d/2} e^{-2 |y|^2}, 
\end{align*}
and
\[
\psi_1(x, y)
= e^{-(|x_{\rR}|^2-|x_{\rI}|^2)/2 + 2 (x_{\rR} \cdot y + i x_{\rI} \cdot y)},
\quad x=x_{\rR}+i x_{\rI} \in \C^d, \quad y \in \R^d.
\]
We see that 
$\Psi_1(x)^2:=\|\psi_1(x, \cdot) \|_{L^2(\R^d, \nu)}^2=e^{|x|^2}$,
$x \in \C^d$. Hence Assumption 3' is satisfied
and then, by Corollary \ref{thm:main3}, 
we obtain the DPP on $\C^d$ with the correlation kernel,
\begin{align*}
K^{(d)}(x, x') 
&= \left( \frac{2}{\pi} \right)^{d/2}
e^{-\{(|x_{\rR}|^2-|x_{\rI}|^2)+(|x^{\prime}_{\rR}|^2-|x^{\prime}_{\rI}|^2)\}/2}
\int_{\R^d} e^{-2 [ |y|^2-\{(x_{\rR}+i x_{\rI}) + (x^{\prime}_{\rR}-i x^{\prime}_{\rI}) \} \cdot y ]}
d y
\nonumber\\
&= \frac{e^{i x_{\rR} \cdot x_{\rI}}}{e^{i x^{\prime}_{\rR} \cdot x^{\prime}_{\rI}}}
K_{\rm Heisenberg}^{(d)}(x, x')
\end{align*}
with 
\[
K^{(d)}_{\rm Heisenberg}(x, x^{\prime})
= e^{x \cdot \overline{x^{\prime}}}, \quad x, x^{\prime} \in \C^d.
\]

The kernels in this form on $\C^d, d \in \N$ have been
studied by Zelditch and his coworkers
(see \cite{Zel00,BSZ00} and references therein),
who identified them with the Szeg\"o kernels
for the reduced Heisenberg group \cite{Fol89,Ste93,Gro01}. 
Here we call the DPPs associated with
the correlation kernels in this form
the {\it Heisenberg family of DPPs} on $\C^d, d \in \N$.
This family includes the Ginibre DPP of type A as the
lowest dimensional case with $d=1$.

\begin{df}
\label{thm:Heisenberg_class}
The Heisenberg family of DPPs is a one-parameter
{\rm ($d \in \N$)} family of \\
$\Big(\Xi, K^{(d)}_{\rm Heisenberg}, \lambda_{\rN(0, 1; \C^d)}(d x) \Big)$
with
\[
K^{(d)}_{\rm Heisenberg}(x, x')
:= e^{x \cdot \overline{x'}},
\quad x, x' \in \C^d.
\]
\end{df}
\vskip 0.3cm
Since
\[
K^{(d)}_{\rm Heisenberg}(x, x) \lambda_{\rN(0, 1; \C^d)}(dx)
= \frac{1}{\pi^d} d x, \quad x \in \C^d,
\]
every DPP in the Heisenberg family 
is uniform on $\C^d$ and the density with respect to
the Lebesgue measure $dx$ is given by $1/\pi^d$.
Hyperuniformity \cite{Tor18} of the Heisenberg family of DPPs
has been studied in \cite{APRT17,MKS21}.

\subsection{Finite DPPs on $\S^d$}
\label{sec:DPP_S_d}

First we recall basic properties of
spherical harmonics on $\S^d$ \cite{Nom18}.
For $d \in \N$, let
$\cP=\cP(\R^{d+1})$ be a vector space of all complex-valued
polynomials on $\R^{d+1}$,
and $\cP_k, k \in \N_0$, be its subspaces
consisting of homogeneous polynomials of degree $k$;
$p(x)=\sum_{|\alpha|=k} c_{\alpha} x^{\alpha}$,
$c_{\alpha} \in \C, x=(x^{(1)}, \dots, x^{(d+1)}) \in \R^{d+1}$,
where we have used the notations
$x^{\alpha} := \prod_{a=1}^{d+1} (x^{(a)})^{\alpha_a}$
with $\alpha:=(\alpha_1, \dots, \alpha_{d+1}) \in \N_0^{d+1}$,
$|\alpha| := \sum_{a=1}^{d+1} \alpha_a$. 
The vector space of all harmonic functions in $\cP$ is
denoted by $\cH =\{p \in \cP : \Delta p=0\}$
and let $\cH_k = \cH \cap \cP_k$, $k \in \N_0$.

Now we consider a unit sphere in $\R^{d+1}$ denoted by
$\S^{d}$, in which we use  the polar coordinates 
for $x=(x^{(1)}, \dots, x^{(d+1)}) \in \S^{d}$, 
\begin{align}
x^{(1)} &= \sin \theta_{d} \cdots \sin \theta_2 \sin \theta_1,
\nonumber\\
x^{(a)} &= \sin \theta_{d} \cdots \sin \theta_a \cos \theta_{a-1},
\quad a=2, \dots, d,
\nonumber\\
x^{(d+1)} &= \cos \theta_{d},
\quad \mbox{with $\theta_1 \in [0, 2 \pi), \quad
\theta_a \in [0, \pi], \quad a=2, \dots, d$}.
\label{eqn:coordinate1}
\end{align}
Note that $\|x\|^2_{\R^{d+1}} := \sum_{a=1}^{d+1} {x^{(a)}}^2 =1$. 
For $d=2$, if we put
$\theta_1=\pi/2-\varphi$ and $\theta_2=\theta$,
the polar coordinates (\ref{eqn:polar1}) 
used in Section \ref{sec:sphere_finite} are obtained.
The standard measure on $\S^{d}$ is given
by the {\it Lebesgue area measure} expressed as
\begin{equation}
d \sigma_{d} (x)
= \sin^{d-1} \theta_{d} \sin^{d-2} \theta_{d-1}
\cdots \sin \theta_2 d \theta_1 \cdots d \theta_{d},
\quad x \in \S^{d}.
\label{eqn:sigma}
\end{equation}
The total measure of $\S^{d}$ is calculated as
\begin{equation}
\omega_{d} := \sigma_{d}(\S^{d})
= \frac{2 \pi^{(d+1)/2}}{\Gamma((d+1)/2)}.
\label{eqn:omega}
\end{equation}
We write the space of harmonic polynomials in $\cH_k$ 
restricted on $\S^{d}$ as
\[
\cY_{(d, k)} = \left\{ h\Big|_{\S^{d}} : h \in \cH_k \right\},
\quad k \in \N_0.
\]
We can see that
\begin{equation}
D(d, k) := \dim \cY_{(d, k)}=
\frac{(d+2k-1)(d+k-2)!}{(d-1)! k!}
=\frac{2}{(d-1)!} k^{d-1}+{\rm o}(k^{d-1}).
\label{eqn:dim}
\end{equation}

Consider an orthonormal basis 
$\{Y^{(d,k)}_j \}_{j=1}^{D(d,k)}$ of $\cY_{(d,k)}$
with respect to $d \sigma_{d}$;
\begin{equation}
\langle Y^{(d, k)}_n, Y^{(d,k)}_m \rangle_{L^2(\S^d, d \sigma_d)}
=
\int_{\S^{d}} Y^{(d, k)}_n(x) \overline{Y^{(d, k)}_m(x)} d \sigma_{d} (x) 
=\delta_{nm}, \quad n, m \in \N_0.
\label{eqn:orthonormal}
\end{equation}
Then, if we put
\[
K^{\cY_{(d, k)}}(x, x')
= \sum_{j=1}^{D(d, k)} Y^{(d, k)}_j(x) \overline{Y^{(d, k)}_j(x')},
\quad x' \in \S^{d}, 
\]
then $\{K^{\cY_{(d, k)}}(x, x')\}_{x, x' \in \S^{d}}$
give the reproducing kernel in $\cY^{(d,k)}$ in the sense that
\[ 
Y(x') = \int_{\S^{d}} Y(x) \overline{K^{\cY_{(d, k)}}(x, x')} d \sigma_{d}(x),
\quad \forall Y \in \cY_{(d, k)}.
\]

For $\lambda > -1/2$, we define
\[
P^{\lambda}_{k}(x)
:= 
{_2}F_1 \left( -k, k+2 \lambda; \lambda+ \frac{1}{2} ; 
\frac{1-x}{2} \right),
\]
where ${_2}F_1$ denotes the {\it Gauss hypergeometric function},
$
{_2}F_1(\alpha, \beta, \gamma; z)
:= \sum_{n=0}^{\infty} \{(\alpha)_n (\beta)_n/\gamma)_n \} z^n/n!
$,
with 
$(\alpha)_n := \alpha(\alpha+1) \cdots (\alpha+n-1)
= \Gamma(\alpha+n)/\Gamma(\alpha)$, $n \in \N$, 
$(\alpha)_0 := 1$. 
Then the following equality is established,
\[
K^{\cY_{(d, k)}}(x, x')
=\frac{D(d, k)}{\omega_{d}} P^{(d-1)/2}_k(x \cdot x'),
\quad x, x' \in \S^{d},
\]
where $\omega_{d}$ and $D(d,k)$ are
given by (\ref{eqn:omega}) and (\ref{eqn:dim}), respectively,
and $x \cdot x' :=\sum_{a=1}^{d+1} x^{(a)} {x'}^{(a)}$. 

We see that $K^{\cY_{(d, k)}}(x, x')$ is $\gO(d+1, \R)$-invariant in the sense
\[
K^{\cY_{(d, k)}}(gx, gx')= K^{\cY_{(d, k)}}(x, x'),
\quad \forall g \in \gO(d+1, \R), \quad
\forall x, x' \in \S^{d}.
\]
Let $\{e_1, \dots, e_{d+1}\}$ be the standard basis of 
$\R^{d+1}$
and $L_0$ be
the stabilizer subgroup of
$\SO(d+1, \R)$ at $e_{d+1}$ represented as
\[
L_0 = \left\{
\left( \begin{array}{ll}
A & {\bf 0} \cr
{^{t}{\bf 0}} & 1
\end{array} 
\right) 
: A \in \SO(d, \R) \right\}.
\]
We define
\[
\cY_{(d, k)}^{L_0} =
\{ Y \in \cY_{(d, k)} : Y(\ell x)=Y(x), \quad
\forall \ell \in L_0, \quad \forall x \in \S^{d} \},
\]
Then $K^{\cY_{(d, k)}} \in \cY_{(d, k)}^{L_0}$.
The function $P^{\lambda}_k(s)$ is called
the {\it ultraspherical polynomial} \cite{NIST10}.
The space $\cY_{(d, k)}^{L_0}$ is a one-dimensional vector space
generated by $P_k^{(d-1)/2}(x \cdot e_{d+1})$.
In general, any $L_0$-invariant function is a constant
for each $L_0$-orbit,
$\cO =\{x \in \S^{d} : x \cdot e_{d+1} = a \}, 
-1 \leq a \leq 1$,
and hence functions in $\cY_{(d, k)}^{L_0}$ is called
the {\it zonal harmonics} of degree $k$. 
Note that, when we set
\[
C_k^{\lambda}(x) := \binom{k+2 \lambda-1}{k} 
P^{\lambda}_k(x),
\]
we call $C_k^{\lambda}(s)$
the {\it Gegenbauer polynomial} of degree $k$ \cite{NIST10}.

Fix $d \in \N$ and $k \in \N_0$.
Then, if we consider the case that
$S_1=\S^{d}$,
$S_2=\N$ with
$\lambda_1(d x)=d \sigma_{d}(x)$,
$L^2(\Gamma, \nu)=\ell^2( \{1, \dots, D(d, k)\}) \subset S_2$,
and $\psi_1(x, n)=Y^{(d, k)}_n(x)$,
then (\ref{eqn:orthonormal}) 
with Remark 3 in Section \ref{sec:OP}
guarantees Assumption 3'.
Hence Corollary \ref{thm:main3} determines
a unique DPP on $\S^{d}$, in which
the correlation kernel is given by \cite{BMOC16}
\begin{align*}
K^{\cY_{(d, k)}}(x, x')
&= \frac{D(d,k)}{\omega_{d}} P^{(d-1)/2}_k(x \cdot x')
\nonumber\\
&= \frac{d-1+2k}{(d-1) \omega_{d}} C^{(d-1)/2}_k(x \cdot x').
\end{align*}
It is obvious that the obtained DPP is
rotationally invariant on $\S^d$, since the kernel
$K^{\cY_{(d, k)}}(x, x')$ depend only on the inner product $x \cdot x'$.
The density of points is uniform on $\S^{d}$
and is given with respect to $\sigma_d(dx)$ by
\begin{align*}
\rho^{\cY_{(d, k)}}
&= K^{\cY_{(d, k)}}(x, x)
= \frac{D(d,k)}{\omega_{d}} P^{(d-1)/2}_k(1)
\nonumber\\
&= \frac{D(d, k)}{\omega_{d}}
= \frac{2 k^{d-1}}{(d-1)! \omega_d}+{\rm o}(k^{d-1}),
\end{align*}
where we have used the fact that
$P^{\lambda}_k(1)
={_2}F_1 \left( -k, k+2 \lambda, \lambda+1/2; 0 \right)
= 1, \lambda> - 1/2$ \cite{NIST10}.

Next we consider the DPP on $\S^{d}$ for
fixed $d \in \N$ and $L \in \N$ such that
the correlation kernel is given 
by the following finite sum \cite{BMOC16},
\begin{align}
K^{(N(d, L))}_{{\rm harmonic}(\S^d)}(x, x')
&: = \sum_{k=0}^{L-1} K^{\cY_{(d, k)}}(x, x')
= \frac{1}{\omega_{d}} 
\sum_{k=0}^{L-1} D(d, k) P^{(d-1)/2}_k(x \cdot x')
\nonumber\\
&= \frac{1}{\omega_{d}} 
\sum_{k=0}^{L-1} \frac{d-1+2k}{d-1} C^{(d-1)/2}_k(x \cdot x'),
\label{eqn:K_sphere2}
\end{align}
where the total number of points on $\S^d$
is given by
\begin{align}
N(d, L) &= \sum_{k=0}^{L-1} D(d, k)
= \frac{2L+d-2}{d} \binom{d+L-2}{L-1}
\nonumber\\
&=\frac{2}{d !} L^d +{\rm o}(L^d).
\label{eqn:N_d_ell}
\end{align}
The DPP $(\Xi, K^{(N)}_{{\rm harmonic}(\S^d)}, d\sigma_d(x))$
is rotationally invariant in $\S^d$ and 
is called the {\it harmonic ensemble}
in $\S^d$ with $N$ points by Beltr\'an {\it et al.} \cite{BMOC16}.
We note the recurrence relation of the Gegenbauer polynomials
(see, Eq.(18.9.7) in \cite{NIST10}),
\[
(n+\lambda) C^{\lambda}_n(x)
=\lambda( C^{\lambda+1}_n(x)-C^{\lambda+1}_{n-2}(x)).
\]
This implies that
\[
\frac{d-1+2k}{d-1}
C^{(d-1)/2}_k(x)
=C^{(d+1)/2}_k(x)-C^{(d+1)/2}_{k-2}(x), \quad k \geq 2.
\]
Since $C^{\lambda}_0(x)=1, C^{\lambda}_1(x)=2 \lambda x$, 
we obtain the following expression for the
correlation kernel,
\[
K^{(N(d, L))}_{{\rm harmonic}(\S^d)}(x, x')
= \frac{1}{\omega_{d}} \left[ C^{(d+1)/2}_{L-1}(x \cdot x')
+ C^{(d+1)/2}_{L-2}(x \cdot x') \right].
\]
If we introduce the {\it Jacobi polynomials}
defined as \cite{NIST10} 
\[
P_n^{(\alpha, \beta)}(x)
:= \frac{(\alpha+1)_n}{n!} 
{_2}F_1 \left(-n, n+\alpha+\beta+1; \alpha+1; \frac{1-x}{2} \right),
\]
and use the contiguous relation, 
$(b-a) \, {_2}F_1(a,b ; c; z)
+a \, {_2}F_1(a+1, b; c; z)-b \, {_2}F_1(a, b+1; c; z)=0$,
the above is written as follows \cite{BMOC16},
\begin{equation}
K^{(N(d, L))}_{{\rm harmonic}(\S^d)}(x, x')
=\frac{1}{\omega_d} \frac{N(d, L)}{\binom{L+d/2-1}{L-1}}
P^{(d/2, (d-2)/2)}_{L-1}(x \cdot x'),
\label{eqn:K_sphere4}
\end{equation}
where $\binom{L+d/2-1}{L-1}:=\Gamma(L+d/2)/\{(L-1)! \Gamma(d/2+1)\} 
=P^{(d/2, (d-2)/2)}_{L-1}(1)$.

In particular, when $d=1$, for
$x=(x^{(1)}, x^{(2)})=(\sin \theta, \cos \theta)$,
$x^{\prime}=({x^{\prime}}^{(1)}, {x^{\prime}}^{(2)})
=(\sin \theta^{\prime}, \cos \theta^{\prime}) \in \S^1 \subset \R^2$,
$\theta, \theta' \in [0, 2 \pi)$, we have
$x \cdot x^{\prime}=\cos(\theta-\theta^{\prime})$
and
\begin{align}
K^{(N(1, L))}_{{\rm harmonic}(\S^1)}(x, x^{\prime}) d \sigma_1(x)
&= \frac{1}{2 \pi} {_2}F_1 \left(
\frac{1-(2L-1)}{2}, \frac{1+(2L-1)}{2}; \frac{3}{2};
\sin^2 \frac{\theta-\theta^{\prime}}{2} \right) d \theta
\nonumber\\
&= \frac{\sin \{(2L-1)(\theta-\theta^{\prime})/2 \}}
{\sin\{(\theta-\theta^{\prime})/2\}}
\frac{d \theta}{2 \pi}
\nonumber\\
&=\frac{\sin\{N(\theta-\theta^{\prime})/2\}}{\sin\{(\theta-\theta^{\prime})/2\}}
\frac{d \theta}{2 \pi},
\label{eqn:S1_CUE}
\end{align}
where we have used the fact that
$N(1, L)=2L-1$ given by (\ref{eqn:N_d_ell}).
This verifies the identification of the 1-sphere case
of the present DPP with the CUE,
$(\Xi, K^{\AN}, \lambda_{[0, 2 \pi)}(d \theta))$,
given in Section \ref{sec:Lie_group}.

On the other hand, when $d=2$,
(\ref{eqn:N_d_ell}) gives
$N(2, L)=L^2$ and
\begin{align*}
K^{(N(2, L))}_{{\rm harmonic}(\S^2)}(x, x^{\prime})
&= \frac{L^2}{4 \pi} 
{_2}F_1 \left( -L+1, L+1; 2; \frac{1-x \cdot x^{\prime}}{2} \right)
\nonumber\\
&= \frac{N}{4 \pi}
{_2}F_1 \left(
-\sqrt{N}+1, \sqrt{N}+1; 2; \frac{\|x-x^{\prime}\|^2_{\R^3}}{4} \right),
\end{align*}
which is different from $K_{\S^2}^{(N)}(x, x^{\prime})$
given by (\ref{eqn:K_S2}) in Section \ref{sec:sphere_finite}.

\subsection{Euclidean family of infinite DPPs on $\R^{d}$}
\label{sec:EuclideanDPP}

We consider the vicinity of the north pole $e_{d+1}=(0, \dots, 0, 1)$ 
on $\S^{d}$
and put $\theta_{d} = R/L$, $R \in [0, \infty)$.
Then, as $L \to \infty$, 
the polar coordinates (\ref{eqn:coordinate1}) behave as
\begin{align*}
x^{(1)} &\sim
 \frac{R}{L} \sin \theta_{d-1} \cdots \sin \theta_2 \sin \theta_1
=: \frac{1}{L} \widetilde{x}^{(1)},
\nonumber\\
x^{(a)} &\sim
\frac{R}{L} \sin \theta_{d-1} \cdots \sin \theta_a \cos \theta_{a-1}
=: \frac{1}{L} \widetilde{x}^{(a)},
\quad a=2, \dots, d,
\nonumber\\
x^{(d+1)} &\sim
1-\frac{1}{2} \left(\frac{R}{L}\right)^2.
\end{align*}
In this case, for $x, x' \in \S^{d}$, 
\[
x \cdot x'
= \sum_{a=1}^{d+1} x^{(a)} {x'}^{(a)}
= 1-\frac{1}{2 L^2} \|\widetilde{x}-\widetilde{x}'\|_{\R^d}^2 + {\rm o} \left( \frac{1}{L^2} \right),
\quad \mbox{as $L \to \infty$}, 
\]
where $\widetilde{x}, \widetilde{x}' \in \R^{d}$ and
$\|\cdot\|_{\R^d}$ denotes the Euclidean norm in $\R^{d}$.
Hence we can conclude that
\begin{equation}
x \cdot x' = \cos \left( \frac{r}{L} \right)
+{\rm o} \left( \frac{1}{L^2} \right),
\quad \mbox{with $r := \|\widetilde{x}-\widetilde{x}'\|_{\R^d}$},
\quad \mbox{ as $L \to \infty$}. 
\label{eqn:settingA}
\end{equation}
In this limit, the measure on $\S^{d}$
given by (\ref{eqn:sigma}) behaves as
\begin{align*}
d \sigma_{d}(x)
&\sim \left( \frac{R}{L} \right)^{d-1} 
\sin^{d-2} \theta_{d-1} \cdots \sin \theta_2 \,
d \theta_1 \cdots d \theta_{d-1} \frac{dR}{L}
\nonumber\\
&= \frac{1}{L^d} 
d \sigma_{d-1}(\widehat{x}) R^{d-1} dR
=\frac{1}{L^{d}} d \widetilde{x},
\quad \widehat{x} \in \S^{d-1}, \quad \widetilde{x} \in \R^{d}.
\end{align*}

The following limit is proved for the correlation kernel
$K^{(N(d, L))}_{{\rm harmonic}(\S^d)}$ given by (\ref{eqn:K_sphere2}). 
\begin{lem}
\label{thm:asymK}
When {\rm (\ref{eqn:settingA})} holds, 
the limit
\[
k^{(d)}(r) =
\lim_{L \to \infty} \frac{1}{L^{d}} K^{(N(d, L))}_{{\rm harmonic}(\S^d)}(x, x')
\]
exists and has the following expressions,
\begin{align}
k^{(d)}(r) 
&= 
\frac{J_{d/2}(r)}{(2 \pi r)^{d/2}}
\label{eqn:lim_K2}
\\
&= \frac{1}{(2 \pi)^{d/2} r^{(d-2)/2}}
\int_0^1 s^{d/2} J_{(d-2)/2}(r s) d s,
\label{eqn:lim_K}
\end{align}
where $J_{\nu}(z)$ is
the Bessel function of the first kind 
with index $\nu$ defined by 
{\rm (\ref{eqn:Bessel_function})}. 
\end{lem}
\noindent{\it Proof} \,
For $d \in \N$ and $r \in (0, \infty)$, 
the following formula of Mehler--Heine type 
is known (Theorem 8.1.1 in \cite{Sze75});
for $\alpha, \beta \in \R$, 
\[
\lim_{n \to \infty} n^{-\alpha}
P^{(\alpha, \beta)}_{n} \left( \cos \frac{r}{n} \right)
= \left( \frac{r}{2} \right)^{-\alpha}
J_{\alpha}(r),
\]
where the limit is uniform on compact subset of $\C$.
Then under (\ref{eqn:settingA}), (\ref{eqn:K_sphere4}) gives
\[
\lim_{L \to \infty} \frac{1}{L^{d}} K^{(N(d, L))}_{{\rm harmonic}(\S^d)}
(x, x')
= \frac{2}{\omega_d d!} \Gamma(d/2+1)
\left( \frac{2}{r} \right)^{d/2}
J_{d/2}(r).
\]
By (\ref{eqn:omega}) and the equality 
$\Gamma(2z)=\{2^{2z}/(2 \sqrt{\pi})\}
\Gamma(z) \Gamma \left( z+1/2 \right)$, 
we can confirm that
$2 \Gamma(d/2+1)/(\omega_d d! )
= 1/\{(2 \sqrt{\pi})^{d} \}$.
Hence (\ref{eqn:lim_K2}) is proved.
By the integral formula (see, for instance, Eq.(10.22.1) in \cite{NIST10}),
\[
\int z^{\nu+1} J_{\nu}(z) dz = z^{\nu+1} J_{\nu+1}(z),
\]
we can derive (\ref{eqn:lim_K}) from (\ref{eqn:lim_K2}).
The proof is complete. \qed
\vskip 0.3cm

This result implies that for each $d \in \N$
we obtain an infinite-dimensional DPP on $\R^{d}$
such that it is uniform and isotropic on $\R^{d}$ and
the correlation kernel is given by
\begin{equation}
K^{(d)}(x, x') = k^{(d)}(\|x-x'\|_{\R^d}),
\quad x, x' \in \R^{d},
\label{eqn:bK1}
\end{equation}
where $k^{(d)}(r)$ is given by (\ref{eqn:lim_K2}) and (\ref{eqn:lim_K}).

We can give the following alternative expression
for $K^{(d)}$. 
\begin{lem}
For $d \in \N$, the correlation kernel
$K^{(d)}$ given by {\rm (\ref{eqn:bK1})}
with {\rm (\ref{eqn:lim_K})}
is written as
\begin{equation}
K^{(d)}(x, x')
= \frac{1}{(2\pi)^{d}}
\int_{\R^{d}} 
{\bf 1}_{\B^d}(y) e^{i (x-x') \cdot y} d y
= \frac{1}{(2\pi)^{d}}
\int_{\B^{d}} e^{i (x-x') \cdot y} d y,
\label{eqn:Fourier}
\end{equation}
where $\B^d$ denotes the unit ball centered at the origin; 
$\B^d := \{ y \in \R^d : |y| \leq 1\}$.
\end{lem}
\noindent{\it Proof} \,
The statement is proved for $d=1$ and 2 by direct calculation as follows.
For $d=1$, (\ref{eqn:lim_K}) gives
\[
k^{(1)}(r)=\sqrt{\frac{r}{2 \pi}}
\int_0^1 s^{1/2} J_{-1/2}(rs) ds.
\]
Here we use the equality
$J_{-1/2}(z)=\sqrt{2/(\pi z)} \cos z$.
Then
\[
k^{(1)}(r) =\frac{1}{\pi} \int_0^1 \cos (rs) ds
=\frac{1}{2 \pi} \int_{-1}^{1} e^{i r y} d y,
\]
which gives (\ref{eqn:Fourier}) with $d=1$,
if we put $r=x-x'$ and regard an interval $[-1, 1] \subset \R$
as $\B^1$. 
For $d=2$, (\ref{eqn:lim_K}) gives
\begin{equation}
k^{(2)}(r) = \frac{1}{2\pi} \int_0^1 s J_0(rs) ds.
\label{eqn:K3_1}
\end{equation}
We use the following integral representation for
$J_0$ given as Eq.(10.9.1) in \cite{NIST10},
\[
J_0(z)= \frac{1}{\pi} \int_0^{\pi} \cos(z \cos \varphi) d \varphi
=\frac{1}{2 \pi} \int_0^{2 \pi} e^{iz \cos \varphi} d \varphi.
\]
Hence (\ref{eqn:K3_1}) is written as
\begin{equation}
k^{(2)}(r)= \frac{1}{(2 \pi)^2}
\int_0^1 ds \, s \int_0^{2 \pi} d \theta \,
e^{i r s \cos \theta}.
\label{eqn:K3_2}
\end{equation}
We can identify the integral variables $(s, \theta)$ in
(\ref{eqn:K3_2}) with the polar coordinates in $\R^2$
and (\ref{eqn:Fourier}) with $d=2$ is obtained,
if we recognize $r=\|x-x'\|_{\R^2}$,
$s=\|y\|_{\R^2}$, and $(x-x') \cdot y=r s \cos \theta$.
Now we assume $d \geq 3$.
In this case RHS of (\ref{eqn:Fourier}) is given by
\begin{align*}
I &:= \frac{1}{(2\pi)^{d}}
\int_0^1 d s \, s^{d-1}
\int_0^{2 \pi} d \theta_1 \int_0^{\pi} d \theta_2 \, \sin \theta_2
\nonumber\\
& \qquad \times
\cdots \times \int_0^{\pi} d \theta_{d-2} \sin^{d-3} \theta_{d-2}
\int_0^{\pi} d \theta_{d-1} \sin^{d-2} \theta_{d-1}
e^{i r s \cos \theta_{d-1}}.
\end{align*}
Since
\[
\int_0^{2 \pi} d \theta_1 \int_0^{\pi} d \theta_2 \, \sin \theta_2
\cdots \int_0^{\pi} d \theta_{d-2} \sin^{d-3} \theta_{d-2}
=\sigma_{d-2}(\S^{d-2})
 = \omega_{d-2},
\]
we have
\[
I=\frac{w_{d-2}}{(2\pi)^{d}}
\int_0^1 d s \, s^{d-1}
\int_0^{\pi} d \theta_{d-1} \, \sin^{d-2} \theta_{d-1}
e^{i r s \cos \theta_{d-1}}.
\]
If we use the following integral representation
of the Bessel function of the first kind, 
\[
J_m(z)= \frac{1}{\sqrt{\pi} \Gamma(m+1/2)}
\left( \frac{z}{2} \right)^m
\int_0^{\pi} \sin^{2m}(\theta) e^{iz \cos \theta} d \theta,
\quad m \in \frac{1}{2} \N,
\]
which is obtained from 
Eq.(10.9.4) in \cite{NIST10}, 
the equivalence between (\ref{eqn:Fourier}) and
(\ref{eqn:lim_K}) is verified. 
Hence the proof is complete. \qed
\vskip 0.3cm

The kernel (\ref{eqn:Fourier}) is obtained as
the correlation kernel $K_{S_1}$ given by
(\ref{eqn:K_simple2}) in Corollary \ref{thm:main3},
if we consider the case such that
$S_1=S_2=\R^d$,
$\lambda_1(d x)=d x$, 
$\lambda_2(d y)=\nu(d y)=d y$,
$\psi_1(x, y)=e^{i x \cdot y}/(2\pi)^{d/2}$, and
$\Gamma=\B^d \subsetneq \R^d$.
We see $\Psi_1(x)^2 := \|\psi_1(x, \cdot)\|_{L^2(\Gamma, d \nu)}^2
\equiv |\B^d|/(2 \pi)^d$, $x \in \R^d$, 
where the volume of $\B^{d}$ is denoted by
$|\B^d| =\pi^{d/2}/\Gamma((d+2)/2)$.

The kernels $K^{(d)}$ on $\R^d, d \in \N$ have been
studied by Zelditch and others 
(see \cite{Zel00,SZ02,Zel09,CH15} and references therein),
who regarded them as the Szeg\"o kernels
for the reduced Euclidean motion group
\cite{Sug90,YY07}.
Here we call the DPPs associated with
the correlation kernels in this form 
the {\it Euclidean family of DPPs} on $\R^d, d \in \N$.

\begin{df}
\label{thm:EuclideanDPP}
The Euclidean family of DPPs is 
a one-parameter {\rm ($d \in \N$)} family of \\
$\Big( \Xi, K^{(d)}_{\rm Euclid}, d x \Big)$
with 
\begin{align*}
K^{(d)}_{\rm Euclid}(x, x')
& := 
\frac{1}{(2 \pi)^{d/2}} \frac{J_{d/2}(\|x-x'\|_{\R^d})}{\|x-x'\|_{\R^d}^{d/2}}
\nonumber\\
&=\frac{1}{(2 \pi)^{d/2}}
\frac{1}{\|x-x'\|_{\R^d}^{(d-2)/2}}
\int_0^1 s^{d/2} J_{(d-2)/2} (\|x-x'\|_{\R^d} s) ds
\nonumber\\
&= \frac{1}{(2 \pi)^d}
\int_{\R^d} {\bf 1}_{\B^d}(y)
e^{i (x-x') \cdot y} d y
= \frac{1}{(2 \pi)^d}
\int_{\B^d} 
e^{i (x-x') \cdot y} d y,
\quad x, x' \in \R^d.
\end{align*}
\end{df}

The above result is summarized as follows \cite{KS1}.
\begin{prop}
\label{thm:limit_Euclid}
The following is established for $d \in \N$,
\[
\left( \frac{d!}{2} \right)^{1/d}
N^{1/d} \circ \Big( \Xi, K^{(N)}_{{\rm harmonic}(\S^d)}, d \sigma_{d}(x) \Big)
\weak \Big( \Xi, K^{(d)}_{\rm Euclid}, d x \Big).
\]
\end{prop}

We see that
\[
K^{(d)}_{\rm Euclid}(x, x)
=\lim_{r \to 0} \frac{1}{(2\pi)^{d/2}} \frac{J_{d/2}(r)}{r^{d/2}}
=\frac{1}{2^d \pi^{d/2} \Gamma((d+2)/2)}.
\]
Then the Euclidean family of DPPs
are uniform on $\R^d$
with densities
\[
\rho^{(d)}_{\rm Euclid}= \frac{1}{2^d \pi^{d/2} \Gamma((d+2)/2)}
\]
with respect to the Lebesgue measures $d x$ of $\R^d$. 

For lower dimensions, the correlation kernels 
and the densities are given as follows,
\begin{align*}
K^{(1)}_{\rm Euclid}(x, x')
&= \frac{\sin(x-x')}{\pi(x-x')}
=K_{\rm sinc}(x, x')
\quad
\mbox{with} \quad \rho^{(1)}_{\rm Euclid} = \frac{1}{\pi},
\nonumber\\
K^{(2)}_{\rm Euclid}(x, x')
&= \frac{J_1(\|x-x'\|_{\R^2})}{2 \pi \|x-x'\|_{\R^2}}
\quad
\mbox{with} \quad \rho^{(2)}_{\rm Euclid} = \frac{1}{4 \pi},
\nonumber\\
K^{(3)}_{\rm Euclid}(x, x')
&= \frac{1}{2 \pi^2 \|x-x'\|_{\R^3}^2}
\left( \frac{\sin \|x-x'\|_{\R^3}}{\|x-x'\|_{\R^3}}
-\cos \|x-x'\|_{\R^3} \right)
\quad
\mbox{with} \quad \rho^{(3)}_{\rm Euclid} = \frac{1}{6 \pi^2}.
\end{align*}
This family of DPPs includes the DPP with the sinc kernel 
$K_{\rm sinc}$ as the
lowest dimensional case with $d=1$.
Since $(\Xi, K^{(N)}_{{\rm harmonic}(\S^1)}, d \sigma_1)$
has been identified with the CUE, $(\Xi, K^{\AN}, \lambda_{[0, 2 \pi)})$,  
by (\ref{eqn:S1_CUE}), Proposition \ref{thm:limit_Euclid} 
can be regarded as the multidimensional extension
of the limit theorem from the CUE to 
the DPP with $K_{\rm sinc}$ given by the first line of
(\ref{eqn:Lie_bulk}). 
Note that, 
if $d$ is odd,
\[
k^{(d)}(r)
=\left(- \frac{1}{2\pi r} \frac{d}{dr} \right)^{(d-1)/2}
\frac{\sin r}{\pi r}.
\]
This is proved by Rayleigh's formula
for the spherical Bessel function of the first kind
(Eq. (10.49.14) in \cite{NIST10});
\[
j_m(x) := \sqrt{\frac{\pi}{2 x}}
J_{m+1/2}(x)
=x^m \left( - \frac{1}{x} 
\frac{d}{dx} \right)^m \frac{\sin x}{x},
\quad m \in \N.
\]

\SSC
{Concluding Remarks} \label{sec:remarks}

In Section \ref{sec:DPP_S_d}, we studied the finite DPPs
$(\Xi, K^{(N)}_{{\rm harmonic}(\S^d)}, d \sigma_d(x))$, $N \in \N$
called the harmonic ensembles on $\S^d$, 
$d \in \N$ \cite{BMOC16}.
Then we proved as Proposition \ref{thm:limit_Euclid} 
in Section \ref{sec:EuclideanDPP}
that their bulk scaling limits are 
given by $(\Xi, K_{\rm Euclid}^{(d)}, dx)$, 
which we call the Euclidean family of DPPs.
On $\S^2$, 
there are two distinct types of uniform and isotropic DPPs,
one of which is the harmonic ensemble
$(\Xi, K^{(N)}_{{\rm harmonic}(\S^2)}, d \sigma_2(x))$
studied in Section \ref{sec:DPP_S_d} \cite{BMOC16} ,
and other of which is the DPP called
the spherical ensemble $(\Xi, K^{(N)}_{\S^2}, d \sigma_2(x))$
studied in Section \ref{sec:sphere_finite} \cite{Kri09,AZ15}.
As mentioned above, the scaling limit of
the former is given by $(\Xi, K^{(2)}_{\rm Euclid}, dx)$,
while as given by Proposition \ref{thm:S2_Ginibre}
the bulk scaling limit of the latter is
$(\Xi, K_{\rm Ginibre}^A, \lambda_{\rN(0,1; \C)}(dx))$,
which is equivalent with
$(\Xi, K^{(1)}_{\rm Heisenberg}, \lambda_{\rN(0, 1; \C)}(dx))$.
The spherical ensemble on $\S^2$ shall be generalized 
to DPPs on the higher dimensional spheres 
$\S^{2d} \simeq \C^d$, $d \geq 2$ so that
they are uniform and isotropic and
their bulk scaling limits are given by
DPPs in the Heisenberg family.
The papers \cite{BE18,BE19} will be useful.

With $L^2(S, \lambda)$ and $L^2(\Gamma, \nu)$, 
we can consider the system of 
{\it biorthonormal functions},
which consists of a pair of distinct families of measurable functions 
$\{\psi(x, \gamma) : x \in S, \gamma \in \Gamma \}$
and $\{\varphi(x, \gamma) : x \in S, \gamma \in \Gamma \}$
satisfying the biorthonormality relations
\begin{equation}
\langle \psi(\cdot, \gamma), \varphi(\cdot, \gamma') \rangle_{L^2(S, \lambda)}
\nu(d \gamma)=\delta(\gamma-\gamma') d \gamma,
\quad \gamma, \gamma' \in \Gamma.
\label{eqn:bi_ortho}
\end{equation}
If the integral kernel defined by
\begin{equation}
K^{\rm bi}(x, x')
=\int_{\Gamma} \psi(x, \gamma) \overline{\varphi(x', \gamma)} \nu(d \gamma),
\quad x, x' \in S, 
\label{eqn:K_bi}
\end{equation}
is of finite rank, we can construct a finite DPP on $S$
whose correlation kernel is given by (\ref{eqn:K_bi})
following a standard method of random matrix theory
(see, for instance, Appendix C in \cite{Kat19a}).
By the biorthonormality (\ref{eqn:bi_ortho}),
it is easy to verify that 
$K^{\rm bi}$ is a projection kernel,
but it is not necessarily an orthogonal projection.
This observation means that such a DPP
is not constructed by the method reported in this paper.
Generalization of the present framework in order to cover
such DPPs associated with biorthonormal systems
is required.
Moreover, the dynamical extensions of DPPs
called {\it determinantal processes}
(see, for instance, \cite{Kat15_Springer})
shall be studied in the context of the present paper.

For finite DPPs, we can readily derive the
systems of {\it stochastic interacting particle systems}
whose stationary states are given by the DPPs.
For example, with $N \in \N$, 
the system of stochastic differential equations (SDEs)
on $\S^1$ \cite{HW96,Kat14},
\begin{equation}
dX_j(t)=dB_j(t)+\frac{1}{2} \sum_{\substack{1 \leq k \leq N, \cr k \not=j}}
\cot \frac{X_j(t)-X_k(t)}{2} dt,
\quad j=1, \dots, N, \quad t \geq 0,
\label{eqn:SDE1}
\end{equation}
driven by independent one-dimensional standard Brownian motions
$B_j(t), j=1, \dots, N, t \geq 0$ has the DPP
$(\Xi, K^{\AN}, \lambda_{[0, 2 \pi)}(dx))$
given in Section \ref{sec:Lie_group} as a stationary probability measure.
Another example is given by
the system of SDEs on $\C$,
\begin{equation}
d Z_j(t) = dB^{\C}_j(t)
- \frac{(N+1) Z_j(t)}{1+|Z_j(t)|^2} dt
+ \sum_{\substack{1 \leq k \leq N, \cr k \not=j}}
\frac{Z_j(t)-Z_k(t)}{|Z_j(t)-Z_k(t)|^2} dt, 
\quad j=1, \dots, N, \quad t \geq 0,
\label{eqn:SDE2}
\end{equation}
$N \in \N$, 
driven by independent complex Brownian motions
$B^{\C}_j(t) :=B^{\rR}_j(t)+i B^{\rI}_j(t)$,  
where $B^{\rR}_j(t), B^{\rI}_j(t)$ are independent
one-dimensional standard Brownian motions, 
$j=1, \dots, N, t \geq 0$, 
does the DPP 
$(\sum_j \delta_{Z_j}, K_{G_1^{-1} G_2}^{(N)}, \lambda(dz))$  
of Krishnapur \cite{Kri09} as a stationary probability measure, 
which is obtained as the stereographic projection
of $(\Xi, K_{\S^2}^{(N)}, d \sigma_2(x))$
as explained in Remark 5 in Section \ref{sec:sphere_finite}.
A general theory has been developed by Osada {\it et al.} for
{\it infinite-dimensional stochastic differential equations} (ISDEs), 
some of which have infinite DPPs as invariant probability measures
\cite{Osa96,Osa12,Osa13a,OT14a,KOT17,OT20}.
We expect to obtain the {\it universal} ISDEs along the limit theorems
given in Propositions \ref{thm:limit_Euclid} and \ref{thm:S2_Ginibre} 
taking account of the fact that (\ref{eqn:SDE1}) and (\ref{eqn:SDE2})
might give useful approximations to characterize 
Osada's Dyson/Ginibre ISDEs.

\vskip 0.5cm
\noindent{\bf Acknowledgements} \,
One of the present authors (MK)  thanks Christian Krattenthaler very much
for his hospitality in Fakult\"{a}t f\"{u}r Mathematik, Universit\"{a}t Wien,
where this study was started
on his sabbatical leave from Chuo University. 
The authors express their gratitude to Michael Schlosser, 
Masatake Hirao, and Shinji Koshida for valuable discussion. 
This work was supported by
the Grant-in-Aid for Scientific Research
(C) (No.26400405), (No.19K03674), 
(B) (No.18H01124), 
(S) (No.16H06338), and
(A) (No.21H04432)
of Japan Society for the Promotion of Science.

\appendix
\SSC{Weyl Denominator Formulas}
\label{sec:Weyl_denominators}

The {\it Weyl denominator formulas} for classical root systems
play a fundamental role in Lie theory and related area.
For reduced root systems they are given in the form,
\[
\sum_{w \in W} \det(w) e^{w(\rho)-\rho}
=\prod_{\alpha \in R_{+}}(1-e^{-\alpha}),
\]
where $W$ is the Weyl group, $R_+$ the set of positive roots
and $\rho=\frac{1}{2} \sum_{\alpha \in R_+} \alpha$.

For classical root systems $\AN, \BN, \CN$ and $\DN$, 
$N \in \N$, 
the explicit forms are given as follows,
\begin{align}
\mbox{ (type $\AN$)}
\quad & \det_{1 \leq j, k \leq N} \Big(z_k^{j-1} \Big)
= \prod_{1 \leq j < k \leq N} (z_k-z_j),
\nonumber\\
\mbox{ (type $\BN$)}
\quad &
\det_{1 \leq, j, k \leq N} 
\Big( z_k^{j-N}-z_k^{N+1-j} \Big)
= \prod_{\ell=1}^{N} z_{\ell}^{1-N} (1-z_{\ell}) 
\prod_{1 \leq j < k \leq N} (z_k-z_j)(1-z_j z_k),
\nonumber\\
\mbox{ (type $\CN$)}
\quad &
\det_{1 \leq, j, k \leq N} 
\Big( z_k^{j-N-1}-z_k^{N+1-j} \Big)
= \prod_{\ell=1}^{N} z_{\ell}^{-N} (1-z_{\ell}^2) 
\prod_{1 \leq j < k \leq N} (z_k-z_j)(1-z_j z_k),
\nonumber\\
\mbox{ (type $\DN$)}
\quad &
\det_{1 \leq, j, k \leq N} 
\Big( z_k^{j-N}+z_k^{N-j} \Big)
= 2 \prod_{\ell=1}^{N} z_{\ell}^{1-N}
\prod_{1 \leq j < k \leq N} (z_k-z_j)(1-z_j z_k), 
\label{eqn:Weyl_denominator}
\end{align}
respectively. See, for instance, \cite{RS06}.

If we change the variables as
\begin{equation}
z_k=e^{-2 i \zeta_k}, \quad \zeta_k \in \C, \quad k=1, \dots, N,
\label{eqn:x_to_z}
\end{equation}
then, the following equalities are derived from
the above.

\begin{lem}
\label{thm:Weyl_trigonometric}
For $\zeta_k \in \C, k=1, \dots,N$, 
the following equalities are established.
\begin{align*}
{\rm (type} \, \AN {\rm )}
\quad &
\det_{1 \leq j, k \leq N} 
\Big[ e^{- i (\cN^{\AN}-2J^{\AN}(j)) \zeta_k} \Big]
= (2i) ^{N(N-1)/2}
\prod_{1 \leq j < k \leq N} 
\sin(\zeta_k-\zeta_j).
\nonumber\\
{\rm (type} \, \BN {\rm )}
\quad &
\det_{1 \leq j, k \leq N} 
\Big[ \sin\{(\cN^{\BN}-2 J^{\BN}(j)) \zeta_k\} \Big]
\nonumber\\
& \qquad
= 2^{N(N-1)}
\prod_{\ell=1}^N \sin \zeta_{\ell}
\prod_{1 \leq j < k \leq N} 
\sin(\zeta_k-\zeta_j) \sin(\zeta_k+\zeta_j),
\nonumber\\
{\rm (type} \, \CN {\rm )}
\quad &
\det_{1 \leq j, k \leq N} 
\Big[ \sin \{ (\cN^{\CN}-2 J^{\CN}(j)) \zeta_k\} \Big]
\nonumber\\
& \qquad
= 2^{N(N-1)}
\prod_{\ell=1}^N \sin(2 z_{\ell}) 
\prod_{1 \leq j < k \leq N} 
\sin (\zeta_k-\zeta_j) \sin (\zeta_k+\zeta_j),
\nonumber\\
{\rm (type} \, \DN {\rm )}
\quad &
\det_{1 \leq j, k \leq N} 
\Big[ \cos \{(\cN^{\DN}-2 J^{\DN}(j)) \zeta_k \} \Big]
\nonumber\\
& \qquad
= 2^{(N-1)^2}
\prod_{1 \leq j < k \leq N} 
\sin(\zeta_k-\zeta_j) \sin (\zeta_k+\zeta_j),
\end{align*}
where $\cN^{\RN}$ and $J^{\RN}(j)$, 
$\RN=\AN, \BN, \CN, \DN$, are given by
{\rm (\ref{eqn:N_R_classic})}
and {\rm (\ref{eqn:J_R_classic})}.
\end{lem}

\SSC{Jacobi Theta Functions}
\label{sec:Jacobi}

Let
\[
z=e^{v \pi i}, \quad q=e^{\tau \pi i},
\]
for $v \in \C$ and $\tau \in \H$. 
The Jacobi theta functions are defined as follows \cite{WW27,NIST10}, 
\begin{align}
\vartheta_0(v; \tau) &= 
\sum_{n \in \Z} (-1)^n q^{n^2} z^{2n}
= 1 + 2 \sum_{n=1}^{\infty} (-1)^n e^{\tau \pi i n^2} \cos (2n \pi v), 
\nonumber\\
\vartheta_1(v; \tau) &= i \sum_{n \in \Z} (-1)^n q^{(n-1/2)^2} z^{2n-1}
= 2 \sum_{n=1}^{\infty} (-1)^{n-1} e^{\tau \pi i (n-1/2)^2} \sin\{ (2n-1) \pi v \},
\nonumber\\
\vartheta_2(v; \tau)
&= \sum_{n \in \Z} q^{(n-1/2)^2} z^{2n-1}
= 2 \sum_{n=1}^{\infty} e^{\tau \pi i (n-1/2)^2} \cos \{ (2n-1) \pi v\},
\nonumber\\
\vartheta_3(v; \tau) 
&=\sum_{n \in \Z} q^{n^2} z^{2n}
= 1 + 2 \sum_{n=1}^{\infty} e^{\tau \pi i n^2} \cos (2n \pi v). 
\label{eqn:theta}
\end{align}
(Note that the present functions 
$\vartheta_{\mu}(v; \tau), \mu=1,2,3$ are denoted by
$\vartheta_{\mu}(\pi v,q)$,
and $\vartheta_0(v;\tau)$ by $\vartheta_4(\pi v,q)$ in \cite{WW27}.)
For $\Im \tau >0$, $\vartheta_{\mu}(v; \tau)$, $\mu=0,1,2,3$
are holomorphic for $|v| < \infty$.
The parity with respect to $v$ is given by
\begin{equation}
\vartheta_1(-v; \tau)=-\vartheta_1(v; \tau),
\quad
\vartheta_{\mu}(-v; \tau)=\vartheta_{\mu}(v; \tau), \quad
\mu=0,2,3,
\label{eqn:even_odd}
\end{equation}
and they have the quasi-double-periodicity; 
\begin{align}
\vartheta_{\mu}(v+1; \tau) &=
\begin{cases}
\vartheta_{\mu}(v; \tau), & \mu=0, 3,
\cr
- \vartheta_{\mu}(v; \tau), & \mu=1, 2,
\end{cases}
\label{eqn:quasi_periodic1}
\\
\vartheta_{\mu}(v+\tau; \tau) &=
\begin{cases}
-e^{-(2v+\tau) \pi i } \vartheta_{\mu}(v; \tau),
& \mu=0, 1,
\cr
e^{-(2v+\tau) \pi i } \vartheta_{\mu}(v; \tau),
& \mu=2,3.
\end{cases}
\label{eqn:quasi_periodic2}
\end{align}
The following relations are derived by (\ref{eqn:theta}),
\begin{align}
\vartheta_0 \left( v+\frac{\tau}{2}; \tau \right)
&= i e^{-(v+\tau/4) \pi i} \vartheta_1(v; \tau),
\nonumber\\
\vartheta_1 \left( v+\frac{\tau}{2}; \tau \right)
&= i e^{-(v+\tau/4) \pi i} \vartheta_0(v; \tau),
\nonumber\\
\vartheta_2 \left( v+\frac{\tau}{2}; \tau \right)
&= e^{-(v+\tau/4) \pi i} \vartheta_3(v; \tau),
\nonumber\\
\vartheta_3 \left( v+\frac{\tau}{2}; \tau \right)
&= e^{-(v+\tau/4) \pi i} \vartheta_2(v; \tau).
\label{eqn:half_tau}
\end{align}
By the definition (\ref{eqn:theta}), 
when $\tau \in \H$, 
\begin{align}
& \vartheta_1(0; \tau)=\vartheta_1(1; \tau)=0, \qquad
\vartheta_1(x; \tau) > 0, \quad x \in (0,1),
\nonumber\\
& \vartheta_2(-1/2; \tau)=\vartheta_2(1/2; \tau)=0, \qquad
\vartheta_2(x; \tau) > 0, \quad x \in (-1/2, 1/2),
\nonumber\\
& \vartheta_0(x; \tau) > 0, \quad \vartheta_3(x; \tau) > 0, \quad x \in \R.
\label{eqn:theta_values}
\end{align}
The asymptotics
\begin{align}
& \vartheta_0(v; \tau) \sim 1, \quad
\vartheta_1(v; \tau) \sim 2 e^{\tau \pi i/4} \sin (\pi v), \quad
\vartheta_2(v; \tau) \sim 2 e^{\tau \pi i/4} \cos(\pi v), \quad
\vartheta_3(v; \tau) \sim 1,
\nonumber\\
&\qquad \qquad \qquad \mbox{in} \quad
\Im \tau \to + \infty \quad
(i.e., \quad q=e^{\tau \pi i} \to 0)
\label{eqn:theta_asym}
\end{align}
are known.
We will use the following functional equations known as
{\it Jacobi's imaginary transformation} \cite{WW27,NIST10},
\begin{align}
\vartheta_0(v; \tau)
&= e^{\pi i/4} \tau^{-1/2} e^{-\pi i v^2/\tau}
\vartheta_2 \left( \frac{v}{\tau}; - \frac{1}{\tau} \right),
\nonumber\\
\vartheta_1(v; \tau)
&= e^{3 \pi i/4} \tau^{-1/2} e^{-\pi i v^2/\tau}
\vartheta_1 \left( \frac{v}{\tau}; - \frac{1}{\tau} \right),
\nonumber\\
\vartheta_2(v; \tau)
&= e^{\pi i/4} \tau^{-1/2} e^{-\pi i v^2/\tau}
\vartheta_0 \left( \frac{v}{\tau}; - \frac{1}{\tau} \right), 
\nonumber\\
\vartheta_3(v; \tau)
&= e^{\pi i/4} \tau^{-1/2} e^{-\pi i v^2/\tau}
\vartheta_3 \left( \frac{v}{\tau}; - \frac{1}{\tau} \right). 
\label{eqn:Jacobi_imaginary}
\end{align}

\SSC{Macdonald Denominators}
\label{sec:Macdonald_denominators}

Assume that $N \in \N$.
As extensions of the Weyl denominators for
classical root systems, Rosengren and Schlosser \cite{RS06}
studied the {\it Macdonald denominators} for the seven
types of irreducible reduced affine root systems \cite{Mac72},
$W^{\RN}(\z; \tau)$, $\z=(z_1, \dots, z_N) \in \C^N$, $\tau \in \H$, 
$\RN=\AN$, $\BN, \BNv, \CN, \CNv, \BCN, \DN$,
$N \in \N$. (See also \cite{War02,Kra05}.)
Up to trivial factors they are written using the Jacobi theta functions
as follows.
\begin{align}
W^{\AN}(\z; \tau) &=
\prod_{1 \leq j < k \leq N} \vartheta_1(z_k-z_j; \tau),
\nonumber\\
W^{\BN}(\z; \tau) &=
\prod_{\ell=1}^N \vartheta_1(z_{\ell}; \tau)
\prod_{1 \leq j < k \leq N} \Big\{
\vartheta_1(z_k-z_j; \tau) \vartheta_1(z_k+z_j; \tau) \Big\},
\nonumber\\
W^{\BNv}(\z; \tau) &=
\prod_{\ell=1}^N \vartheta_1(2 z_{\ell}; 2 \tau)
\prod_{1 \leq j < k \leq N} \Big\{
\vartheta_1(z_k-z_j; \tau) \vartheta_1(z_k+z_j; \tau) \Big\},
\nonumber\\
W^{\CN}(\z; \tau) &=
\prod_{\ell=1}^N \vartheta_1(2 z_{\ell}; \tau)
\prod_{1 \leq j < k \leq N} \Big\{
\vartheta_1(z_k-z_j; \tau) \vartheta_1(z_k+z_j; \tau) \Big\},
\nonumber\\
W^{\CNv}(\z; \tau) &=
\prod_{\ell=1}^N \vartheta_1 \left(z_{\ell}; \frac{\tau}{2} \right)
\prod_{1 \leq j < k \leq N} \Big\{
\vartheta_1(z_k-z_j; \tau) \vartheta_1(z_k+z_j; \tau) \Big\},
\nonumber\\
W^{\BCN}(\z; \tau) &=
\prod_{\ell=1}^N \Big\{ \vartheta_1(z_{\ell}; \tau) 
\vartheta_0(2 z_{\ell}; 2 \tau) 
\Big\}
\prod_{1 \leq j < k \leq N} \Big\{
\vartheta_1(z_k-z_j; \tau) \vartheta_1(z_k+z_j; \tau) \Big\},
\nonumber\\
W^{\DN}(\z; \tau) &=
\prod_{1 \leq j < k \leq N} \Big\{
\vartheta_1(z_k-z_j; \tau) \vartheta_1(z_k+z_j; \tau) \Big\}.
\label{eqn:Macdonald_denominators}
\end{align}


\end{document}